\documentclass[9pt,letterpaper]{extarticle}

\usepackage[
letterpaper,
textwidth=6.5in,
textheight=9in,
centering
]{geometry}

\usepackage[UTF8]{ctex}
\usepackage{float}
\usepackage{graphicx}
\usepackage{mathrsfs}
\usepackage{amsmath}
\usepackage{amssymb}
\usepackage{amsfonts}
\usepackage{bm}
\usepackage{setspace}
\usepackage{abstract}
\usepackage{caption}
\usepackage{diagbox}
\usepackage{multirow}
\usepackage[numbers,sort&compress]{natbib}
\usepackage{color}
\usepackage{makecell}
\usepackage{threeparttable}
\usepackage{booktabs}
\usepackage{subcaption}
\usepackage{placeins}
\usepackage{array}
\usepackage{enumitem}
\usepackage{titlesec}
\usepackage{pdflscape}
\usepackage[hidelinks]{hyperref}
\usepackage[title]{appendix}
\usepackage{dashbox}
\usepackage{tikz}
\usepackage{etoolbox}
\usepackage{changepage}
\usepackage{siunitx}

\usepackage{indentfirst}

\graphicspath{{./}{figures/}{pics/}}

\titleformat{\section}
{\large\bfseries}
{\thesection}
{0.75em}
{}

\titleformat{\subsection}
{\normalsize\bfseries}
{\thesubsection}
{0.75em}
{}

\titleformat{\subsubsection}
{\normalsize\bfseries}
{\thesubsubsection}
{0.75em}
{}

\titleformat{\paragraph}[hang]
{\normalsize\bfseries}
{\theparagraph}
{1em}
{}

\titlespacing{\section}
{0pt}
{2.0ex plus 0.4ex minus 0.2ex}
{0.8ex}

\titlespacing{\subsection}
{0pt}
{1.6ex plus 0.3ex minus 0.2ex}
{0.6ex}

\titlespacing{\subsubsection}
{0pt}
{1.4ex plus 0.3ex minus 0.2ex}
{0.5ex}

\titlespacing{\paragraph}
{0pt}
{1.2ex plus 0.2ex minus 0.2ex}
{0.6ex plus 0.2ex}

\newcommand{\figplaceholder}[2][]{%
	\IfFileExists{#2}{%
		\includegraphics[#1]{#2}%
	}{%
		\fbox{%
			\parbox[c][0.20\textheight][c]{0.88\linewidth}{%
				\centering
				图像占位：\\[0.5em]
				\texttt{\detokenize{#2}}%
			}%
		}%
	}%
}

\newcommand{\figmaybe}[2][]{%
	\IfFileExists{#2}{%
		\includegraphics[#1]{#2}%
	}{%
		\fbox{%
			\parbox[c][0.20\textheight][c]{0.88\linewidth}{%
				\centering
				图像占位：\\[0.5em]
				\texttt{\detokenize{#2}}%
			}%
		}%
	}%
}

\AtBeginEnvironment{tabular}{\small}
\AtBeginEnvironment{tabular*}{\small}

\AtBeginEnvironment{tablenotes}{\footnotesize}

\begin{document}

	\title{\LARGE Two-Step MV-DeepONet: Probabilistic Operator Learning for Uncertainty Propagation Driven by Random Input Fields}
	
	\author{
		\begin{tabular}{@{}c@{}}
			\small Yupei Nie$^{a}$ \quad
			Lei Wang$^{a,*}$ \quad
			Jiasen Liu$^{b}$\\[0.35em]
			\small\emph{$^{a}$School of Mathematics and Physics, North China Electric Power University, Beijing 102206, P.R. China} \\
			\small\emph{$^{b}$School of Energy Power and Mechanical Engineering, North China Electric Power University, Beijing 102206, P.R. China}
		\end{tabular}
	}
	
	\date{}
	
	\maketitle
	\vspace{-1.5em}
	
	\footnotetext[1]{Corresponding author: \href{mailto:50901924@ncepu.edu.cn}{50901924@ncepu.edu.cn}}
	
	\captionsetup[figure]{
		labelfont={bf},
		labelformat={default},
		labelsep=period,
		name={\small Fig.}
	}
	
	\captionsetup[table]{
		labelfont={bf},
		labelformat={default},
		labelsep=period,
		name={\small Table}
	} 
	
	\noindent\textbf{Abstract:}\par
	\begingroup
	\small
	\setlength{\parindent}{2em}
	\setlength{\parskip}{0.25em}
	
	Forward uncertainty propagation in complex physical systems can induce structured covariance across field-valued outputs. For a probabilistic surrogate, the total predictive covariance comprises the covariance of conditional means across input realizations and the average conditional predictive covariance. Probabilistic DeepONet (Prob-DeepONet) provides lightweight uncertainty quantification by predicting pointwise Gaussian means and variances in a single forward pass, but its conditional predictive covariance is restricted to a diagonal form. To represent cross-location conditional dependence without explicitly parameterizing a full high-dimensional covariance matrix, we develop a two-step mean–variance DeepONet (two-step MV-DeepONet) through two principal modifications. First, two-step training is used to decouple output-basis learning from the input-to-coefficient mapping, together with basis orthogonalization and subspace rotation. Second, Gaussian probabilistic modeling is transferred from the high-dimensional physical output space to the low-dimensional rotated coefficient space. Mapping these probabilistic coefficients through the shared basis induces a generally non-diagonal conditional predictive covariance in the physical output space while retaining single-pass inference. A Frobenius-norm error decomposition and corresponding upper bound identify low-rank covariance compressibility, trunk-subspace approximation, finite-sample statistical error, and coefficient-space covariance estimation as the principal factors governing covariance recovery. Numerical experiments on three representative problems governed by partial differential equations (PDEs) and a hypersonic blunt-body aerothermal problem show improved generalization, more structured uncertainty bands, and accurate recovery of off-diagonal correlation patterns compared with Prob-DeepONet. Conformal calibration further improves empirical interval coverage, while two-step MV-DeepONet maintains sharper prediction intervals.
	\endgroup
	
	\vspace{0.5em}
	\noindent\textbf{Keywords:} Input uncertainty propagation; Two-step MV-DeepONet; Structured covariance recovery; Probabilistic operator learning; Uncertainty quantification
	\vspace{0.3em}
	
	
	\newpage

\section{Introduction}


Uncertainty quantification is essential for credible modeling and reliable design of complex engineering systems. A standard workflow generally includes uncertainty identification and characterization, uncertainty propagation, and uncertainty analysis~\cite{UP_Mohammadi2022}. Among these steps, input uncertainty propagation characterizes the variability of model outputs induced by random inputs and is therefore central to engineering design under uncertainty. In complex physical systems, uncertain inputs may arise from material properties, boundary conditions, external loads, and other sources of variability, and can be represented as random variables, stochastic processes, or random fields~\cite{RFEM_Liu1986}. For a probabilistic surrogate under random inputs, the total predictive covariance consists of the covariance of the predicted means across input realizations and the average conditional predictive covariance over inputs. The former captures cross-location dependence induced by random inputs~\cite{StochasticHeat_XiuKarniadakis2003}, while the latter may also contain structured dependence that is not explicitly represented under a diagonal covariance assumption~\cite{Simpson2022}.

Classical methods for input uncertainty propagation can be broadly classified into simulation methods, local approximation methods, most probable point methods, functional expansion methods, and numerical integration methods~\cite{UP_LeeChen2009}. Among these approaches, polynomial chaos expansion and generalized polynomial chaos are widely used functional expansion methods for uncertainty propagation~\cite{UP_Review_Abdi2024}. These methods can provide efficient approximations for smooth problems with moderate stochastic dimensions. However, their computational cost may increase rapidly with the number of uncertain inputs and the approximation order~\cite{SFEM_GhanemSpanos1991,gPC_XiuKarniadakis2002,SC_XiuHesthaven2005}. With the development of data-driven modeling, machine-learning approaches have increasingly been adopted for uncertainty propagation, including analytical approximations through trained neural networks~\cite{SDP_Petersen2024,PNN_ShekhovtsovFlach2018}, interval-based propagation~\cite{IntervalDL_Betancourt2022}, and offline sampling combined with deterministic surrogates~\cite{RSM_Sofi2020,GPDimReduction_Tripathy2016,ANNFatigue_Giannella2023}. In addition to input-induced uncertainty, surrogate predictions may also be affected by limited or noisy training data, model inadequacy, and variability in the training procedure~\cite{Gawlikowski2023}. Representative approaches for quantifying such predictive uncertainty include Bayesian neural networks~\cite{BNN_Neal1996}, Monte Carlo dropout~\cite{MCDropout_Gal2016}, and deep ensembles~\cite{DeepEnsemble_Lakshminarayanan2017}. Other approaches directly parameterize predictive distributions, including evidential deep learning~\cite{EvidentialDL_Sensoy2018} and mean--variance estimation (MVE)~\cite{MVE_NixWeigend1994}. Among these, MVE provides a lightweight formulation that assumes an input-dependent Gaussian target distribution and jointly predicts its conditional mean and variance by minimizing the Gaussian negative log-likelihood (NLL)~\cite{MVE_NixWeigend1994}. This formulation enables input-dependent predictive variance to be estimated in a single forward pass. However, joint mean--variance optimization may be unstable~\cite{PNN_Seitzer2022}; a mean-only warm-up with fixed variance can therefore be used to improve optimization stability~\cite{MVE_Sluijterman2024}.

Both classical input-uncertainty propagation methods and conventional machine-learning approaches to predictive uncertainty are commonly formulated in finite-dimensional spaces. For physical systems driven by random input fields, however, both the inputs and the corresponding responses are naturally function-valued, motivating operator-learning methods that directly approximate mappings between function spaces~\cite{Subedi2026OperatorLearning,Kovachki2023NeuralOperator}. Among representative neural operator architectures, Deep Operator Network (DeepONet) provides a flexible framework grounded in the universal approximation theory of nonlinear operators~\cite{ChenChen1995,Lu_DeepONet_2021}. DeepONet represents an operator through a branch network and a trunk network. The branch network maps sensor observations of the input function to input-dependent coefficients, whereas the trunk network maps output coordinates to output-domain basis functions. Their inner product reconstructs the operator output, thereby separating input-function encoding from output-space representation. DeepONet has demonstrated strong capability in a variety of scientific and engineering problems involving high-dimensional field prediction~\cite{PIDeepONet_Goswami2023,HeatTransferDeepONet_Sahin2024,CrackDeepONet_Goswami2022,DeepMMNet_Electroconvection_Cai2021,DeepMMNet_Hypersonics_Mao2021}. Nevertheless, standard deterministic DeepONet provides only point predictions and does not explicitly quantify predictive uncertainty~\cite{NeuralOperatorComparison_Lu2022}, motivating extensions of DeepONet that incorporate uncertainty quantification for more reliable predictions.

Uncertainty quantification for DeepONet has been explored through Bayesian methods, randomized-prior ensembles, and calibration-based approaches. Bayesian methods provide a principled probabilistic framework but often require iterative inference and remain sensitive to prior specification and posterior approximation~\cite{BDeepONet_Lin2023,VBDeepONet_Garg2023,AlphaVIDeepONet_Lone2026}. Randomized-prior ensembles provide an effective means of quantifying epistemic uncertainty, while their computational cost grows with ensemble size and their performance depends on careful tuning~\cite{UQDeepONet_Yang2022}. Calibration-based methods can improve empirical interval coverage, whereas they primarily operate as post-processing procedures and do not directly modify the underlying uncertainty representation~\cite{ConformalizedDeepONet_Moya2025,ConformalNeuralOperator_Kobayashi2025}. Recent extensions further broaden
the scope to conditional field generation, stochastic operators, and bounded uncertainty~\cite{DeepSetOperatorUQ_Ma2026,DiffusionStochasticOperator_Huynh2026,IntervalPropagation_Faza2026}. Among lightweight alternatives, Probabilistic DeepONet (Prob-DeepONet) directly parameterizes the predictive distribution using a formulation analogous to MVE, jointly predicting the output mean and variance through a Gaussian NLL~\cite{MVE_NixWeigend1994,DeepONetGridUQ_Moya2023}. Representative studies have employed Prob-DeepONet for post-fault trajectory prediction, uncertainty-guided sample selection, and active learning~\cite{DeepONetGridUQ_Moya2023,ActiveOperatorLearning_Winovich2026}. 
Unlike Bayesian posterior sampling and model ensembles, Prob-DeepONet estimates predictive means and variances in a single forward pass, providing an efficient approach to uncertainty quantification in operator learning~\cite{DeepONetGridUQ_Moya2023}. However, its commonly adopted pointwise Gaussian formulation models marginal means and variances independently at individual output locations, resulting in a diagonal conditional covariance structure.  Extending Prob-DeepONet to represent cross-location conditional dependence while retaining its lightweight formulation would therefore provide a more expressive probabilistic operator-learning framework for engineering applications.

Accordingly, this study considers probabilistic operator learning under random input fields and focuses on recovering the total predictive covariance of high-dimensional output fields. The principal methodological target is the conditional predictive component of this covariance, for which the pointwise conditional-independence assumption in Prob-DeepONet limits the representation of cross-location dependence. To address this limitation, we develop a two-step MV-DeepONet framework through two principal modifications. First, the two-step training strategy of Lee and Shin~\cite{TwoStepDeepONet_LeeShin2024} is incorporated to decouple the learning of output-space basis functions from the input-to-coefficient mapping. Second, probabilistic modeling is transferred from the high-dimensional physical output space to the low-dimensional modal coefficient space. Although a diagonal Gaussian model is retained in the coefficient space, the probabilistic modal coefficients are jointly mapped to multiple output locations through the shared basis functions, allowing a generally non-diagonal conditional predictive covariance to be represented in the physical output space. The resulting formulation avoids explicit learning or storage of the full high-dimensional covariance matrix while retaining the lightweight structure and single-pass inference of Prob-DeepONet, thereby providing a structured representation of the conditional component of the total predictive covariance.
The decoupled training strategy follows the two-step DeepONet formulation of Lee and Shin~\cite{TwoStepDeepONet_LeeShin2024}, in which the output-space basis is first learned and orthogonalized, followed by the learning of the input-to-coefficient mapping. This separation facilitates the construction of a linearly independent and numerically stable output basis and provides a low-dimensional representation for the subsequent coefficient-space probabilistic regression. 
The two-step formulation has subsequently been applied to a range of physical and engineering problems, including brittle fracture, three-dimensional electromagnetic field prediction, stochastic mechanical metamaterials, poroelasticity, Riemann problems, and hypersonic aerothermal analysis~\cite{CrackTwoStepDeepONet_Kiyani2025,ComplexDeepONet_Jiang2026,MetamaterialNeuralOperator_Jin2025,PoroelasticityDeepONet_Park2025,RiemannONets_Peyvan2024,WaveriderDeepONet_Shukla2024}. In this study, this decoupling provides a methodological foundation for probabilistic modeling in the coefficient space and structured output covariance recovery.

This study proposes a two-step MV-DeepONet framework for forward uncertainty propagation driven by random input fields and focuses on structured predictive uncertainty characterized by the total predictive covariance of high-dimensional output fields. The main contributions are summarized as follows:

\noindent\textbf{(1)}
A two-step MV-DeepONet framework is developed for forward uncertainty propagation driven by random input fields. By decoupling output-space basis learning from the input-to-coefficient mapping and performing probabilistic regression in the low-dimensional coefficient space, the proposed method relaxes the pointwise conditional-independence assumption of Prob-DeepONet. Through the shared basis functions, each probabilistic modal coefficient simultaneously affects multiple output locations, allowing a generally non-diagonal conditional predictive covariance to be represented in the physical output space. As a supplementary extension, deep ensembles and post-hoc calibration are used to account for epistemic uncertainty and improve empirical interval coverage.

\noindent\textbf{(2)}
A theoretical analysis of output covariance recovery is established through an error decomposition and a corresponding upper bound. The analysis identifies the low-rank compressibility of the target covariance, the approximation quality of the learned trunk subspace, finite-sample statistical error, and coefficient-space covariance estimation as the principal factors governing recovery accuracy. The resulting low-rank formulation avoids explicit learning or storage of the full high-dimensional covariance matrix while preserving lightweight, single-pass inference.

\noindent\textbf{(3)}
The proposed framework is validated on the reaction--diffusion, Burgers, and Darcy equations, together with a hypersonic aerothermal prediction problem. Compared with Prob-DeepONet, the proposed method exhibits improved generalization to unseen inputs, produces tighter and more structured uncertainty bands, and recovers correlation patterns consistent with the underlying physical mechanisms. These results highlight the complementary roles of the two methodological components: the two-step strategy helps reduce the mutual compensation between basis functions and coefficients, while coefficient-space probabilistic modeling propagates uncertainty jointly to multiple output locations through the shared basis functions. 

The remainder of this paper is organized as follows. Section~\ref{sec:method} reviews the relevant DeepONet formulations and develops the two-step MV-DeepONet framework. Section~\ref{sec:results} evaluates the proposed method in terms of mean prediction, predictive uncertainty, and output covariance recovery using several representative examples. Section~\ref{sec:conclusions} summarizes the paper and discusses directions for future work.

\section{Method}
\label{sec:method}

This section first reviews the standard DeepONet architecture and its two-step training strategy, followed by the pointwise probabilistic formulation of Prob-DeepONet. Building on these foundations, we develop a two-step MV-DeepONet that performs Gaussian modeling in a rotated low-dimensional coefficient space and represents structured conditional predictive covariance in the physical output space through a shared basis representation. Finally, a Frobenius-norm error decomposition and the corresponding upper bound are derived to characterize the principal factors governing total predictive covariance recovery.

\subsection{Deep Operator Network}

DeepONet~\cite{Lu_DeepONet_2021} is a neural operator architecture that learns a mapping between infinite-dimensional function spaces, i.e., an operator $\mathcal{G}:\mathcal{U}\to\mathcal{S}$ that maps $u\in\mathcal{U}$ to $s=\mathcal{G}(u)\in\mathcal{S}$, where $u$ and $s$ denote the input and output functions, respectively. It approximates $\mathcal{G}$ through two subnetworks. Given the sampled input
\begin{equation}
	\mathbf{u}
	=
	[u(x_1),u(x_2),\ldots,u(x_m)]^{\top}
	\in\mathbb{R}^{m},
\end{equation}
the branch network produces an input-dependent coefficient vector
\begin{equation}
	\mathbf{b}_{\theta}(\mathbf{u})
	=
	[b_1(\mathbf{u}),\ldots,b_p(\mathbf{u})]^{\top}
	\in\mathbb{R}^{p}.
\end{equation}

For the discrete output coordinates
\(\{y_i\}_{i=1}^{M}\), the trunk network generates the basis matrix
\begin{equation}
	\Phi_{\omega}
	=
	[\phi_j(y_i)]_{i=1,\ldots,M;\,j=1,\ldots,p}
	\in\mathbb{R}^{M\times p}.
\end{equation}
The discretized output prediction is then expressed as
\begin{equation}
	\hat{\mathbf{s}}(\mathbf{u})
	=
	\Phi_{\omega}\mathbf{b}_{\theta}(\mathbf{u})
	\in\mathbb{R}^{M}.
	\label{eq:deeponet_vector_prediction}
\end{equation}

Thus, the trunk network learns basis functions over the output domain, while the branch network predicts their input-dependent coefficients; $\theta$ and $\omega$ denote their respective trainable parameters.

For \(K\) training samples, the standard DeepONet is trained
end-to-end by minimizing
\begin{equation}
	\min_{\theta,\omega}
	\frac{1}{K}
	\sum_{k=1}^{K}
	\left\|
	\Phi_{\omega}
	\mathbf{b}_{\theta}(\mathbf{u}^{(k)})
	-
	\mathbf{s}^{(k)}
	\right\|_{2}^{2}.
	\label{eq:standard_deeponet_training}
\end{equation}

\subsection{Two-Step Training of DeepONet}
\label{subsec:twostep_deeponet_base}

For high-dimensional output fields, end-to-end DeepONet training requires the joint optimization of the branch and trunk networks in a large coupled problem. The two-step strategy decomposes this problem into two sequential subproblems, decoupling output-basis learning from input-dependent coefficient learning~\cite{TwoStepDeepONet_LeeShin2024}, thereby simplifying optimization and providing a stable coordinate system for subsequent probabilistic modeling in the coefficient space.

\vspace{4pt}
\noindent\textbf{Step 1: Trunk Basis Learning.}
The first step learns the trunk network alone. To isolate basis learning from the nonlinear branch mapping, the branch output is temporarily replaced by a trainable coefficient matrix $\mathbf{A} \in \mathbb{R}^{p \times K}$, and the trunk is trained by
\begin{equation}
	(\omega^{*},\mathbf{A}^{*})
	=
	\arg\min_{\omega,\mathbf{A}}
	\left\|
	\Phi_{\omega}\mathbf{A}-\mathbf{S}
	\right\|_{F}^{2},
	\label{eq:twostep_trunk_training}
\end{equation}
where $\Phi^{*}:=\Phi_{\omega^{*}}$ and $\mathbf{S}=[\mathbf{s}^{(1)},\ldots,\mathbf{s}^{(K)}] \in\mathbb{R}^{M\times K}$ collects the $K$ output snapshots columnwise. This yields the low-rank approximation $\widehat{\mathbf{S}}^{*}=\Phi^{*}\mathbf{A}^{*}\approx\mathbf{S}$, in which the columns of $\Phi^{*}$ span the learned low-dimensional subspace of the output field.

\vspace{4pt}
\noindent\textbf{Intermediate Step: Basis Orthogonalization.}
The optimized trunk basis matrix is orthogonalized through the QR
factorization
\begin{equation}
	\Phi^{*}
	=
	\mathbf{Q}^{*}\mathbf{R}^{*},
	\qquad
	(\mathbf{Q}^{*})^{\top}\mathbf{Q}^{*}
	=
	\mathbf{I}_{p},
	\label{eq:qr_decomposition}
\end{equation}
where \(\mathbf{Q}^{*}\in\mathbb{R}^{M\times p}\) is column-orthonormal and $\mathbf{R}^{*}\in\mathbb{R}^{p\times p}$ is upper triangular. The
corresponding coefficient targets are defined as
\begin{equation}
	\mathbf{C}^{*}
	=
	\mathbf{R}^{*}\mathbf{A}^{*}
	=
	[\mathbf{c}_1^{*},\ldots,\mathbf{c}_K^{*}]
	\in\mathbb{R}^{p\times K}.
	\label{eq:coeff_target_matrix}
\end{equation}
Accordingly, the first-step reconstruction can be written as
\(\Phi^{*}\mathbf{A}^{*}=\mathbf{Q}^{*}\mathbf{C}^{*}\).

\vspace{4pt}
\noindent\textbf{Step 2: Branch Coefficient Regression.}
With $\mathbf{Q}^{*}$ fixed, the branch network is trained to predict the coefficient vector associated with each input function. Denoting its prediction by
$\hat{\mathbf{c}}_{\theta}(\mathbf{u})\in\mathbb{R}^{p}$, the training objective is
\begin{equation}
	\theta^{*} = \arg\min_{\theta} \frac{1}{K} \sum_{k=1}^{K} \left\| \hat{\mathbf{c}}_{\theta}(\mathbf{u}^{(k)})-\mathbf{c}_k^{*} \right\|_{2}^{2}.
	\label{eq:twostep_branch_training}
\end{equation}

The resulting output prediction is reconstructed as
\begin{equation}
	\hat{\mathbf{s}}(\mathbf{u})
	=
	\mathbf{Q}^{*}\hat{\mathbf{c}}_{\theta^{*}}(\mathbf{u}).
	\label{eq:twostep_output_reconstruction}
\end{equation}
The original operator-learning problem is thereby reduced to a finite-dimensional supervised regression problem in coefficient space.

\subsection{Probabilistic Deep Operator Network}
\label{sec:single_step_mvdeeponet}

To further quantify predictive uncertainty, Prob-DeepONet~\cite{DeepONetGridUQ_Moya2023} extends the deterministic DeepONet by imposing pointwise Gaussian distributions at each physical output location. For a discrete input representation $\mathbf{u}$, the Prob-DeepONet prediction is modeled as
\begin{equation}
	\hat{\mathbf{s}}(\mathbf{u}) \sim \mathcal{N}\left( \hat{\boldsymbol{\mu}}_{\eta}(\mathbf{u}), \operatorname{diag}\left( \hat{\boldsymbol{\sigma}}_{\eta}^{2}(\mathbf{u}) \right) \right),
	\label{eq:single_joint_gaussian}
\end{equation}
where
$\eta=\{\theta,\omega\}$ denotes all trainable parameters, with
$\theta=\{\theta_{\mu},\theta_{\sigma}\}$ and
$\omega=\{\omega_{\mu},\omega_{\sigma}\}$.
Prob-DeepONet adopts separate branch--trunk paths for the mean and
log-variance of the output field:
\begin{equation}
	\hat{\boldsymbol{\mu}}_{\eta}(\mathbf{u})
	=
	\Phi_{\mu,\omega_{\mu}}
	\mathbf{b}_{\mu,\theta_{\mu}}(\mathbf{u})
	\in\mathbb{R}^{M},
	\qquad
	\hat{\boldsymbol{\ell}}_{\eta}(\mathbf{u})
	=
	\Phi_{\sigma,\omega_{\sigma}}
	\mathbf{b}_{\sigma,\theta_{\sigma}}(\mathbf{u})
	\in\mathbb{R}^{M},
	\qquad
	\hat{\boldsymbol{\sigma}}_{\eta}^{2}(\mathbf{u})
	=
	\exp\left(
	\hat{\boldsymbol{\ell}}_{\eta}(\mathbf{u})
	\right)
	\in\mathbb{R}^{M}.
	\label{eq:single_mean_var_output}
\end{equation}
Here, the exponential is applied elementwise, ensuring strictly positive predictive variances. For $K$ training samples, the model is trained by minimizing the pointwise Gaussian NLL:
\begin{equation}
	\mathcal{L}_{\mathrm{single}} = \frac{1}{2K} \sum_{k=1}^{K} \sum_{i=1}^{M} \left[ \frac{\left( s_{i}^{(k)} - \hat{\mu}_{\eta,i}(\mathbf{u}^{(k)}) \right)^2}{\hat{\sigma}_{\eta,i}^{2}(\mathbf{u}^{(k)})} + \log \hat{\sigma}_{\eta,i}^{2}(\mathbf{u}^{(k)}) \right].
	\label{eq:single_nll}
\end{equation}
This formulation learns an input-dependent marginal distribution at each output location and provides pointwise prediction intervals for assessing marginal
predictive reliability~\cite{DeepONetGridUQ_Moya2023}.

\subsection{Two-Step MV-DeepONet for Covariance Recovery}
\label{sec:twostep_mvdeeponet}

To model the total predictive covariance of a probabilistic operator surrogate under random input functions, the proposed framework combines two-step basis learning~\cite{TwoStepDeepONet_LeeShin2024} with Gaussian probabilistic modeling in a low-dimensional coefficient space. The proposed extension focuses on relaxing the pointwise conditional-independence assumption in Prob-DeepONet while retaining a computationally tractable covariance representation.

\subsubsection{Structured Output Covariance under Random Inputs}
\label{sec:covariance_motivation}

Let $\mathbf{U}$ denote a random input function and let
$\hat{\mathbf{s}}$ denote the field-valued prediction of a
probabilistic operator surrogate. By the law of total covariance, the predictive covariance can be decomposed as
\begin{equation}
	\operatorname{Cov}(\hat{\mathbf{s}})
	=
	\operatorname{Cov}_{\mathbf{U}}
	\left(
	\mathbb{E}[\hat{\mathbf{s}}\mid\mathbf{U}]
	\right)
	+
	\mathbb{E}_{\mathbf{U}}
	\left[
	\operatorname{Cov}
	(\hat{\mathbf{s}}\mid\mathbf{U})
	\right].
	\label{eq:total_predictive_covariance}
\end{equation}
The first term characterizes the variation of the predicted mean across
random input realizations, and the second term represents the
average conditional predictive covariance of the surrogate.

To illustrate why the first term in
Eq.~\eqref{eq:total_predictive_covariance} may exhibit cross-location
dependence, consider a random perturbation \(\delta u\) applied to the input
function \(u\). The corresponding output field is $s = \mathcal{G}(u+\delta u)$. A first-order expansion around \(u\) gives
\begin{equation}
	s
	\approx
	\mathcal{G}(u)
	+
	\nabla\mathcal{G}(u)\,\delta u,
	\label{eq:taylor_operator}
\end{equation}
where \(\nabla\mathcal{G}(u)\) denotes the linearized sensitivity of the
operator. The resulting output covariance can therefore be approximated
as
\begin{equation}
	\operatorname{Cov}(s)
	\approx
	\nabla\mathcal{G}(u)
	\operatorname{Cov}(\delta u)
	\nabla\mathcal{G}(u)^{\top}.
	\label{eq:operator_cov_taylor}
\end{equation}
Equation~\eqref{eq:operator_cov_taylor} shows that the covariance of the input perturbation is transformed by the sensitivity of the physical operator. For operators defined by partial differential equations or complex engineering systems, the same input perturbation can affect multiple output locations through differential operators, boundary conditions, or global physical constraints. Consequently, the sensitivities at different output locations generally share common components of the input perturbation, resulting in a generally non-diagonal \(\operatorname{Cov}(s)\). 
Regarding the second term, the pointwise Gaussian assumption in Prob-DeepONet yields a diagonal covariance:
\begin{equation}
	\widehat{\boldsymbol{\Sigma}}_{\mathrm{Prob}}(\mathbf{u})
	=
	\operatorname{diag}
	\left(
	\hat{\sigma}_{\eta,1}^{2}(\mathbf{u}),
	\ldots,
	\hat{\sigma}_{\eta,M}^{2}(\mathbf{u})
	\right),
	\qquad
	\left[
	\widehat{\boldsymbol{\Sigma}}_{\mathrm{Prob}}(\mathbf{u})
	\right]_{ij}
	=
	0,
	\quad i\neq j.
	\label{eq:prob_diag_cov}
\end{equation}
which estimates the marginal uncertainty at each output location but does not represent covariance between distinct locations for a
fixed input. Consequently, the conditional predictive covariance, corresponding to the second term in Eq.~\eqref{eq:total_predictive_covariance}, is restricted to a diagonal form and therefore does not explicitly represent cross-location dependence within this component, motivating the two-step MV-DeepONet developed below.

\subsubsection{Probabilistic Extension Based on the Two-Step Strategy}

\vspace{4pt}
\noindent\textbf{Step 1: Trunk Basis Learning.}
The trunk training follows the original two-step DeepONet strategy~\cite{TwoStepDeepONet_LeeShin2024} without structural
modification. This step provides a stable low-dimensional output subspace for the subsequent probabilistic extension.

\vspace{4pt}
\noindent\textbf{Intermediate Step: Orthogonalization and Subspace
	Rotation.}
The branch network models the coefficients with a diagonal Gaussian covariance, which implicitly assumes uncorrelated coefficients. Since basis orthogonalization does not generally diagonalize the empirical coefficient covariance, an additional orthogonal rotation is introduced within the learned subspace.

After first-step trunk training and QR orthogonalization, the orthonormal basis $\mathbf{Q}^{*}$ yields the reconstruction
\begin{equation}
	\widehat{\mathbf{S}}^{*}
	=
	\Phi^{*}\mathbf{A}^{*}
	=
	\mathbf{Q}^{*}
	\left(
	\mathbf{R}^{*}\mathbf{A}^{*}
	\right)
	\approx
	\mathbf{S},
\end{equation}
where
\(\mathbf{Q}^{*}\in\mathbb{R}^{M\times p}\) is column-orthonormal and
\(\mathbf{R}^{*}\in\mathbb{R}^{p\times p}\) is upper triangular.

The training outputs are then centered as
\begin{equation}
	\bar{\mathbf{s}} = \frac{1}{K}\sum_{k=1}^{K}\mathbf{s}^{(k)}, \qquad \mathbf{x}^{(k)} = \mathbf{s}^{(k)} - \bar{\mathbf{s}}, \qquad \mathbf{X} = \left[ \mathbf{x}^{(1)},\ldots,\mathbf{x}^{(K)} \right] \in \mathbb{R}^{M\times K}.
\end{equation}
Their coefficients under the orthonormal basis \(\mathbf{Q}^{*}\) are
\begin{equation}
	\mathbf{C}_{Q}
	=
	\left(\mathbf{Q}^{*}\right)^{\top}\mathbf{X}
	=
	\left[
	\mathbf{c}_{Q}^{(1)},\ldots,
	\mathbf{c}_{Q}^{(K)}
	\right]
	\in
	\mathbb{R}^{p\times K},
	\label{eq:cq_star}
\end{equation}
with empirical covariance
\begin{equation}
	\widehat{\mathbf{B}}_{Q}
	=
	\frac{1}{K}
	\mathbf{C}_{Q}\mathbf{C}_{Q}^{\top}
	\in
	\mathbb{R}^{p\times p}.
	\label{eq:bq_star}
\end{equation}

Because $\widehat{\mathbf{B}}_{Q}$ is generally non-diagonal, an eigendecomposition is performed to rotate the coefficient coordinates and diagonalize their empirical covariance, thereby providing a representation better aligned with the subsequent diagonal Gaussian approximation:
\begin{equation}
	\widehat{\mathbf{B}}_{Q}
	=
	\widehat{\mathbf{V}}_{Q}
	\widehat{\boldsymbol{\Lambda}}_{Q}
	\widehat{\mathbf{V}}_{Q}^{\top},
	\qquad
	\widehat{\mathbf{V}}_{Q}^{\top}
	\widehat{\mathbf{V}}_{Q}
	=
	\mathbf{I}_{p},
	\label{eq:rotation_eig}
\end{equation}
where \(\widehat{\boldsymbol{\Lambda}}_{Q}
=
\operatorname{diag}
\left(
\widehat{\lambda}_{Q,1},\ldots,\widehat{\lambda}_{Q,p}
\right)\)
contains the eigenvalues of \(\widehat{\mathbf{B}}_{Q}\), and the rotated
basis is defined as
\(\widetilde{\mathbf{Q}}^{*}
=
\mathbf{Q}^{*}\widehat{\mathbf{V}}_{Q}\).
Since $\widehat{\mathbf{V}}_{Q}$ is orthogonal, $\widetilde{\mathbf{Q}}^{*}$ and $\mathbf{Q}^{*}$ span the same subspace and differ only in the orientation of the coordinate axes. The corresponding coefficient matrix is
\begin{equation}
	\mathbf{C}_{\widetilde Q}
	=
	\left(
	\widetilde{\mathbf{Q}}^{*}
	\right)^{\top}
	\mathbf{X}
	=
	\widehat{\mathbf{V}}_{Q}^{\top}
	\mathbf{C}_{Q},
	\label{eq:ctilde_star}
\end{equation}
and its empirical covariance satisfies
\begin{equation}
	\widehat{\mathbf{B}}_{\widetilde Q}
	=
	\frac{1}{K}
	\mathbf{C}_{\widetilde Q}
	\mathbf{C}_{\widetilde Q}^{\top}
	=
	\widehat{\mathbf{V}}_{Q}^{\top}
	\widehat{\mathbf{B}}_{Q}
	\widehat{\mathbf{V}}_{Q}
	=  
	\widehat{\boldsymbol{\Lambda}}_{Q}.
	\label{eq:bq_rotated}
\end{equation}

The rotation preserves the learned output subspace while diagonalizing the empirical coefficient covariance, thereby providing coefficient coordinates better suited to the subsequent diagonal Gaussian modeling.

\vspace{4pt}
\noindent\textbf{Step 2: Probabilistic Regression for Branch
	Coefficients.}
The second stage extends deterministic coefficient regression to
probabilistic modeling. For a given input \(\mathbf{u}\), the branch
network predicts the coefficient mean and log-variance:
\begin{equation}
	\hat{\boldsymbol{\mu}}_{c,\theta}(\mathbf{u})
	=
	\left[
	\hat{\mu}_{c,1}(\mathbf{u}),\ldots,
	\hat{\mu}_{c,p}(\mathbf{u})
	\right]^{\top}
	\in\mathbb{R}^{p},
	\qquad
	\hat{\boldsymbol{\ell}}_{c,\theta}(\mathbf{u})
	=
	\left[
	\hat{\ell}_{c,1}(\mathbf{u}),\ldots,
	\hat{\ell}_{c,p}(\mathbf{u})
	\right]^{\top}
	\in\mathbb{R}^{p}.
\end{equation}
The predicted coefficient variances are defined by
\begin{equation}
	\hat{\sigma}_{c,j}^{2}(\mathbf{u})
	=
	\exp
	\left(
	\hat{\ell}_{c,j}(\mathbf{u})
	\right),
	\qquad
	j=1,\ldots,p,
\end{equation}
which ensures their strict positivity. The corresponding coefficient
covariance matrix $\widehat{\boldsymbol{\Lambda}}_{c,\theta}(\mathbf{u})$
is diagonal, with its $j$th diagonal entry given by
$\hat{\sigma}_{c,j}^{2}(\mathbf{u})$. The modal coefficient vector is therefore modeled as
\begin{equation}
	\hat{\mathbf{c}}_{\theta}(\mathbf{u})
	\sim
	\mathcal{N}
	\left(
	\hat{\boldsymbol{\mu}}_{c,\theta}(\mathbf{u}),
	\widehat{\boldsymbol{\Lambda}}_{c,\theta}(\mathbf{u})
	\right).
	\label{eq:coef_distribution}
\end{equation}
The independence assumption in
Eq.~\eqref{eq:coef_distribution} is imposed on the modal coefficients
in the rotated coefficient space rather than on the physical output
variables.

For the \(k\)-th training sample, let
\(\mathbf{c}_{\widetilde Q}^{(k)}\) denote the empirical coefficient
target under \(\widetilde{\mathbf{Q}}^{*}\). The branch network is trained using the Gaussian NLL to predict the input-dependent coefficient means and variances:
\begin{equation}
	\mathcal{L}_{\mathrm{coef}}
	=
	\frac{1}{2K}
	\sum_{k=1}^{K}
	\sum_{j=1}^{p}
	\left[
	\frac{
		\left(
		c_{\widetilde Q,j}^{(k)}
		-
		\hat{\mu}_{c,j}
		\left(
		\mathbf{u}^{(k)}
		\right)
		\right)^2
	}{
		\hat{\sigma}_{c,j}^{2}
		\left(
		\mathbf{u}^{(k)}
		\right)
	}
	+
	\log
	\hat{\sigma}_{c,j}^{2}
	\left(
	\mathbf{u}^{(k)}
	\right)
	\right].
	\label{eq:coef_nll}
\end{equation}

\subsubsection{Output-Space Covariance Recovery}

For a given input \(\mathbf{u}\), the centered output field is reconstructed from the rotated orthonormal basis and the probabilistic modal coefficients as
\begin{equation}
	\hat{\mathbf{x}}_{\theta}(\mathbf{u})
	=
	\widetilde{\mathbf{Q}}^{*}
	\hat{\mathbf{c}}_{\theta}(\mathbf{u}).
\end{equation}

By the covariance propagation rule for linear transformations, the corresponding conditional predictive covariance in the physical output space is
\begin{equation}
	\widehat{\boldsymbol{\Sigma}}_{s}(\mathbf{u})
	=
	\operatorname{Cov}
	\left(
	\hat{\mathbf{x}}_{\theta}(\mathbf{u})
	\mid \mathbf{u}
	\right)
	=
	\widetilde{\mathbf{Q}}^{*}
	\widehat{\boldsymbol{\Lambda}}_{c,\theta}(\mathbf{u})
	\left(
	\widetilde{\mathbf{Q}}^{*}
	\right)^{\top}.
	\label{eq:single_input_covariance}
\end{equation}

Although
\(\widehat{\boldsymbol{\Lambda}}_{c,\theta}(\mathbf{u})\) is diagonal,
the output covariance in
Eq.~\eqref{eq:single_input_covariance} is generally non-diagonal.
Each probabilistic modal coefficient affects multiple output locations through the shared basis, thereby inducing structured output covariance without directly predicting a full covariance matrix. For two output locations \(y_i\) and \(y_j\), the corresponding
covariance is
\begin{equation}
	\left[
	\widehat{\boldsymbol{\Sigma}}_{s}(\mathbf{u})
	\right]_{ij}
	=
	\sum_{m=1}^{p}
	\widetilde{q}_m^{*}(y_i)
	\hat{\sigma}_{c,m}^{2}(\mathbf{u})
	\widetilde{q}_m^{*}(y_j),
	\label{eq:single_input_point_covariance}
\end{equation}
where \(\widetilde{q}_m^{*}(y_i)\) is the value of the \(m\)-th rotated
basis function at \(y_i\), and
\(\hat{\sigma}_{c,m}^{2}(\mathbf{u})\) is the predicted variance of the corresponding modal coefficient. For a given pair of distinct output locations, this sum may vanish under special configurations, including: (i) all basis functions vanish simultaneously at $y_i$ or $y_j$, i.e., $\widetilde{q}_m^{*}(y_\ell)=0$ for all $m=1,\ldots,p$ and some $\ell\in\{i,j\}$, a strong coincidence condition uncommon for neural-network bases on general output fields; (ii) the signed modal contributions cancel exactly, $\sum_{m=1}^{p}\widetilde{q}_m^{*}(y_i)\hat{\sigma}_{c,m}^{2}(\mathbf{u})\widetilde{q}_m^{*}(y_j)=0$, which requires the branch-predicted variances and the basis values to satisfy a strict algebraic relation that is generally unstable across inputs and locations; or (iii) $\hat{\sigma}_{c,m}^{2}(\mathbf{u})=0$ for all $m$, corresponding to a deterministic limiting case. This exact degeneracy is excluded by the adopted log-variance parameterization. 
Excluding the aforementioned cases, there exists at least one mode $m$ with $\widetilde{q}_m^{*}(y_i)\widetilde{q}_m^{*}(y_j)\hat{\sigma}_{c,m}^{2}(\mathbf{u})\neq 0$. Thus the fluctuation of a probabilistic modal coefficient propagates through its associated shared basis function to multiple output locations, so random fluctuations at different locations are not independent but jointly driven by the same low-dimensional random coefficients. 

\begin{figure}[!htbp]
	\centering
	\makebox[\textwidth][c]{%
		\includegraphics[
		width=1.13\textwidth,
		page=1,
		trim=18pt 18pt 18pt 18pt,
		clip
		]{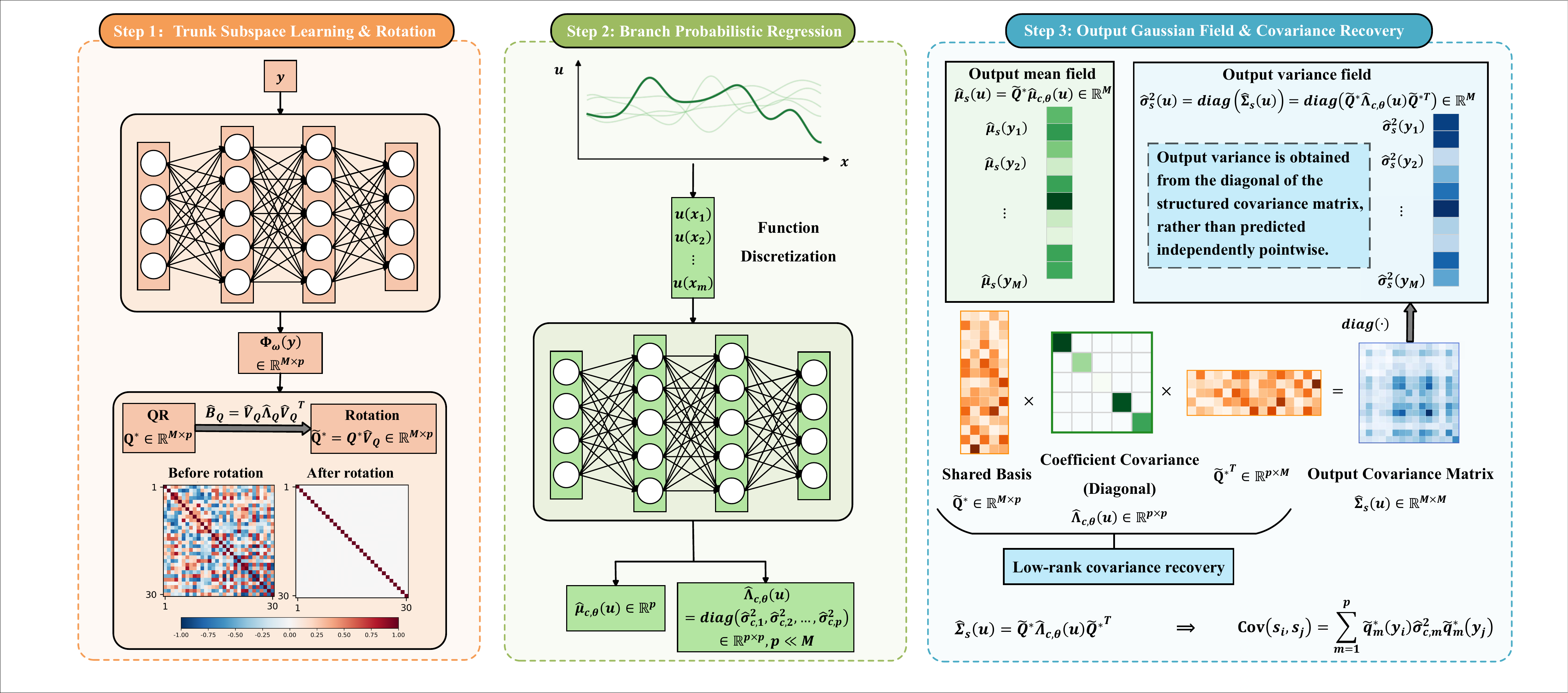}
	}
	\vspace{-4pt}
	\caption{Schematic of the two-step MV-DeepONet framework.
		Step 1 learns the output-space basis and applies QR
		orthogonalization and the subspace rotation
		\(\widetilde{\mathbf{Q}}^{*}
		=
		\mathbf{Q}^{*}\widehat{\mathbf{V}}_{Q}\).
		Step 2 predicts the coefficient mean
		\(\hat{\boldsymbol{\mu}}_{c,\theta}(\mathbf{u})\)
		and diagonal coefficient covariance
		\(\widehat{\boldsymbol{\Lambda}}_{c,\theta}(\mathbf{u})\).
		Finally, the resulting coefficient distribution is mapped back to the
		physical output space through the shared rotated basis, producing the
		output mean and structured conditional predictive covariance.}
	\label{fig:mvdeeponet_framework}
\end{figure}

Accordingly, the conditional mean and variance of the centered output at each output location are given by the following expressions, with the resulting network architecture shown in Fig.~\ref{fig:mvdeeponet_framework}.
\begin{equation}
	\mathbb{E}\left[
	\hat{x}_{\theta}(\mathbf{u})(y_i)
	\mid \mathbf{u}
	\right]
	=
	\sum_{m=1}^{p}
	\widetilde{q}_{m}^{*}(y_i)
	\hat{\mu}_{c,m}(\mathbf{u}),
	\qquad
	\operatorname{Var}\left[
	\hat{x}_{\theta}(\mathbf{u})(y_i)
	\mid \mathbf{u}
	\right]
	=
	\sum_{m=1}^{p}
	\left(
	\widetilde{q}_{m}^{*}(y_i)
	\right)^{2}
	\hat{\sigma}_{c,m}^{2}(\mathbf{u}).
\end{equation}

\FloatBarrier

By contrast, Prob-DeepONet represents predictive randomness independently at the physical output locations:
\begin{equation}
	\hat{s}_{i}(\mathbf{u})
	=
	\hat{\mu}_{\eta,i}(\mathbf{u})
	+
	\hat{\sigma}_{\eta,i}(\mathbf{u})\varepsilon_i,
	\qquad
	\varepsilon_i
	\mathrel{\overset{\mathrm{i.i.d.}}{
			\scalebox{1.25}{$\sim$}
	}}
	\mathcal{N}(0,1).
\end{equation}
Consequently,
\begin{equation}
	\operatorname{Cov}
	\left(
	\hat{s}_{i}(\mathbf{u}),
	\hat{s}_{j}(\mathbf{u})
	\mid \mathbf{u}
	\right)
	=
	0,
	\qquad
	i\neq j.
\end{equation}
The principal distinction between the two models therefore lies in the space where randomness is represented. Prob-DeepONet models conditionally independent randomness at individual output locations, whereas the proposed framework models low-dimensional random coefficients whose fluctuations are shared across the output domain through the basis functions.

At the dataset level, the law of total covariance decomposes the coefficient covariance into the covariance of the predicted coefficient means across inputs and the average conditional coefficient covariance. For \(K\) input samples, let \(\overline{\hat{\boldsymbol{\mu}}}_{c,\theta}
=
\frac{1}{K}\sum_{k=1}^{K}
\hat{\boldsymbol{\mu}}_{c,\theta}(\mathbf{u}^{(k)})\) denote the sample average of the predicted coefficient mean vectors.
The covariance induced by the branch mean head is
\begin{equation}
	\widehat{\mathbf{B}}_{\mu}
	=
	\frac{1}{K}
	\sum_{k=1}^{K}
	\left[
	\hat{\boldsymbol{\mu}}_{c,\theta}
	\left(
	\mathbf{u}^{(k)}
	\right)
	-
	\overline{\hat{\boldsymbol{\mu}}}_{c,\theta}
	\right]
	\left[
	\hat{\boldsymbol{\mu}}_{c,\theta}
	\left(
	\mathbf{u}^{(k)}
	\right)
	-
	\overline{\hat{\boldsymbol{\mu}}}_{c,\theta}
	\right]^{\top},
	\label{eq:b_mu_hat}
\end{equation}
whereas the average covariance predicted by the branch variance head is
\begin{equation}
	\widehat{\mathbf{B}}_{\sigma}
	=
	\frac{1}{K}
	\sum_{k=1}^{K}
	\widehat{\boldsymbol{\Lambda}}_{c,\theta}
	\left(
	\mathbf{u}^{(k)}
	\right).
	\label{eq:b_sigma_hat}
\end{equation}

The total predictive covariance in the coefficient space is therefore
\(\widehat{\mathbf{B}}_{\mathrm{2step}}
=
\widehat{\mathbf{B}}_{\mu}
+
\widehat{\mathbf{B}}_{\sigma}\), 
and the corresponding covariance in the physical output space is
\begin{equation}
	\widehat{\boldsymbol{\Sigma}}_{\mathrm{2step}}
	=
	\widetilde{\mathbf{Q}}^{*}
	\widehat{\mathbf{B}}_{\mathrm{2step}}
	\left(
	\widetilde{\mathbf{Q}}^{*}
	\right)^{\top}
	=
	\widetilde{\mathbf{Q}}^{*}
	\left(
	\widehat{\mathbf{B}}_{\mu}
	+
	\widehat{\mathbf{B}}_{\sigma}
	\right)
	\left(
	\widetilde{\mathbf{Q}}^{*}
	\right)^{\top}.
	\label{eq:dataset_total_covariance}
\end{equation}
Thus, dataset-level covariance recovery depends jointly on the predicted coefficient means and variances.

From a computational perspective, the proposed method does not require
direct parameterization of the full \(M\times M\) covariance matrix.
Instead, it predicts \(p\) conditional modal variances and constructs a
low-dimensional coefficient-space covariance, which is mapped to the
physical output space through the shared basis. When
\(p\ll M\), this yields a low-rank covariance representation at reduced
cost. Its diagonal entries provide pointwise uncertainty bands, while
the factorized form supports cross-location correlation analysis.

\subsubsection{Error Decomposition and Upper Bound for Output Covariance Recovery}
\label{sec:cov-error-decomposition}

To characterize the covariance recovery capability of the proposed model, this subsection derives a Frobenius-norm decomposition and a corresponding upper bound for the output covariance recovery error.

\paragraph{Basic Notation and Covariance Estimator}
\label{sec:cov-basic-notation}

Let \(\mathbf{U}\) denote the random input function and \(\mathbf{u}\) a generic realization of \(\mathbf{U}\), with corresponding discretized output field
\(\mathbf{s}(\mathbf{u})\in\mathbb{R}^{M}\), mean field \(\boldsymbol{\mu}_{s}=\mathbb{E}_{\mathbf{U}}[\mathbf{s}(\mathbf{U})]\),
and centered field \(\mathbf{x}(\mathbf{u})=\mathbf{s}(\mathbf{u})-\boldsymbol{\mu}_{s}\). The true input-induced output covariance is
\begin{equation}
	\boldsymbol{\Sigma}
	=
	\mathbb{E}_{\mathbf{U}}
	\left[
	\mathbf{x}(\mathbf{U})
	\mathbf{x}(\mathbf{U})^{\top}
	\right]
	\in
	\mathbb{R}^{M\times M}.
\end{equation}
Let the rotated trunk basis form the column-orthonormal matrix $\widetilde{\mathbf{Q}}^{*}=[\widetilde{\mathbf{q}}_{1}^{*},\ldots,\widetilde{\mathbf{q}}_{p}^{*}]\in\mathbb{R}^{M\times p}$ with $(\widetilde{\mathbf{Q}}^{*})^{\top}\widetilde{\mathbf{Q}}^{*}=\mathbf{I}_{p}$ and $p\ll M$, and let $\mathbf{P}=\widetilde{\mathbf{Q}}^{*}(\widetilde{\mathbf{Q}}^{*})^{\top}$ be the associated orthogonal projector. In this learned subspace, the true modal coefficients and their covariance are
\begin{equation}
	\mathbf{c}(\mathbf{u})
	=
	\left(
	\widetilde{\mathbf{Q}}^{*}
	\right)^{\top}
	\mathbf{x}(\mathbf{u})
	\in
	\mathbb{R}^{p},
	\qquad
	\mathbf{B}
	=
	\operatorname{Cov}(\mathbf{c})
	=
	\left(
	\widetilde{\mathbf{Q}}^{*}
	\right)^{\top}
	\boldsymbol{\Sigma}
	\widetilde{\mathbf{Q}}^{*}
	\in
	\mathbb{R}^{p\times p}.
\end{equation}

The dataset-level total predictive covariance follows from the law of total covariance:
\begin{equation}
	\operatorname{Cov}
	\left(
	\hat{\mathbf{x}}_{\theta}
	\right)
	=
	\operatorname{Cov}_{\mathbf{U}}
	\left(
	\mathbb{E}
	\left[
	\hat{\mathbf{x}}_{\theta}
	\mid
	\mathbf{U}
	\right]
	\right)
	+
	\mathbb{E}_{\mathbf{U}}
	\left[
	\operatorname{Cov}
	\left(
	\hat{\mathbf{x}}_{\theta}
	\mid
	\mathbf{U}
	\right)
	\right].
\end{equation}
With the trained basis $\widetilde{\mathbf{Q}}^{*}$ fixed, linearity of conditional expectation and the covariance propagation
rule for linear transformations give
\begin{equation}
	\mathbb{E}[\hat{\mathbf{x}}_{\theta}\mid\mathbf{u}^{(k)}]
	=
	\widetilde{\mathbf{Q}}^{*}\hat{\boldsymbol{\mu}}_{c,\theta}(\mathbf{u}^{(k)}),
	\qquad
	\operatorname{Cov}(\hat{\mathbf{x}}_{\theta}\mid\mathbf{u}^{(k)})
	=
	\widetilde{\mathbf{Q}}^{*}\widehat{\boldsymbol{\Lambda}}_{c,\theta}(\mathbf{u}^{(k)})(\widetilde{\mathbf{Q}}^{*})^{\top},
\end{equation}
from the coefficient distribution
$\hat{\mathbf{c}}_{\theta}(\mathbf{u}^{(k)})\sim
\mathcal{N}(\hat{\boldsymbol{\mu}}_{c,\theta}(\mathbf{u}^{(k)}),
\widehat{\boldsymbol{\Lambda}}_{c,\theta}(\mathbf{u}^{(k)}))$, with predicted mean
$\hat{\boldsymbol{\mu}}_{c,\theta}(\mathbf{u}^{(k)})\in\mathbb{R}^{p}$ and
diagonal covariance
$\widehat{\boldsymbol{\Lambda}}_{c,\theta}(\mathbf{u}^{(k)})
=\operatorname{diag}(\hat{\sigma}_{c,1}^{2}(\mathbf{u}^{(k)}),\ldots,
\hat{\sigma}_{c,p}^{2}(\mathbf{u}^{(k)}))$.

Taking the sample covariance of the conditional means and the sample average of the conditional covariances over the $K$ input samples gives, as established in the
preceding subsection, the dataset-level estimators
\begin{equation}
	\widehat{\mathbf{B}}_{\mathrm{2step}}
	=
	\widehat{\mathbf{B}}_{\mu}
	+
	\widehat{\mathbf{B}}_{\sigma},
	\qquad
	\widehat{\boldsymbol{\Sigma}}_{\mathrm{2step}}
	=
	\widetilde{\mathbf{Q}}^{*}
	\widehat{\mathbf{B}}_{\mathrm{2step}}
	\left(\widetilde{\mathbf{Q}}^{*}\right)^{\top},
\end{equation}
where $\widehat{\mathbf{B}}_{\mu}$ and $\widehat{\mathbf{B}}_{\sigma}$ are the mean-head and variance-head coefficient covariances defined in Eqs.~\eqref{eq:b_mu_hat}--\eqref{eq:b_sigma_hat}. Their sum defines the dataset-level coefficient covariance analyzed below.

\paragraph{Basic Error Decomposition}
\label{sec:cov-basic-error-decomposition}

The recovery error is measured in the Frobenius norm as
$\|\boldsymbol{\Sigma}-\widehat{\boldsymbol{\Sigma}}_{\mathrm{2step}}\|_{F}$.
Introducing the projected covariance
$\mathbf{P}\boldsymbol{\Sigma}\mathbf{P}
=\widetilde{\mathbf{Q}}^{*}\mathbf{B}(\widetilde{\mathbf{Q}}^{*})^{\top}$ as an
intermediate term gives
\begin{equation}
	\boldsymbol{\Sigma}
	-
	\widehat{\boldsymbol{\Sigma}}_{\mathrm{2step}}
	=
	\left(
	\boldsymbol{\Sigma}
	-
	\mathbf{P}
	\boldsymbol{\Sigma}
	\mathbf{P}
	\right)
	+
	\widetilde{\mathbf{Q}}^{*}
	\left(
	\mathbf{B}
	-
	\widehat{\mathbf{B}}_{\mathrm{2step}}
	\right)
	\left(
	\widetilde{\mathbf{Q}}^{*}
	\right)^{\top}.
\end{equation}
Because $\widetilde{\mathbf{Q}}^{*}$ has orthonormal columns, the triangle inequality yields
\begin{equation}
	\left\|
	\boldsymbol{\Sigma}
	-
	\widehat{\boldsymbol{\Sigma}}_{\mathrm{2step}}
	\right\|_{F}
	\leq
	\left\|
	\boldsymbol{\Sigma}
	-
	\mathbf{P}
	\boldsymbol{\Sigma}
	\mathbf{P}
	\right\|_{F}
	+
	\left\|
	\mathbf{B}
	-
	\widehat{\mathbf{B}}_{\mathrm{2step}}
	\right\|_{F}.
\end{equation}
Here,
\(\left\|
\boldsymbol{\Sigma}
-
\mathbf{P}\boldsymbol{\Sigma}\mathbf{P}
\right\|_{F}\)
denotes the subspace projection error, while
\(\left\|
\mathbf{B}
-
\widehat{\mathbf{B}}_{\mathrm{2step}}
\right\|_{F}\)
denotes the total coefficient covariance estimation error.

Consider the eigendecomposition
\begin{equation}
	\boldsymbol{\Sigma}
	=
	\mathbf{W}\boldsymbol{\Lambda}\mathbf{W}^{\top},
	\qquad
	\boldsymbol{\Lambda}
	=
	\operatorname{diag}(\lambda_{1},\ldots,\lambda_{M}),
	\qquad
	\lambda_{1}\geq\cdots\geq\lambda_{M}\geq0.
\end{equation}
Let $\mathbf{W}_{p}\in\mathbb{R}^{M\times p}$ collect the leading $p$ eigenvectors, and $\mathbf{P}^{\star}=\mathbf{W}_{p}\mathbf{W}_{p}^{\top}$
be the projector onto the optimal rank-$p$ principal subspace. By the Eckart--Young--Mirsky theorem~\cite{EckartYoung1936},
\begin{equation}
	\left\|
	\boldsymbol{\Sigma}
	-
	\mathbf{P}^{\star}\boldsymbol{\Sigma}\mathbf{P}^{\star}
	\right\|_{F}
	=
	\sqrt{
		\sum_{j>p}\lambda_{j}^{2}
	}.
\end{equation}
Introducing \(\mathbf{P}^{\star}\boldsymbol{\Sigma}\mathbf{P}^{\star}\) as an intermediate term, the projection error can be decomposed as
\begin{equation}
	\left\|
	\boldsymbol{\Sigma}
	-
	\mathbf{P}\boldsymbol{\Sigma}\mathbf{P}
	\right\|_{F}
	\leq
	\sqrt{
		\sum_{j>p}\lambda_{j}^{2}
	}
	+
	\left\|
	\mathbf{P}^{\star}
	\boldsymbol{\Sigma}
	\mathbf{P}^{\star}
	-
	\mathbf{P}
	\boldsymbol{\Sigma}
	\mathbf{P}
	\right\|_{F}.
\end{equation}
Applying the mixed-norm inequality to the second term, we obtain
\begin{equation}
	\left\|
	\mathbf{P}^{\star}\boldsymbol{\Sigma}\mathbf{P}^{\star}
	-
	\mathbf{P}\boldsymbol{\Sigma}\mathbf{P}
	\right\|_{F}
	\leq
	\left\|
	\mathbf{P}^{\star}-\mathbf{P}
	\right\|_{2}
	\left\|
	\boldsymbol{\Sigma}
	\right\|_{F}
	\left\|
	\mathbf{P}^{\star}
	\right\|_{2}
	+
	\left\|
	\mathbf{P}
	\right\|_{2}
	\left\|
	\boldsymbol{\Sigma}
	\right\|_{F}
	\left\|
	\mathbf{P}^{\star}-\mathbf{P}
	\right\|_{2}
	\leq
	2
	\left\|
	\boldsymbol{\Sigma}
	\right\|_{F}
	\left\|
	\mathbf{P}-\mathbf{P}^{\star}
	\right\|_{F}.
	\label{eq:projected_covariance_difference}
\end{equation}
Substituting this estimate into the preceding bound gives
\begin{equation}
	\left\|
	\boldsymbol{\Sigma}
	-
	\mathbf{P}
	\boldsymbol{\Sigma}
	\mathbf{P}
	\right\|_{F}
	\leq
	\sqrt{
		\sum_{j>p}\lambda_{j}^{2}
	}
	+
	2
	\left\|
	\boldsymbol{\Sigma}
	\right\|_{F}
	\left\|
	\mathbf{P}-\mathbf{P}^{\star}
	\right\|_{F}.
\end{equation}
Consequently,
\begin{equation}
	\left\|
	\boldsymbol{\Sigma}
	-
	\widehat{\boldsymbol{\Sigma}}_{\mathrm{2step}}
	\right\|_{F}
	\leq
	\sqrt{
		\sum_{j>p}\lambda_{j}^{2}
	}
	+
	2\|\boldsymbol{\Sigma}\|_{F}
	\|\mathbf{P}-\mathbf{P}^{\star}\|_{F}
	+
	\left\|
	\mathbf{B}
	-
	\widehat{\mathbf{B}}_{\mathrm{2step}}
	\right\|_{F}.
	\label{eq:app-cov-basic-bound}
\end{equation}
For convenience, let \(T_{1}\), \(T_{2}\), and \(T_{3}\) represent the first, second, and third terms on the right-hand side, respectively. Specifically, \(T_{1}\) is the low-rank truncation error of the true covariance, \(T_{2}\) denotes the discrepancy between the trunk-learned subspace and the optimal rank-\(p\) principal subspace, and \(T_{3}\) quantifies the error in the dataset-level total coefficient covariance.

\paragraph{Trunk Subspace Error Term}
\label{sec:cov-trunk-subspace-error}

We next analyze \(T_{2}\), which is governed by the subspace discrepancy \(\|\mathbf{P}-\mathbf{P}^{\star}\|_{F}\), with \(\|\boldsymbol{\Sigma}\|_{F}\) serving as a problem-dependent scale factor. Let $\mathbf{X}=[\mathbf{x}^{(1)},\ldots,\mathbf{x}^{(K)}]\in\mathbb{R}^{M\times K}$ be the centered training-output matrix, and let $\widehat{\boldsymbol{\Sigma}}=\frac{1}{K}\mathbf{X}\mathbf{X}^{\top}$ be the
empirical covariance, with eigendecomposition
\begin{equation}
	\widehat{\boldsymbol{\Sigma}}
	=
	\widehat{\mathbf{W}}
	\widehat{\boldsymbol{\Lambda}}
	\widehat{\mathbf{W}}^{\top},
	\qquad
	\widehat{\boldsymbol{\Lambda}}
	=
	\operatorname{diag}
	\left(
	\widehat{\lambda}_{1},\ldots,\widehat{\lambda}_{M}
	\right),
	\qquad
	\widehat{\lambda}_{1}\geq\cdots\geq\widehat{\lambda}_{M}\geq0.
\end{equation}
Let $\widehat{\mathbf{W}}_{p}\in\mathbb{R}^{M\times p}$ collect the
leading $p$ eigenvectors of
$\widehat{\boldsymbol{\Sigma}}$,
and define
\(\widehat{\mathbf{P}}^{\star}
=\widehat{\mathbf{W}}_{p}\widehat{\mathbf{W}}_{p}^{\top}\)
as the projector onto the empirically optimal \(p\)-dimensional principal subspace. Then
\begin{equation}
	\|\mathbf{P}-\mathbf{P}^{\star}\|_{F}
	\leq
	\|\mathbf{P}-\widehat{\mathbf{P}}^{\star}\|_{F}
	+
	\|\widehat{\mathbf{P}}^{\star}-\mathbf{P}^{\star}\|_{F},
\end{equation}
where the two terms represent trunk learning error and finite-sample
statistical error, respectively.

\subparagraph{Trunk Learning Error}
\label{sec:cov-trunk-learning-error}

Define the reconstruction risk of a rank-$p$ projector $\mathbf{P}$ as
\begin{equation}
	\mathcal{R}_{T}(\mathbf{P})
	=
	\frac{1}{K}\|\mathbf{X}-\mathbf{P}\mathbf{X}\|_{F}^{2}
	=
	\operatorname{tr}
	\left[
	(\mathbf{I}_{M}-\mathbf{P})
	\widehat{\boldsymbol{\Sigma}}
	\right].
\end{equation}
The reconstruction risk of the empirically optimal principal subspace is similarly expressed through its orthogonal projector
$\widehat{\mathbf{P}}^{\star}$ as
\begin{equation}
	\mathcal{R}_{T}
	\left(
	\widehat{\mathbf{P}}^{\star}
	\right)
	=
	\operatorname{tr}
	\left[
	\left(
	\mathbf{I}_{M}
	-
	\widehat{\mathbf{P}}^{\star}
	\right)
	\widehat{\boldsymbol{\Sigma}}
	\right].
\end{equation}
Define the excess reconstruction risk of the trunk-learned subspace relative to the empirical optimum as
\begin{equation}
	\Delta_{T}
	=
	\mathcal{R}_{T}(\mathbf{P})
	-
	\mathcal{R}_{T}(\widehat{\mathbf{P}}^{\star})
	=
	\operatorname{tr}
	\left(
	\widehat{\mathbf{P}}^{\star}
	\widehat{\boldsymbol{\Sigma}}
	\right)
	-
	\operatorname{tr}
	\left(
	\mathbf{P}\widehat{\boldsymbol{\Sigma}}
	\right)
	\geq0.
\end{equation}

To relate $\Delta_{T}$ to the discrepancy between the trunk-learned and empirical principal subspaces, express $\mathbf{P}$ in the eigenbasis of
$\widehat{\boldsymbol{\Sigma}}$ by defining $\widetilde{\mathbf{P}}=\widehat{\mathbf{W}}^{\top}\mathbf{P}\widehat{\mathbf{W}}$. Using the cyclic property of the trace gives
\begin{equation}
	\operatorname{tr}\left(\mathbf{P}\widehat{\boldsymbol{\Sigma}}\right)
	=
	\operatorname{tr}\left(\mathbf{P}\widehat{\mathbf{W}}\widehat{\boldsymbol{\Lambda}}\widehat{\mathbf{W}}^{\top}\right)
	=
	\operatorname{tr}\left(\widehat{\mathbf{W}}^{\top}\mathbf{P}\widehat{\mathbf{W}}\widehat{\boldsymbol{\Lambda}}\right)
	=
	\operatorname{tr}\left(\widehat{\boldsymbol{\Lambda}}\widetilde{\mathbf{P}}\right)
	=
	\sum_{j=1}^{M}\widehat{\lambda}_{j}\widetilde{p}_{jj}.
\end{equation}
Let \(\widehat{\gamma}_{p}=\widehat{\lambda}_{p}-\widehat{\lambda}_{p+1}>0\) denote the empirical spectral gap at the selected dimension. 

The empirically optimal projector can then be represented in the eigenbasis of \(\widehat{\boldsymbol{\Sigma}}\) as
\(\widehat{\mathbf{P}}^{\star}
=
\widehat{\mathbf{W}}\mathbf{E}_{p}\widehat{\mathbf{W}}^{\top}\),
with
\(\mathbf{E}_{p}
=
\left[
\begin{smallmatrix}
	\mathbf{I}_{p} & \mathbf{0} \\
	\mathbf{0} & \mathbf{0}
\end{smallmatrix}
\right]
\in\mathbb{R}^{M\times M}\).
Consequently,
\begin{equation}
	\operatorname{tr}\left(\widehat{\mathbf{P}}^{\star}\widehat{\boldsymbol{\Sigma}}\right)
	=
	\operatorname{tr}\left(\widehat{\mathbf{W}}\mathbf{E}_{p}\widehat{\mathbf{W}}^{\top}\widehat{\mathbf{W}}\widehat{\boldsymbol{\Lambda}}\widehat{\mathbf{W}}^{\top}\right)
	=
	\operatorname{tr}\left(\mathbf{E}_{p}\widehat{\boldsymbol{\Lambda}}\right)
	=
	\sum_{j=1}^{p}\widehat{\lambda}_{j}.
\end{equation}
It follows that
\(\Delta_{T}
=
\sum_{j=1}^{p}\widehat{\lambda}_{j}
(1-\widetilde{p}_{jj})
-
\sum_{j=p+1}^{M}\widehat{\lambda}_{j}\widetilde{p}_{jj}\).
Using the eigenvalue ordering of \(\widehat{\boldsymbol{\Sigma}}\) yields
\begin{equation}
	\Delta_{T}
	\geq
	\widehat{\lambda}_{p}
	\sum_{j=1}^{p}
	\left(
	1-\widetilde{p}_{jj}
	\right)
	-
	\widehat{\lambda}_{p+1}
	\sum_{j=p+1}^{M}
	\widetilde{p}_{jj}
	=
	\left(
	\widehat{\lambda}_{p}
	-
	\widehat{\lambda}_{p+1}
	\right)
	\sum_{j=1}^{p}
	\left(
	1-\widetilde{p}_{jj}
	\right)
	=
	\widehat{\gamma}_{p}
	\sum_{j=1}^{p}
	\left(
	1-\widetilde{p}_{jj}
	\right).
\end{equation}
Since $\mathbf{P}$ and $\widehat{\mathbf{P}}^{\star}$ are rank-$p$ orthogonal projectors,
\begin{equation}
	\left\|
	\mathbf{P}
	-
	\widehat{\mathbf{P}}^{\star}
	\right\|_{F}^{2}
	=
	\operatorname{tr}
	\left[
	\left(
	\mathbf{P}
	-
	\widehat{\mathbf{P}}^{\star}
	\right)^{\top}
	\left(
	\mathbf{P}
	-
	\widehat{\mathbf{P}}^{\star}
	\right)
	\right]
	=
	2
	\left[
	p
	-
	\operatorname{tr}
	\left(
	\mathbf{P}
	\widehat{\mathbf{P}}^{\star}
	\right)
	\right].
\end{equation}
Moreover,
\(\operatorname{tr}
(\mathbf{P}\widehat{\mathbf{P}}^{\star})
=
\operatorname{tr}
(\mathbf{P}\widehat{\mathbf{W}}\mathbf{E}_{p}
\widehat{\mathbf{W}}^{\top})
=
\operatorname{tr}
(\widehat{\mathbf{W}}^{\top}\mathbf{P}
\widehat{\mathbf{W}}\mathbf{E}_{p})
=
\operatorname{tr}
(\widetilde{\mathbf{P}}\mathbf{E}_{p})
=
\sum_{j=1}^{p}\widetilde{p}_{jj}\).
Therefore,
\begin{equation}
	\left\|
	\mathbf{P}
	-
	\widehat{\mathbf{P}}^{\star}
	\right\|_{F}^{2}
	=
	2
	\sum_{j=1}^{p}
	\left(
	1-\widetilde{p}_{jj}
	\right)
	\leq
	\frac{2\Delta_{T}}{\widehat{\gamma}_{p}}.
\end{equation}
Equivalently,
\[
\left\|
\mathbf{P}
-
\widehat{\mathbf{P}}^{\star}
\right\|_{F}
\leq
\sqrt{
	\frac{2\Delta_{T}}
	{\widehat{\gamma}_{p}}
}.
\]

This inequality shows that a small excess reconstruction risk implies a small discrepancy between the trunk-learned and empirical principal subspaces, provided that \(\widehat{\gamma}_{p}\) does not vanish. This conclusion is further supported by Theorem~3.5 of Lee and Shin~\cite{TwoStepDeepONet_LeeShin2024}, which establishes that a sufficiently expressive trunk network can attain the best rank-\(p\) approximation error of the training-output matrix. Together with the spectral-gap bound above, this result provides theoretical justification for the trunk network's capacity to recover the empirically optimal \(p\)-dimensional principal subspace.

\subparagraph{Finite-Sample Statistical Error}
\label{sec:cov-finite-sample-error}

The finite-sample statistical error arises from the covariance perturbation \(\widehat{\boldsymbol{\Sigma}}-\boldsymbol{\Sigma}\), where \(\widehat{\boldsymbol{\Sigma}}
=
\boldsymbol{\Sigma}
+
(\widehat{\boldsymbol{\Sigma}}-\boldsymbol{\Sigma})\).
Under a non-vanishing true spectral gap \(\gamma_{p}=\lambda_{p}-\lambda_{p+1}>0\), Theorem~2 of Yu \textit{et al.}~\cite{YuWangSamworth2015}, a variant of the Davis--Kahan \(\sin\Theta\) theorem based on the population spectral gap, gives the following bound for the leading \(p\)-dimensional eigenspaces:
\begin{equation}
	\left\|\sin\Theta\left(\mathbf{W}_{p},\widehat{\mathbf{W}}_{p}\right)\right\|_{F}
	\leq
	\frac{2}{\gamma_{p}}
	\min\left(
	\sqrt{p}\left\|\widehat{\boldsymbol{\Sigma}}-\boldsymbol{\Sigma}\right\|_{2},
	\left\|\widehat{\boldsymbol{\Sigma}}-\boldsymbol{\Sigma}\right\|_{F}
	\right)
	\leq
	\frac{2}{\gamma_{p}}
	\left\|\widehat{\boldsymbol{\Sigma}}-\boldsymbol{\Sigma}\right\|_{F}.
\end{equation}

Combined with the identity relating the Frobenius distance between the orthogonal projectors to the principal angles,
\(\|\widehat{\mathbf{P}}^{\star}-\mathbf{P}^{\star}\|_{F}
=
\sqrt{2}\|\sin\Theta(\mathbf{W}_{p},\widehat{\mathbf{W}}_{p})\|_{F}\),
this yields
\begin{equation}
	\left\|
	\widehat{\mathbf{P}}^{\star}
	-
	\mathbf{P}^{\star}
	\right\|_{F}
	\leq
	\frac{
		2\sqrt{2}
		\left\|
		\widehat{\boldsymbol{\Sigma}}
		-
		\boldsymbol{\Sigma}
		\right\|_{F}
	}{
		\gamma_{p}
	}.
\end{equation}
Substituting this estimate and the preceding trunk-learning bound into $T_{2}$ gives
\begin{equation}
	T_{2}
	\leq
	2
	\left\|
	\boldsymbol{\Sigma}
	\right\|_{F}
	\sqrt{
		\frac{
			2\Delta_{T}
		}{
			\widehat{\gamma}_{p}
		}
	}
	+
	\frac{
		4\sqrt{2}
		\left\|
		\boldsymbol{\Sigma}
		\right\|_{F}
	}{
		\gamma_{p}
	}
	\left\|
	\widehat{\boldsymbol{\Sigma}}
	-
	\boldsymbol{\Sigma}
	\right\|_{F}.
\end{equation}
Thus, $T_{2}$ comprises the trunk-learning error relative to the empirical principal subspace and the finite-sample statistical error.

\paragraph{Learning Error in the Branch-Induced Total Coefficient Covariance}
\label{sec:cov-branch-covariance-error}

To bound \(T_{3}\), we introduce the empirical coefficient covariance
$\widehat{\mathbf{B}}=(\widetilde{\mathbf{Q}}^{*})^{\top}\widehat{\boldsymbol{\Sigma}}\widetilde{\mathbf{Q}}^{*}\in\mathbb{R}^{p\times p}$
as an intermediate quantity. The triangle inequality gives
\begin{equation}
	\left\|
	\mathbf{B}
	-
	\widehat{\mathbf{B}}_{\mathrm{2step}}
	\right\|_{F}
	\leq
	\left\|
	\mathbf{B}
	-
	\widehat{\mathbf{B}}
	\right\|_{F}
	+
	\left\|
	\widehat{\mathbf{B}}
	-
	\widehat{\mathbf{B}}_{\mathrm{2step}}
	\right\|_{F}.
\end{equation}
Here, the first term represents the finite-sample statistical error, whereas the second term quantifies the branch learning error. 
Since $\mathbf{B}-\widehat{\mathbf{B}}=(\widetilde{\mathbf{Q}}^{*})^{\top}(\boldsymbol{\Sigma}-\widehat{\boldsymbol{\Sigma}})\widetilde{\mathbf{Q}}^{*}$, the first term satisfies $\|\mathbf{B}-\widehat{\mathbf{B}}\|_{F}\leq\|\boldsymbol{\Sigma}-\widehat{\boldsymbol{\Sigma}}\|_{F}$. 

Let
$\widehat{\mathbf{B}}_{\mathrm{2step}}^{\mathrm{NN,best}}$ denote the
dataset-level covariance induced by the branch model in the considered
function class that best approximates $\widehat{\mathbf{B}}$, and define
\begin{equation}
	\eta_{\mathrm{app}}^{B}
	=
	\left\|
	\widehat{\mathbf{B}}
	-
	\widehat{\mathbf{B}}_{\mathrm{2step}}^{\mathrm{NN,best}}
	\right\|_{F},
	\qquad
	\eta_{\mathrm{opt}}^{B}
	=
	\left\|
	\widehat{\mathbf{B}}_{\mathrm{2step}}^{\mathrm{NN,best}}
	-
	\widehat{\mathbf{B}}_{\mathrm{2step}}
	\right\|_{F}.
\end{equation}
Here, \(\eta_{\mathrm{app}}^{B}\) measures the approximation error of the
branch-network function class relative to the empirical target, whereas
\(\eta_{\mathrm{opt}}^{B}\) measures the optimization-related discrepancy
between this benchmark and the trained model. The former depends on network expressivity and target complexity; the
latter is influenced by training configuration. Therefore,
\begin{equation}
	T_{3}
	\leq
	\|\boldsymbol{\Sigma}-\widehat{\boldsymbol{\Sigma}}\|_{F}
	+
	\eta_{\mathrm{app}}^{B}
	+
	\eta_{\mathrm{opt}}^{B}.
	\label{eq:branch-covariance-bound}
\end{equation}
Thus, the total coefficient covariance error is governed jointly by finite-sample covariance estimation and branch probabilistic regression error.

\paragraph{Combined Error Bound}
\label{sec:cov-combined-error-bound}

Combining the preceding results gives
\begin{equation}
	\left\|
	\boldsymbol{\Sigma}
	-
	\widehat{\boldsymbol{\Sigma}}_{\mathrm{2step}}
	\right\|_{F}
	\leq
	\sqrt{
		\sum_{j>p}\lambda_{j}^{2}
	}
	+
	2\|\boldsymbol{\Sigma}\|_{F}
	\sqrt{
		\frac{2\Delta_{T}}{\widehat{\gamma}_{p}}
	}
	+
	\left(
	\frac{4\sqrt{2}\|\boldsymbol{\Sigma}\|_{F}}{\gamma_{p}}
	+1
	\right)
	\left\|
	\widehat{\boldsymbol{\Sigma}}
	-
	\boldsymbol{\Sigma}
	\right\|_{F}
	+
	\eta_{\mathrm{app}}^{B}
	+
	\eta_{\mathrm{opt}}^{B}.
	\label{eq:cov-final-bound-main}
\end{equation}
Under standard covariance concentration conditions, the empirical
covariance satisfies the spectral-norm convergence rate $\|\widehat{\boldsymbol{\Sigma}}-\boldsymbol{\Sigma}\|_{2}
=\mathcal{O}_{\mathbb{P}}(K^{-1/2})$~\cite{KoltchinskiiLounici2017}. For any
\(\mathbf{A}\in\mathbb{R}^{M\times M}\), \(\|\mathbf{A}\|_{F}\leq\sqrt{M}\|\mathbf{A}\|_{2}\). Therefore,
\begin{equation}
	\left\|
	\widehat{\boldsymbol{\Sigma}}
	-
	\boldsymbol{\Sigma}
	\right\|_{F}  
	=
	\mathcal{O}_{\mathbb{P}}
	\left(
	\sqrt{\frac{M}{K}} 
	\right).
\end{equation}
Hence,
\begin{equation}
	\left\|
	\boldsymbol{\Sigma}
	-
	\widehat{\boldsymbol{\Sigma}}_{\mathrm{2step}}
	\right\|_{F}
	\leq
	\sqrt{
		\sum_{j>p}\lambda_{j}^{2}
	}
	+
	2\|\boldsymbol{\Sigma}\|_{F}
	\sqrt{
		\frac{2\Delta_{T}}{\widehat{\gamma}_{p}}
	}
	+
	\mathcal{O}_{\mathbb{P}}
	\left[
	\left(
	\frac{4\sqrt{2}\|\boldsymbol{\Sigma}\|_{F}}{\gamma_{p}}
	+1
	\right)
	\sqrt{\frac{M}{K}}
	\right]
	+
	\eta_{\mathrm{app}}^{B}
	+
	\eta_{\mathrm{opt}}^{B}.
	\label{eq:cov-final-bound-rate}
\end{equation}
  
The bound identifies four principal factors governing covariance recovery.
First, the spectral-tail term \(\sqrt{\sum_{j>p}\lambda_{j}^{2}}\) characterizes the low-rank compressibility of the true covariance. Second, the trunk-related contribution is determined by the excess reconstruction risk \(\Delta_{T}\), with \(\widehat{\gamma}_{p}\) linking this risk to the discrepancy between the learned and empirical principal subspaces. Third, finite-sample statistical error depends on the number of training samples, the output dimension, and the distribution of the training outputs. Fourth, \(\eta_{\mathrm{app}}^{B}+\eta_{\mathrm{opt}}^{B}\) accounts for the branch network's approximation and optimization errors in learning the empirical total coefficient covariance. Additionally, the spectral gaps $\widehat{\gamma}_{p}$ and $\gamma_{p}$ control the stability of the corresponding bounds. The true spectral gap is assumed to be non-vanishing, whereas the empirical spectral gap is verified numerically in each example. The factor $\|\boldsymbol{\Sigma}\|_{F}$ reflects the overall covariance scale.

In summary, accurate and stable covariance recovery relies on rapid spectral decay of the true covariance, non-vanishing spectral gaps, controlled trunk and branch errors, and sufficiently accurate finite-sample covariance estimation.

\section{Results and Discussions}
\label{sec:results}

To ensure consistent evaluation across the case studies, a unified set of metrics is used to assess mean prediction, covariance recovery, spatial correlation, and prediction intervals. As a supplementary enhancement module, deep ensembles and post-hoc uncertainty calibration are introduced to account for epistemic uncertainty and improve interval coverage. The low-dimensional coefficient-space formulation keeps the additional cost of finite ensembles computationally manageable.

\begin{enumerate}[label=(\arabic*),leftmargin=2em,itemsep=0.6em]
	
	\item \textbf{Normalized Frobenius covariance error.}
	Covariance estimation error is commonly evaluated under the operator and
	Frobenius norms, which characterize different aspects of estimation
	accuracy and lead to different optimal procedures~\cite{CovEst_Cai2010}.
	Since this study focuses on the overall covariance structure, the
	normalized Frobenius error is adopted:
	\begin{equation}
		E_{\Sigma}
		=
		\frac{
			\left\|
			\widehat{\boldsymbol{\Sigma}}_{\mathrm{pred}}
			-
			\widehat{\boldsymbol{\Sigma}}_{\mathrm{ref}}
			\right\|_{F}
		}{
			\left\|
			\widehat{\boldsymbol{\Sigma}}_{\mathrm{ref}}
			\right\|_{F}
		},
		\label{eq:metric-cov-frobenius}
	\end{equation}
	where \(\widehat{\boldsymbol{\Sigma}}_{\mathrm{ref}}\) and
	\(\widehat{\boldsymbol{\Sigma}}_{\mathrm{pred}}\) denote the reference
	empirical covariance and model prediction, respectively. Following the decomposition in Eq.~\eqref{eq:app-cov-basic-bound}, three
	error curves are compared over \(p\):
	
	\begin{enumerate}[label=\alph*),leftmargin=2em,itemsep=0.4em]
		
		\item \textbf{Oracle singular value decomposition (SVD) curve.}
		The optimal rank-\(p\) approximation
		\(\widehat{\boldsymbol{\Sigma}}_{\mathrm{SVD}}^{(p)}\) retains the
		leading \(p\) modes of
		\(\widehat{\boldsymbol{\Sigma}}_{\mathrm{ref}}\) and represents
		\(T_{1}\).
		
		\item \textbf{Learned-subspace oracle curve.}
		The reference covariance is projected onto the trunk-learned subspace as
		\(\widetilde{\mathbf{Q}}^{*}
		[(\widetilde{\mathbf{Q}}^{*})^{\top}
		\widehat{\boldsymbol{\Sigma}}_{\mathrm{ref}}
		\widetilde{\mathbf{Q}}^{*}]
		(\widetilde{\mathbf{Q}}^{*})^{\top}\).
		Using the reference coefficient covariance excludes branch error, so
		this curve represents \(T_{1}+T_{2}\).
		
		\item \textbf{Final model prediction curve.}
		The complete model estimates
		\(\widehat{\boldsymbol{\Sigma}}_{\mathrm{2step}}
		=
		\widetilde{\mathbf{Q}}^{*}
		(\widehat{\mathbf{B}}_{\mu}
		+
		\widehat{\mathbf{B}}_{\sigma})
		(\widetilde{\mathbf{Q}}^{*})^{\top}\),
		which incorporates \(T_{1}+T_{2}+T_{3}\).
		
	\end{enumerate}
		
	\item \textbf{Multi-anchor correlation plots.}
	These plots assess off-diagonal dependence across the output domain. The
	correlation coefficient is
	\begin{equation}
		\rho_{ij}
		=
		\frac{\Sigma_{ij}}
		{\sqrt{\Sigma_{ii}\Sigma_{jj}}}.
		\label{eq:metric-correlation}
	\end{equation}
	Reference and predicted correlation profiles are compared at several
	anchor points to evaluate the recovered correlation range, decay, and
	local structure.
	
	\item \textbf{Prediction interval coverage and width.}
	Prediction interval coverage probability (PICP) and mean prediction
	interval width (MPIW) measure calibration and sharpness,
	respectively~\cite{PIReview_Khosravi2011}:
	\begin{equation}
		\begin{aligned}
			\mathrm{PICP}_{\gamma}
			&=
			\frac{1}{N_{\mathrm{test}}M}
			\sum_{i=1}^{N_{\mathrm{test}}}
			\sum_{m=1}^{M}
			\mathbb{I}
			\left(
			s_m^{(i)}
			\in
			\left[
			L_{\gamma,m}^{(i)},
			U_{\gamma,m}^{(i)}
			\right]
			\right),
			\qquad
			\mathrm{MPIW}_{\gamma}
			&=
			\frac{1}{N_{\mathrm{test}}M}
			\sum_{i=1}^{N_{\mathrm{test}}}
			\sum_{m=1}^{M}
			\left(
			U_{\gamma,m}^{(i)}
			-
			L_{\gamma,m}^{(i)}
			\right).
		\end{aligned}
		\label{eq:metric-picp-mpiw}
	\end{equation}
	Here, \(N_{\mathrm{test}}\) is the number of test samples,
	\(\mathbb{I}(\cdot)\) is the indicator function, and the interval is
	defined by
	\(L_{\gamma,m}^{(i)}
	=
	\hat{\mu}_{m}^{(i)}
	-
	\kappa_{\gamma}\hat{\sigma}_{m}^{(i)}\) and
	\(U_{\gamma,m}^{(i)}
	=
	\hat{\mu}_{m}^{(i)}
	+
	\kappa_{\gamma}\hat{\sigma}_{m}^{(i)}\).
	Under the Gaussian assumption,
	\(\kappa_{\gamma}=z_{(1+\gamma)/2}\); under split conformal prediction (CP),
	it is determined by the calibration-score quantile. PICP should approach
	\(\gamma\), while a smaller MPIW indicates sharper intervals at comparable
	coverage.
	
\end{enumerate}

\subsection{Reaction--Diffusion Equation}
\label{sec:dr-introduction}

Reaction--diffusion equations describe the spatiotemporal evolution of concentration fields governed by diffusion and local reactions~\cite{RD_Fife1979}, with applications in chemical reaction engineering~\cite{RD_Aris1975}, biological morphogenesis~\cite{Morphogenesis_Turing1990}, and spatial ecology~\cite{SpatialEcology_CantrellCosner2004}. We consider
\begin{equation}
	\begin{aligned}
		&\frac{\partial s}{\partial t}
		=
		\nu\frac{\partial^2 s}{\partial x^2}
		+
		ks^2  
		+
		u(x),
		&&
		(x,t)\in[0,1]\times[0,1],
		\\
		&s(x,0)=0,
		&&
		x\in[0,1],
		\\
		&s(0,t)=s(1,t)=0,
		&&
		t\in[0,1].
	\end{aligned}
	\label{eq:dr-pde}
\end{equation}
where $\nu$ and $k$ are the diffusion and reaction coefficients, respectively, and are both set to $0.01$. The random source function \(u(x)\) is mapped to the spatiotemporal solution field \(s(x,t)\) through
\begin{equation}
	\mathcal{G}:u(x)\mapsto s(x,t).
	\label{eq:dr-operator}
\end{equation}

\vspace{4pt}
\par\noindent\textbf{Sources of Uncertainty.}
Random source terms may arise from fluctuating generation rates, spatial heterogeneity, and environmental perturbations. Examples include
stochastic gene expression in morphogen formation~\cite{NoisyMorphogen_EnglandCardy2005} and randomly distributed catalytic sources in systems with disordered kinetics~\cite{RandomChannel_Vlad2010}. Their propagation through the governing dynamics induces uncertainty and correlation across the output domain. Accordingly, this example evaluates the recovery of the mean response and output correlation structure under a random input function \(u(x)\).

\vspace{4pt}
\par\noindent\textbf{Data Generation.}
The source function is sampled from a Gaussian process with the squared-exponential covariance kernel
\begin{equation}
	K(x,x')
	=
	\exp\left(
	-\frac{|x-x'|^2}{2\ell^2}
	\right),
	\label{eq:dr-kernel}
\end{equation}
where \(\ell\) is the correlation length. For each realization, Eq.~\eqref{eq:dr-pde} is solved by a finite-difference scheme with the spatial and temporal domains each discretized into 100 uniformly spaced points, producing a \(100\times100\) response field. A total of 2,200 samples are generated, including 1,000 training samples, 1,000 test samples, and 200 calibration samples.

\subsubsection{Prediction Accuracy and Generalization Performance}
\label{sec:dr-accuracy}

\vspace{4pt}
\noindent\textbf{Training and In-Distribution Accuracy}
\label{subsec:dr-accuracy}

\begin{figure}[!htbp]
	\centering
	\includegraphics[
	width=0.5\linewidth
	]{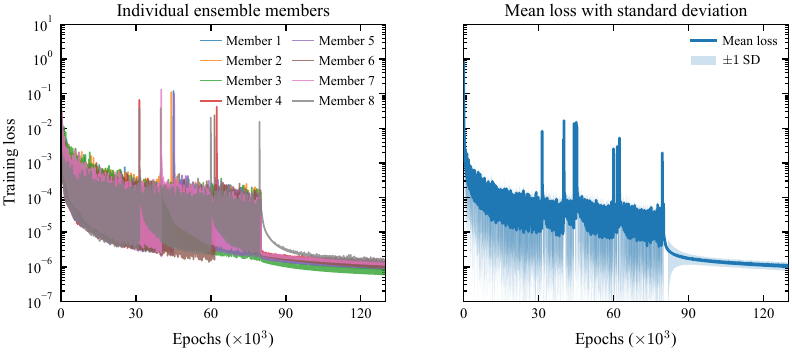}
	\caption{Training loss curves of Prob-DeepONet: individual ensemble members (left) and the ensemble mean with a $\pm 1$ standard deviation band (right).}
	\label{fig:dr-trainloss-onestep}
\end{figure}

\begin{figure}[!htbp]
	\centering
	
	\captionsetup[subfigure]{
		skip=2pt,
		justification=centering,
		singlelinecheck=true
	}
	
	\begin{subfigure}[t]{0.48\linewidth}
		\centering
		\begin{tikzpicture}
			\node[
			draw=black,
			dashed,
			line width=0.8pt,
			inner sep=3pt
			] {
				\includegraphics[width=\linewidth]{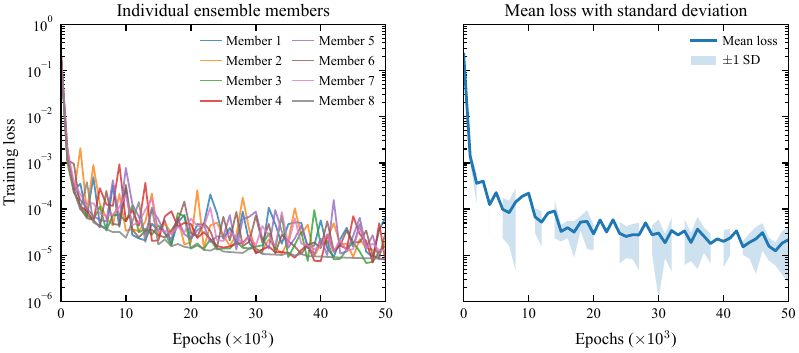}
			};
		\end{tikzpicture}
		\caption{Trunk network}
		\label{fig:dr-trainloss-twostep-trunk}
	\end{subfigure}
	\hfill
	\begin{subfigure}[t]{0.48\linewidth}
		\centering
		\begin{tikzpicture}
			\node[
			draw=black,
			dashed,
			line width=0.8pt,
			inner sep=3pt
			] {
				\includegraphics[width=\linewidth]{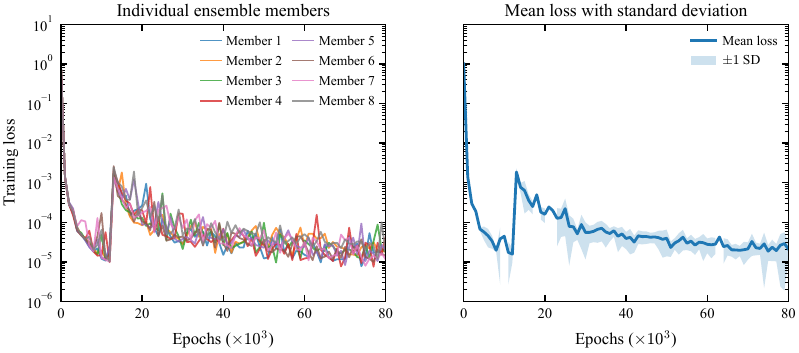}
			};
		\end{tikzpicture}
		\caption{Branch network}
		\label{fig:dr-trainloss-twostep-branch}
	\end{subfigure}
	\caption{Training loss curves of two-step MV-DeepONet for the trunk (a) and branch (b) networks: individual ensemble members (left) and the ensemble mean with a $\pm 1$ standard deviation band (right).}
	\label{fig:dr-trainloss-twostep}
\end{figure}

We first evaluate the basic fitting capability and generalization stability of the two methods within the training distribution. Both the training set and the in-distribution test set were generated with the correlation length $\ell=0.2$. The training loss curves of the two methods are shown in
Figs.~\ref{fig:dr-trainloss-onestep} and \ref{fig:dr-trainloss-twostep}. Both methods exhibit stable convergence, ensuring that the subsequent accuracy comparison is conducted under adequate training. For probabilistic branch training, a mean-only warm-up with fixed variance is adopted before jointly optimizing the mean and variance~\cite{MVE_Sluijterman2024}; accordingly, the branch loss curve exhibits a distinct transition at the end of the warm-up
stage and subsequently converges steadily under joint NLL optimization. The same training strategy is adopted in all subsequent examples. 

\vspace{-0.4em}

The quantitative errors under the same input distribution are summarized in Table~\ref{tab:dr-l2-train-interp}. Prob-DeepONet achieves lower training errors, whereas two-step MV-DeepONet yields lower mean and maximum in-distribution test errors. The relative degradation in predictive accuracy from the training set to unseen in-distribution samples is quantified by the ratio of the mean in-distribution test error to the mean training error. This ratio is \(0.0129/0.0010=12.9\) for Prob-DeepONet and \(0.0060/0.0037\approx1.62\) for two-step MV-DeepONet. These results indicate that the fixed orthogonal basis and subsequent coefficient regression provide stronger structural constraints and more stable in-distribution generalization. The out-of-distribution (OOD) results below further examine this behavior.
\begin{table}[H]
	\centering
	\caption{Mean and maximum relative $L_2$ errors on the training set and the in-distribution test set.}
	\label{tab:dr-l2-train-interp}
	\begin{tabular}{lcc}
		\toprule
		Metric & Prob-DeepONet & two-step MV-DeepONet \\
		\midrule
		Mean relative \(L_2\) error (training)
		& \textbf{0.0010} & 0.0037 \\
		Maximum relative \(L_2\) error (training)
		& \textbf{0.0046} & 0.0099 \\
		Mean relative \(L_2\) error (in-distribution)
		& 0.0129 & \textbf{0.0060} \\
		Maximum relative \(L_2\) error (in-distribution)
		& 0.1025 & \textbf{0.0391} \\
		\bottomrule
	\end{tabular}
\end{table}

\FloatBarrier

\vspace{4pt}
\noindent\textbf{Out-of-Distribution Generalization Performance}
\label{subsec:dr-extrapolation}

Table~\ref{tab:dr-extrapolation-l2} compares the OOD performance over input correlation lengths \(\ell\in[0.001,1.0]\), with \(\ell=0.2\) used for training. Both methods maintain low errors for moderate and large correlation lengths (\(\ell\geq0.18\)), while their difference becomes evident for \(\ell\leq0.15\). Under the most severe OOD shifts, Prob-DeepONet yields errors of \(44.15\%\) at \(\ell=0.005\) and \(109.64\%\) at \(\ell=0.001\). The corresponding errors of two-step MV-DeepONet are \(12.26\%\) and \(33.09\%\), respectively. Across all correlation lengths, two-step MV-DeepONet reduces the mean error from \(13.43\%\) to \(5.82\%\), corresponding to a reduction of approximately \(56.7\%\).

\begin{table}[!htbp]
	\centering
	\caption{Relative $L_2$ errors of Prob-DeepONet and two-step MV-DeepONet for out-of-distribution samples at different input correlation lengths.}
	\label{tab:dr-extrapolation-l2}
	\begin{tabular}{lcc}
		\toprule
		Correlation length \(\ell\) & Prob-DeepONet & two-step MV-DeepONet \\
		\midrule
		1.0   & 0.488\%   & \textbf{0.234\%} \\
		0.8   & 0.461\%   & \textbf{0.246\%} \\
		0.6   & 0.413\%   & \textbf{0.282\%} \\
		0.5   & 0.435\%   & \textbf{0.319\%} \\
		0.4   & \textbf{0.318\%} & 0.343\% \\
		0.3   & 0.368\%   & \textbf{0.309\%} \\
		0.25  & 0.542\%   & \textbf{0.300\%} \\
		0.22  & 0.513\%   & \textbf{0.290\%} \\
		0.18  & 0.708\%   & \textbf{0.273\%} \\
		0.15  & 2.932\%   & \textbf{0.530\%} \\
		0.1   & 7.391\%   & \textbf{4.916\%} \\
		0.08  & \textbf{6.401\%} & 6.462\% \\
		0.05  & 9.882\%   & \textbf{8.762\%} \\
		0.03  & 20.576\%  & \textbf{10.734\%} \\
		0.01  & 19.245\%  & \textbf{12.243\%} \\
		0.005 & 44.152\%  & \textbf{12.260\%} \\
		0.002 & 17.289\%  & \textbf{13.125\%} \\
		0.001 & 109.639\% & \textbf{33.094\%} \\
		\midrule
		Mean error & 13.431\% & \textbf{5.818\%} \\
		\bottomrule
	\end{tabular}
\end{table}

\begin{figure}[!htbp]
	\centering
	\captionsetup[subfigure]{
		skip=2pt,
		justification=centering,
		singlelinecheck=true
	}
	
	\begin{subfigure}[t]{0.48\textwidth}
		\centering
		\includegraphics[width=\linewidth]{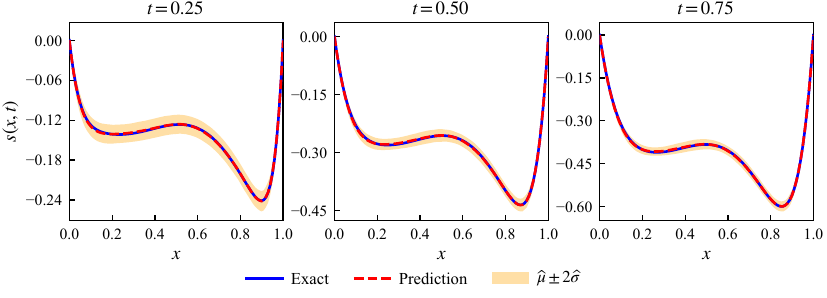}
		\caption{Prob-DeepONet, \(\ell=0.5\).}
		\label{fig:dr-extrap-onestep-L05}
	\end{subfigure}
	\hfill
	\begin{subfigure}[t]{0.48\textwidth}
		\centering
		\includegraphics[width=\linewidth]{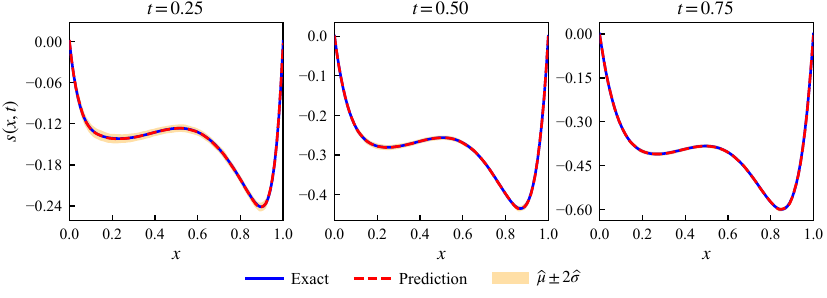}
		\caption{two-step MV-DeepONet, \(\ell=0.5\).}
		\label{fig:dr-extrap-twostep-L05}
	\end{subfigure}
	
	\par\vspace{0.3em}
	
	\begin{subfigure}[t]{0.48\textwidth}
		\centering
		\includegraphics[width=\linewidth]{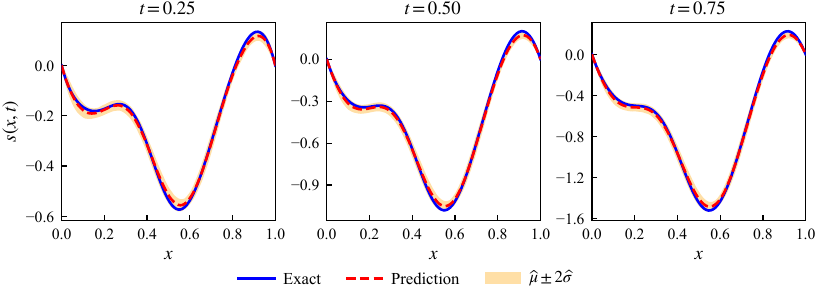}
		\caption{Prob-DeepONet, \(\ell=0.15\).}
		\label{fig:dr-extrap-onestep-L015}
	\end{subfigure}
	\hfill
	\begin{subfigure}[t]{0.48\textwidth}
		\centering
		\includegraphics[width=\linewidth]{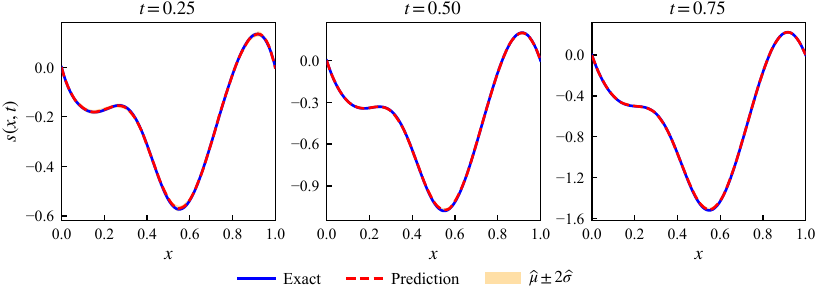}
		\caption{two-step MV-DeepONet, \(\ell=0.15\).}
		\label{fig:dr-extrap-twostep-L015}
	\end{subfigure}
	
	\par\vspace{0.3em}
	
	\begin{subfigure}[t]{0.48\textwidth}
		\centering
		\includegraphics[width=\linewidth]{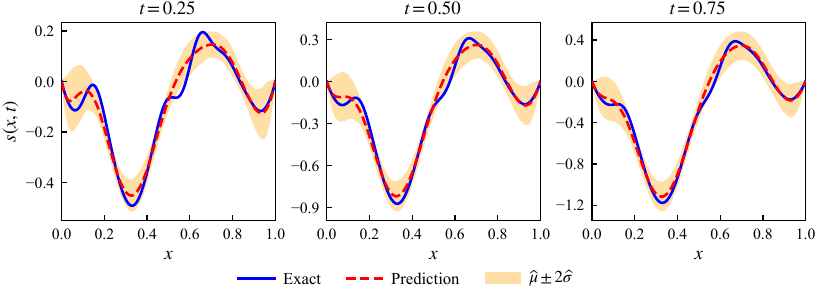}
		\caption{Prob-DeepONet, \(\ell=0.05\).}
		\label{fig:dr-extrap-onestep-L005}
	\end{subfigure}
	\hfill
	\begin{subfigure}[t]{0.48\textwidth}
		\centering
		\includegraphics[width=\linewidth]{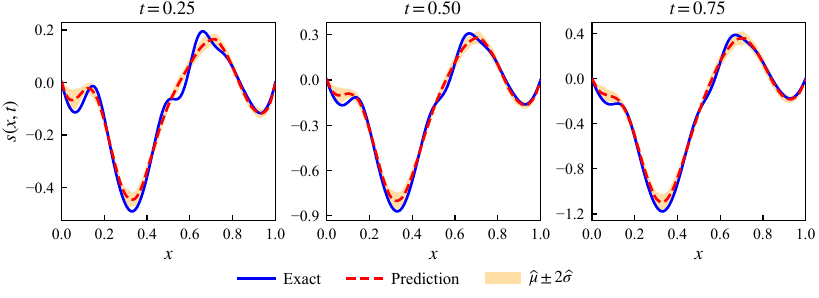}
		\caption{two-step MV-DeepONet, \(\ell=0.05\).}
		\label{fig:dr-extrap-twostep-L005}
	\end{subfigure}
	
	\par\vspace{0.3em}
	
	\begin{subfigure}[t]{0.48\textwidth}
		\centering
		\includegraphics[width=\linewidth]{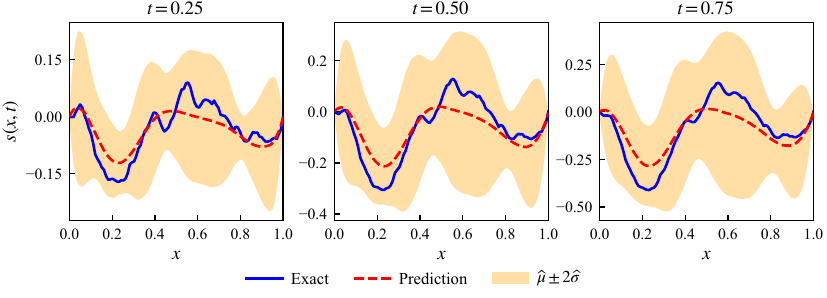}
		\caption{Prob-DeepONet, \(\ell=0.005\).}
		\label{fig:dr-extrap-onestep-L0005}
	\end{subfigure}
	\hfill
	\begin{subfigure}[t]{0.48\textwidth}
		\centering
		\includegraphics[width=\linewidth]{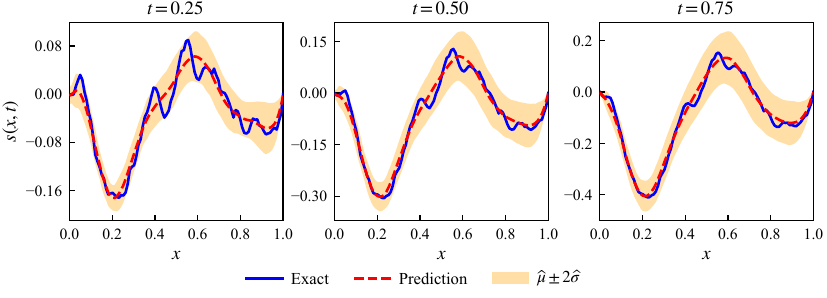}
		\caption{two-step MV-DeepONet, \(\ell=0.005\).}
		\label{fig:dr-extrap-twostep-L0005}
	\end{subfigure}
	
	\caption{Predictions and predictive uncertainty bands of Prob-DeepONet (left column) and two-step MV-DeepONet (right column) for representative out-of-distribution samples at different input correlation lengths.}
	\label{fig:dr-extrap-comparison}
\end{figure}

To illustrate the predictive behavior under varying degrees of distribution shift, Fig.~\ref{fig:dr-extrap-comparison} presents representative samples at four correlation lengths: $\ell=0.5$ corresponds to a smooth-input regime with a large correlation length; $\ell=0.15$ represents a mild distribution shift relatively close to the training distribution; $\ell=0.05$ represents a moderate OOD condition; and $\ell=0.005$ represents an extreme OOD condition.

Both methods remain accurate under mild distribution shifts. As \(\ell\) decreases, the input source becomes increasingly oscillatory and induces more complex local variations in the solution. Under severe shifts, Prob-DeepONet exhibits pronounced prediction deviations, whereas two-step MV-DeepONet retains the principal solution structure with lower prediction errors.

The contrasting OOD behavior of the two methods is related to differences in both their learned output representations and their formulations of the input-to-coefficient mapping. Prob-DeepONet jointly optimizes the trunk basis and input-dependent coefficients, allowing them to compensate for each other
on the training samples. This coupling facilitates training-data fitting but may produce a less stable representation when the input distribution changes. 
In two-step MV-DeepONet, a shared output basis is first learned from the output snapshots, orthogonalized, and then fixed, after which the branch network learns the corresponding modal coefficients. This decoupling provides a better-conditioned and more stable coordinate system for coefficient regression while reducing the co-adaptation between the basis functions and coefficients.

The role of the orthogonal basis can be expressed through the output error decomposition
\begin{align}
	\hat{\mathbf{s}}-\mathbf{s}
	&=
	\widetilde{\mathbf{Q}}^{*}\hat{\mathbf{c}}
	-
	\mathbf{P}\mathbf{s}
	+
	\mathbf{P}\mathbf{s}
	-
	\mathbf{s}
	\notag\\
	&=
	\widetilde{\mathbf{Q}}^{*}
	\left(
	\hat{\mathbf{c}}
	-
	\mathbf{c}_{\widetilde Q}
	\right)
	-
	\left(
	\mathbf{I}-\mathbf{P}
	\right)\mathbf{s}.
	\label{eq:ood-output-error-decomposition}
\end{align}
Here, \(\hat{\mathbf{s}}
=
\widetilde{\mathbf{Q}}^{*}\hat{\mathbf{c}}\) is the reconstructed
prediction, \(\hat{\mathbf{c}}\) is the coefficient vector predicted
by the branch network, and
\(\mathbf{c}_{\widetilde Q}
=
(\widetilde{\mathbf{Q}}^{*})^{\top}\mathbf{s}\) denotes the projection
coefficients of the true output. Moreover,
\(\mathbf{P}
=
\widetilde{\mathbf{Q}}^{*}
(\widetilde{\mathbf{Q}}^{*})^{\top}\) is the orthogonal projector onto the learned output subspace. The first term represents the coefficient regression error within the learned subspace, whereas the second is the component of the true output that cannot be represented by that subspace. Since these two
terms lie in mutually orthogonal subspaces,
\begin{equation}
	\left\|
	\hat{\mathbf{s}}-\mathbf{s}
	\right\|_{2}^{2}
	=
	\left\|
	\widetilde{\mathbf{Q}}^{*}
	\left(
	\hat{\mathbf{c}}
	-
	\mathbf{c}_{\widetilde Q}
	\right)
	\right\|_{2}^{2}
	+
	\left\|
	\left(
	\mathbf{I}-\mathbf{P}
	\right)\mathbf{s}
	\right\|_{2}^{2}.
	\label{eq:ood-output-error-pythagorean}
\end{equation}
Moreover, the orthonormality of
\(\widetilde{\mathbf{Q}}^{*}\) gives
\begin{equation}
	\left\|
	\widetilde{\mathbf{Q}}^{*}
	\left(
	\hat{\mathbf{c}}
	-
	\mathbf{c}_{\widetilde Q}
	\right)
	\right\|_{2}
	=
	\left\|
	\hat{\mathbf{c}}
	-
	\mathbf{c}_{\widetilde Q}
	\right\|_{2},
	\label{eq:ood-coefficient-error-isometry}
\end{equation}
so coefficient prediction errors are not amplified during output reconstruction.

For a general nonorthogonal trunk matrix \(\boldsymbol{\Phi}\), as used in the jointly trained representation of Prob-DeepONet, only
\begin{equation}
	\left\|
	\boldsymbol{\Phi}
	\left(
	\hat{\mathbf{c}}
	-
	\mathbf{c}_{\Phi}
	\right)
	\right\|_{2}
	\leq
	\left\|
	\boldsymbol{\Phi}
	\right\|_{2}
	\left\|
	\hat{\mathbf{c}}
	-
	\mathbf{c}_{\Phi}
	\right\|_{2}
	\label{eq:ood-nonorthogonal-error-bound}
\end{equation}
is guaranteed. Strongly correlated basis functions may therefore lead to an ill-conditioned coordinate representation and amplify coefficient errors during reconstruction.

Fixing the shared basis also introduces a structural regularization effect, encouraging the model to preserve dominant solution structures shared across the training samples rather than overadapting the output representation to specific local details. Consequently, when OOD inputs introduce previously unseen high-frequency components or local structures, two-step MV-DeepONet may not fully reconstruct these new details, but it can generally retain the principal structure of the solution field and exhibit more controlled error growth. By contrast, the training-specific compensation learned by Prob-DeepONet through
joint optimization may become ineffective under severe distribution shifts, leading to larger prediction deviations.

A similar pattern is observed in the uncertainty bands. As shown in
Fig.~\ref{fig:dr-extrap-comparison}, the bands of Prob-DeepONet
become increasingly broad as \(\ell\) decreases and expand over much
of the response profile under severe distribution shifts. In contrast,
two-step MV-DeepONet generally produces tighter bands whose widths vary
more coherently with the local solution structure.

This behavior is consistent with the covariance parameterizations of the two methods. Prob-DeepONet estimates the marginal variance separately at each output point under a pointwise diagonal covariance assumption and therefore does not explicitly couple uncertainty across different points in the output domain. Local prediction discrepancies may consequently be accommodated through pointwise variance inflation, producing broad uncertainty bands with limited organization across the output domain. In contrast, two-step MV-DeepONet predicts the variances of the modal coefficients and maps them to the
output domain through the shared basis functions. Each modal coefficient variance contributes to the marginal uncertainty at multiple output points through the shared basis. This shared contribution induces coordinated variations in marginal uncertainty and off-diagonal covariance across the output domain.

\FloatBarrier

\subsubsection{Analysis of Output Covariance Recovery Accuracy}
\label{sec:dr-cov-recovery}

Output covariance recovery is evaluated using the normalized Frobenius error defined in Eq.~\eqref{eq:metric-cov-frobenius}. Multi-anchor correlation plots are further used to examine local correlation structures across the output domain. Since Prob-DeepONet assumes a pointwise diagonal covariance and does not explicitly represent off-diagonal dependence, this analysis focuses on two-step MV-DeepONet.

\vspace{4pt}
\noindent\textbf{Covariance Error Decomposition and Quantitative Accuracy Analysis}

\begin{figure}[!htbp]
	\centering
	\includegraphics[width=0.5\textwidth]
	{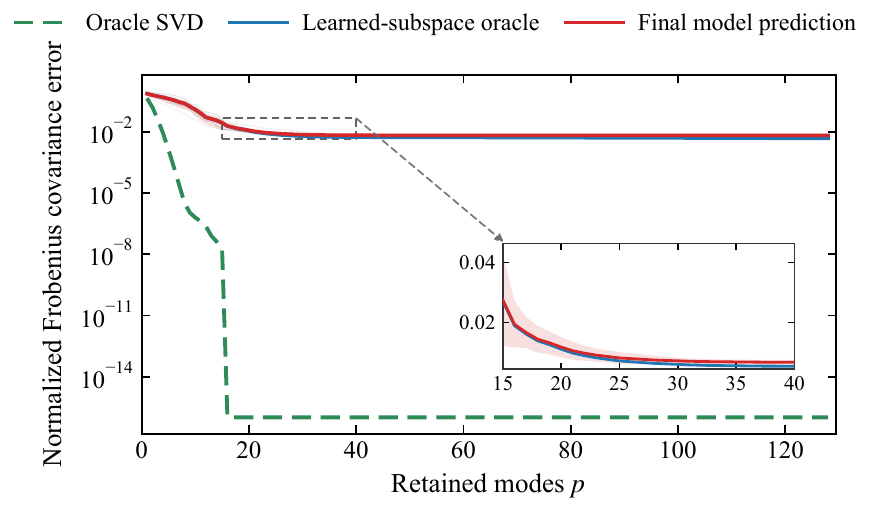}
	\caption{Normalized Frobenius covariance error as a function of the
		number of retained modes.}
	\label{fig:dr-covariance-error-decomposition}
\end{figure}

Following the decomposition in Eq.~\eqref{eq:app-cov-basic-bound}, the Oracle SVD curve characterizes the intrinsic rank-\(p\) truncation error. The separation between the Oracle SVD and Learned-subspace oracle curves primarily reflects the trunk-subspace learning error, whereas the separation between the Learned-subspace oracle and Final model prediction curves primarily reflects the branch coefficient-covariance prediction error. Fig.~\ref{fig:dr-covariance-error-decomposition} presents the three error curves, and Table~\ref{tab:dr-covariance-thresholds} reports the minimum number of modes required to reach each target normalized Frobenius error.

\begin{table}[!htbp]
	\centering
	\caption{Minimum numbers of retained modes required to achieve different target normalized Frobenius covariance-error levels for the three covariance-recovery curves.}
	\label{tab:dr-covariance-thresholds}
	\begin{tabular}{
			>{\centering\arraybackslash}p{0.24\textwidth}
			>{\centering\arraybackslash}p{0.20\textwidth}
			>{\centering\arraybackslash}p{0.25\textwidth}
			>{\centering\arraybackslash}p{0.25\textwidth}
		}
		\toprule
		Target error
		& Oracle SVD
		& Learned-subspace oracle
		& Final model prediction \\
		\midrule
		Error \(\leq 20.00\%\) & \(p\geq 2\) & \(p\geq 9\)  & \(p\geq 9\)  \\
		Error \(\leq 10.00\%\) & \(p\geq 3\) & \(p\geq 11\) & \(p\geq 11\) \\
		Error \(\leq 5.00\%\)  & \(p\geq 3\) & \(p\geq 13\) & \(p\geq 13\) \\
		Error \(\leq 2.00\%\)  & \(p\geq 4\) & \(p\geq 16\) & \(p\geq 16\) \\
		Error \(\leq 1.00\%\)  & \(p\geq 4\) & \(p\geq 21\) & \(p\geq 22\) \\
		\bottomrule
	\end{tabular}
\end{table}

\textbf{(1) Oracle SVD curve.}
The Oracle SVD error decreases rapidly and reaches approximately \(10^{-13}\) for \(p\geq15\), indicating that the low-rank compressibility assumption introduced in Section~\ref{sec:cov-error-decomposition} is naturally satisfied in the reaction--diffusion problem, and the abrupt decrease near $p\approx15$ also supports the nondegenerate spectral-gap assumption $\gamma_{p}>0$. Once the effective numerical rank is reached, \(T_{1}\) becomes negligible and is no longer the principal limitation on covariance recovery.

\textbf{(2) Comparison between the Learned-subspace oracle curve and the Final model prediction curve.}
The two curves are nearly indistinguishable over the full range of \(p\). For \(p\in[15,40]\), their separation remains approximately \(1\times10^{-3}\), while both errors decrease from approximately \(2\times10^{-1}\) to \(6\times10^{-3}\). These results support the following observations:

\begin{enumerate}[leftmargin=2.2em,itemsep=0.4em]
	\item
	The contribution of \(T_{3}\) is substantially smaller than that of \(T_{2}\). This indicates that the branch approximation and optimization errors have been controlled to a relatively low level. Branch coefficient-covariance prediction is therefore not the principal limitation on covariance recovery.
	
	\item
	The dominant contribution is \(T_{2}\), which is associated with the difference between the learned projector
	\(\mathbf P=\widetilde{\mathbf Q}^{*}
	(\widetilde{\mathbf Q}^{*})^{\top}\)
	and the true principal-subspace projector
	\(\mathbf P^{\star}\). As discussed in Section~\ref{sec:cov-error-decomposition}, this difference contains two main components. The first is the excess trunk reconstruction risk \(\Delta_T\) relative to the empirically optimal projector \(\widehat{\mathbf P}^{\star}\). The second arises from finite-sample covariance estimation, characterized by \(\left\|\widehat{\boldsymbol{\Sigma}}-\boldsymbol{\Sigma}\right\|_{F} =\mathcal{O}_{\mathbb{P}}\left(\sqrt{M}K^{-1/2}\right)\), where \(M\) is the output dimension and \(K\) is the number of training samples. Its influence on the estimated subspace is further controlled by the corresponding spectral gap. For this problem, \(M=10^{4}\) and \(K=10^{3}\). The resulting relative statistical-error scale is
	\begin{equation}
		\frac{
			\sqrt{M}K^{-1/2}
		}{
			\left\|
			\widehat{\boldsymbol{\Sigma}}
			\right\|_{F}
		}
		\approx
		2.0\times10^{-3}.
		\label{eq:dr-statistical-error-scale}
	\end{equation}
	Because \(\mathcal{O}_{\mathbb{P}}(\cdot)\) specifies only an asymptotic order, this value should be interpreted as an approximate scale rather than an exact error. It is about one third of the observed \(T_{2}\approx6\times10^{-3}\), suggesting that \(\Delta_T\) is the larger contribution within \(T_{2}\).
	
	This conclusion concerns only the relative importance of the two contributions within $T_{2}$ and does not imply inadequate trunk training. The trunk loss in Fig.~\ref{fig:dr-trainloss-twostep} stabilizes at a low level after approximately \(50{,}000\) epochs. Moreover, $T_{2}\approx6\times10^{-3}$ is of the same order as the $0.60\%$ mean relative $L_{2}$ error of two-step MV-DeepONet on the in-distribution test set (Table~\ref{tab:dr-l2-train-interp}). In principle, $\Delta_{T}$ could be further reduced by increasing the trunk-network capacity, extending the training duration, or adopting more refined optimization strategies. However, because the current covariance recovery error is already of the same order as the mean prediction error, further hyperparameter refinement would provide limited improvement while substantially increasing the training cost. The objective here is to determine whether the covariance matrix can be stably recovered to a practically useful level of accuracy, rather than to pursue the lowest attainable numerical error. 
\end{enumerate}

In summary, covariance recovery for the reaction--diffusion problem shows the clear hierarchy
\(
T_{1}\ll T_{2},
T_{3}\ll T_{2}.
\)
Within $T_{2}$, the excess trunk reconstruction risk $\Delta_{T}$ is the principal source, whereas the finite-sample statistical contribution is comparatively small. These results are consistent with the theoretical decomposition in Section~\ref{sec:cov-error-decomposition} and demonstrate that two-step MV-DeepONet recovers the output covariance with an error comparable in magnitude to its mean prediction error.

\vspace{4pt}
\noindent\textbf{Recovery of Local Correlation Structures:
	Multi-Anchor Correlation Map Assessment}

\begin{figure}[!htbp]
	\centering
	\includegraphics[width=0.88\textwidth]
	{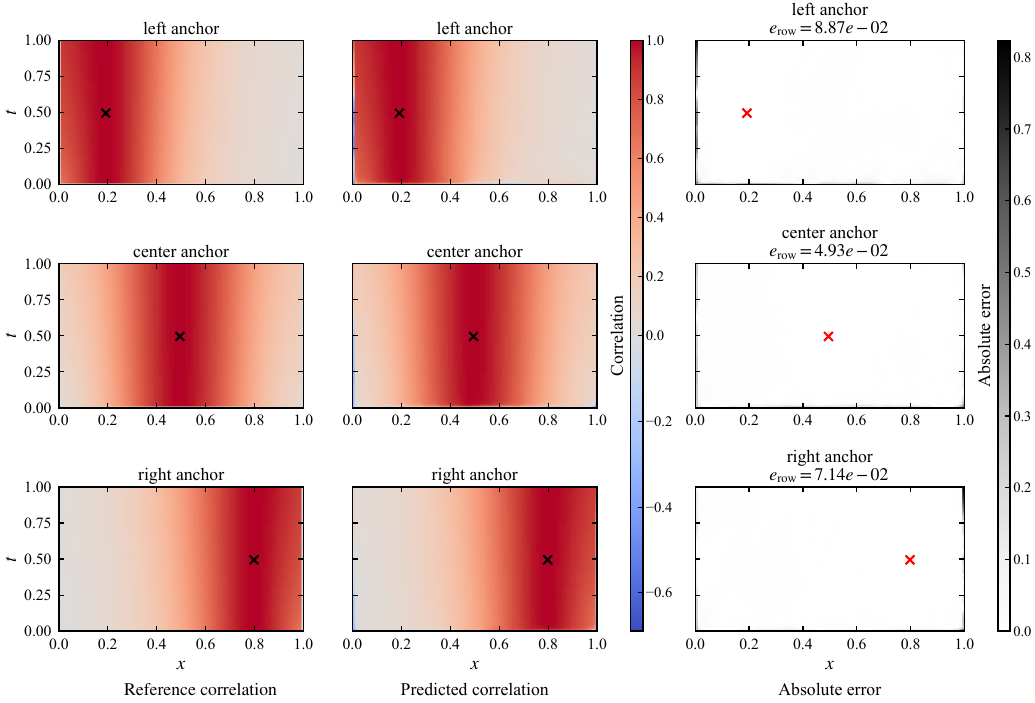}
	\caption{Multi-anchor correlation maps showing the reference correlations, predicted correlations, and corresponding absolute errors for three anchor points, where $e_{\mathrm{row}}$ denotes the relative $L_2$ error of the predicted correlation map for each anchor point.}
	\label{fig:dr-multi-anchor-correlation}
\end{figure}

The normalized Frobenius error measures the overall covariance discrepancy but does not reveal whether local correlations across the output domain are accurately recovered. Multi-anchor correlation maps are therefore used for a more detailed structural assessment.
Three anchor points are selected from the left, central, and right regions of the output domain \((x,t)\in[0,1]^2\). For each anchor, the reference and predicted correlation maps are compared together with their absolute error. The resulting maps are shown in Fig.~\ref{fig:dr-multi-anchor-correlation}.

\vspace{4pt}
\noindent\textbf{Overall pattern consistency.}
For all three anchors, the predicted maps reproduce the principal vertical correlation bands, the smooth decay away from the anchor location, and the weak negative correlations in the far field. The mean row-wise correlation error is approximately \(6.97\times10^{-2}\), indicating a relatively small structural
discrepancy. These results show that two-step MV-DeepONet recovers not only the global covariance magnitude but also its local correlation patterns.

\vspace{4pt}
\noindent\textbf{Physical interpretation of the correlation structure.}
The vertical bands arise primarily from the time-independent random source \(u(x)\), which continuously drives the system throughout its evolution. At a fixed spatial location \(x=x_a\), responses at different times are influenced by similar local components of the same random source and therefore remain strongly correlated. This shared influence weakens as the spatial distance from \(x_a\) increases, leading to correlation decay away from the anchor. The combination of persistent temporal correlation and spatial decay produces the observed vertical-band structure.

The recovery of this physical structure highlights a core advantage of two-step MV-DeepONet. The trunk stage first learns the shared output-space basis $\widetilde{\mathbf Q}^{*}$, which captures the dominant structures of the output field. The branch stage then models the input-dependent variances of the modal coefficients in the low-dimensional coefficient space. Mapping these coefficient uncertainties back through the shared basis induces an output covariance that naturally inherits the spatial patterns encoded by the columns of $\widetilde{\mathbf Q}^{*}$.

\subsubsection{Prediction Interval Calibration and Coverage Analysis}
\label{sec:dr-calibration}

Accurate covariance recovery does not necessarily ensure nominal prediction interval coverage. For engineering uncertainty quantification, prediction intervals should provide reliable coverage while remaining sufficiently narrow. We therefore apply CP as a post-hoc calibration procedure and evaluate its effectiveness using the PICP and MPIW defined in Eq.~\eqref{eq:metric-picp-mpiw}.

\begin{table}[!htbp]
	\centering
	\caption{PICP and MPIW of Prob-DeepONet and two-step MV-DeepONet before and after CP calibration.}
	\label{tab:dr-calibration}
	
	\renewcommand{\arraystretch}{1.25}
	\setlength{\tabcolsep}{6pt}
	\footnotesize
	
	\sisetup{
		detect-all,
		table-number-alignment = center
	}
	
	\begin{tabular}{
			l
			c
			S[table-format=3.2]
			S[table-format=1.2e-2]
			S[table-format=3.2]
			S[table-format=1.2e-2]
		}
		\toprule
		\multirow{2}{*}{Method}
		& \multirow{2}{*}{\makecell{Nominal\\coverage (\%)}}
		& \multicolumn{2}{c}{Before CP calibration}
		& \multicolumn{2}{c}{After CP calibration} \\
		
		\cmidrule(lr){3-4}
		\cmidrule(lr){5-6}
		
		&
		& {PICP (\%)}
		& {MPIW}
		& {PICP (\%)}
		& {MPIW} \\
		\midrule
		
		\multirow{3}{*}{Prob-DeepONet}
		& 90
		& 99.96
		& 4.62e-2
		& 90.52
		& 1.60e-2 \\
		
		& 95
		& 99.99
		& 5.51e-2
		& 95.30
		& 2.01e-2 \\
		
		& 99
		& 100.00
		& 7.24e-2
		& 98.97
		& 2.88e-2 \\
		
		\midrule
		
		\multirow{3}{*}{\makecell[l]{two-step\\MV-DeepONet}}
		& 90
		& 98.99
		& 1.60e-2
		& 89.23
		& 5.10e-3 \\
		
		& 95
		& 99.41
		& 1.90e-2
		& 94.21
		& 6.95e-3 \\
		
		& 99
		& 99.81
		& 2.50e-2
		& 98.44
		& 1.34e-2 \\
		
		\bottomrule
	\end{tabular}
\end{table}

Before calibration, both methods exhibit pronounced overcoverage, with PICP values close to \(100\%\) at all nominal levels. Their original
intervals are therefore overly conservative. Prob-DeepONet also produces substantially wider intervals. At the \(95\%\) nominal level, its MPIW is approximately \(2.9\) times that of two-step MV-DeepONet. The latter thus provides more compact uncertainty intervals even before calibration.
After CP calibration, the PICP values move closer to their nominal targets, while the MPIW values decrease substantially. This confirms that CP effectively reduces excessive coverage and improves interval sharpness. The calibrated PICP values of Prob-DeepONet are slightly closer to the nominal levels, but its intervals remain considerably wider. At the \(95\%\) level, its calibrated MPIW is again approximately \(2.9\) times that of two-step MV-DeepONet.

Overall, both methods achieve coverage close to the prescribed levels after calibration. Two-step MV-DeepONet maintains markedly narrower intervals and thus
achieves a better balance between coverage and interval sharpness in the reaction--diffusion problem.

\subsection{Burgers Equation}
\label{sec:burgers-introduction}

The Burgers equation is a canonical nonlinear convection--diffusion model for studying nonlinear transport, shock formation, and simplified turbulence mechanisms~\cite{LinearNonlinearWaves_Whitham1974}. Its quadratic convective term is analogous to that in the Navier--Stokes equations, making it a widely used benchmark in fluid dynamics, stochastic wave propagation, and numerical method validation~\cite{BurgersTurbulence_Bec2007}.
We consider the one-dimensional viscous Burgers equation on \([0,1]\):
\begin{equation}
	\frac{\partial s}{\partial t}
	+
	s\frac{\partial s}{\partial x}
	=
	\nu\frac{\partial^2 s}{\partial x^2},
	\qquad
	(x,t)\in[0,1]\times[0,1],
	\label{eq:burgers-equation}
\end{equation}
subject to the periodic boundary conditions and the initial condition
\begin{equation}
	s(0,t)=s(1,t),
	\qquad
	\frac{\partial s}{\partial x}(0,t)
	=
	\frac{\partial s}{\partial x}(1,t),
	\qquad
	s(x,0)=u(x).
\end{equation}
The viscosity coefficient is set to \(\nu=0.01\). Here, the random initial condition \(u(x)\) is the input function, and the corresponding solution \(s(x,t)\) is the output field. The target operator is
\begin{equation}
	\mathcal{G}:u(x)\mapsto s(x,t).
\end{equation}

\vspace{4pt}
\noindent\textbf{Source of Uncertainty.}
The uncertainty arises from the random initial condition \(u(x)\).
Random initial velocity fields commonly occur in fluid
transport~\cite{VehicularTraffic_Buendia2008} and the early evolution of
turbulence~\cite{BurgersTurbulence_Bec2007}. Through nonlinear convection, these perturbations can alter waveform
translation, steepening, and shock location. The resulting changes in
the output distribution and correlation structure make this problem a
suitable benchmark for nonlinear uncertainty propagation and covariance
recovery.

\vspace{4pt}
\noindent\textbf{Data Generation.}
The initial condition \(u(x)\) is modeled as a Gaussian random field
with covariance operator
\begin{equation}
	25^2\left(-\Delta+\tau^2 I\right)^{-2},
	\label{eq:burgers-grf}
\end{equation}
where \(\tau\) controls the correlation scale and spectral composition of the random field. A larger \(\tau\) corresponds to a shorter characteristic length scale and richer high-frequency content. The training distribution is defined by \(\tau=25\). Consistent with the periodic boundary conditions, the random fields are generated using a Fourier spectral expansion and evaluated at \(m=101\) uniformly spaced spatial locations as inputs to the branch network. For each initial condition, Eq.~\eqref{eq:burgers-equation} is solved on \((x,t)\in[0,1]^2\) using a Fourier spectral method implemented in Chebfun. The solution is evaluated at \(101\) spatial and \(101\) temporal points, yielding a \(101\times101\) spatiotemporal output matrix. A total of \(2{,}000\) samples are generated, comprising \(1{,}000\) training samples, \(500\) test samples, and \(500\) calibration samples.

\FloatBarrier

\subsubsection{Prediction Accuracy and Generalization Performance}
\label{sec:burgers-accuracy-generalization}

\vspace{4pt}
\noindent\textbf{Training and In-Distribution Accuracy}
\label{subsec:burgers-train-interp}

\begin{figure}[!htbp]
	\centering
	\includegraphics[
	width=0.5\linewidth
	]{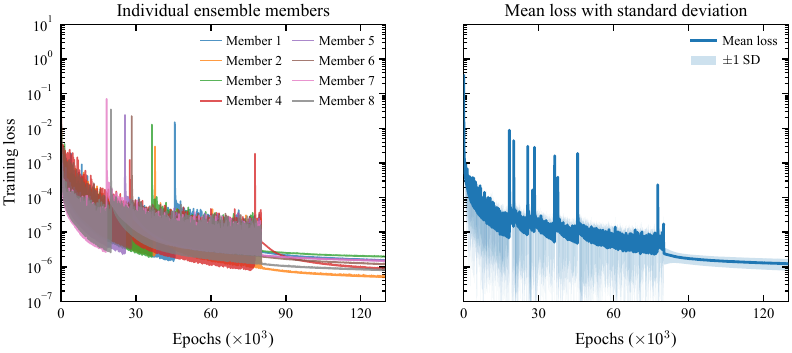}
	\caption{Training loss curves of Prob-DeepONet: individual ensemble members (left) and the ensemble mean with a $\pm 1$ standard deviation band (right).}
	\label{fig:burgers-trainloss-onestep}
\end{figure}

\begin{figure}[!htbp]
	\centering
	
	\captionsetup[subfigure]{
		skip=2pt,
		justification=centering,
		singlelinecheck=true
	}
	
	\begin{subfigure}[t]{0.48\linewidth}
		\centering
		\begin{tikzpicture}
			\node[
			draw=black,
			dashed,
			line width=0.8pt,
			inner sep=3pt
			] {
				\includegraphics[
				width=\linewidth
				]{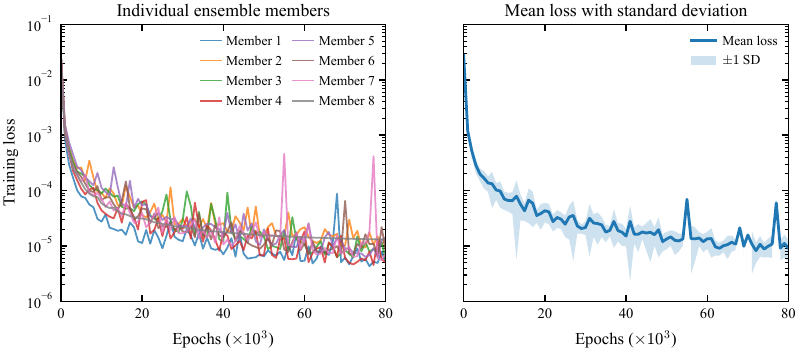}
			};
		\end{tikzpicture}
		\caption{Trunk network}
		\label{fig:burgers-trainloss-twostep-trunk}
	\end{subfigure}
	\hfill
	\begin{subfigure}[t]{0.48\linewidth}
		\centering
		\begin{tikzpicture}
			\node[
			draw=black,
			dashed,
			line width=0.8pt,
			inner sep=3pt
			] {
				\includegraphics[
				width=\linewidth
				]{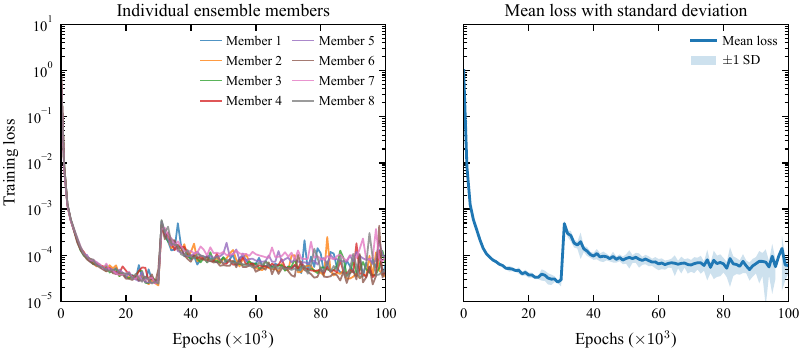}
			};
		\end{tikzpicture}
		\caption{Branch network}
		\label{fig:burgers-trainloss-twostep-branch}
	\end{subfigure}
	\caption{Training loss curves of two-step MV-DeepONet for the trunk (a) and branch (b) networks: individual ensemble members (left) and the ensemble mean with a $\pm 1$ standard deviation band (right).}
	\label{fig:burgers-trainloss-twostep}
\end{figure}

Figures~\ref{fig:burgers-trainloss-onestep} and~\ref{fig:burgers-trainloss-twostep} show stable convergence for both methods. Their training and in-distribution errors are reported in Table~\ref{tab:burgers-l2-train-interp}.
Prob-DeepONet achieves lower training errors, whereas two-step MV-DeepONet yields lower mean and maximum in-distribution errors. In particular, the maximum in-distribution error decreases from \(9.5806\%\) to \(5.6790\%\), indicating more stable performance on challenging in-distribution samples.

\begin{table}[htbp]
	\centering
	\caption{Mean and maximum relative $L_2$ errors on the training set and the in-distribution test set.}
	\label{tab:burgers-l2-train-interp}
	\begin{tabular}{lcc}
		\toprule
		Metric & Prob-DeepONet & two-step MV-DeepONet \\
		\midrule
		Mean relative \(L_2\) error (training)
		& \textbf{0.4440\%} & 0.6040\% \\
		Maximum relative \(L_2\) error (training)
		& \textbf{2.6840\%} & 3.8693\% \\
		Mean relative \(L_2\) error (in-distribution)
		& 1.6643\% & \textbf{1.5891\%} \\
		Maximum relative \(L_2\) error (in-distribution)
		& 9.5806\% & \textbf{5.6790\%} \\
		\bottomrule
	\end{tabular}
\end{table}

\FloatBarrier

\vspace{4pt}
\noindent\textbf{Out-of-Distribution Generalization Performance}
\label{subsec:burgers-extrapolation}

In this example, OOD generalization is evaluated by varying the random-field parameter \(\tau\) from its training value of \(25\). A larger \(\tau\) corresponds to a shorter correlation scale and richer high-frequency content. The resulting relative \(L_2\) errors are reported in Table~\ref{tab:burgers-extrapolation-l2}.

Prob-DeepONet performs better only at \(\tau=5\) and \(8\), whereas two-step MV-DeepONet achieves lower errors for all remaining cases. Under the extreme shift \(\tau=400\), their errors are \(293.932\%\) and \(68.931\%\), respectively. Across all conditions, two-step MV-DeepONet reduces the mean error from \(59.2977\%\) to \(27.3800\%\), corresponding to a reduction of approximately \(53.8\%\).
Fig.~\ref{fig:burgers-ood-bands-comparison} shows representative results at \(\tau=20\), \(64\), \(160\), and \(400\), spanning mild to extreme OOD conditions.

The qualitative differences observed for the challenging in-distribution samples become more pronounced under OOD shifts. Prob-DeepONet produces increasingly broad, nearly horizontal tubular uncertainty bands with limited adaptation to the local waveform, while its predictive mean deviates substantially under strong shifts. In contrast, two-step MV-DeepONet better preserves the overall solution profile and produces tighter, more structure-adaptive uncertainty bands that remain narrow in smooth regions and widen moderately near sharp transitions. Although both methods deteriorate under severe OOD conditions, the two-step model exhibits substantially greater predictive robustness.

\begin{figure}[!htbp]
	\centering
	
	\captionsetup[subfigure]{
		skip=2pt,
		justification=centering,
		singlelinecheck=true
	}
	
	\begin{subfigure}[t]{0.48\textwidth}
		\centering
		\includegraphics[
		width=\linewidth,
		height=0.145\textheight,
		keepaspectratio
		]{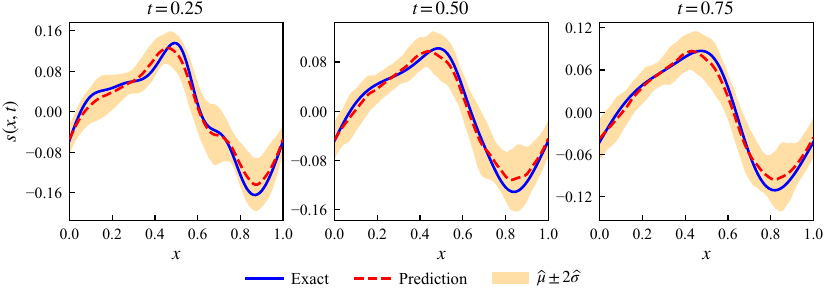}
		\caption{Prob-DeepONet, \(\tau=20\).}
		\label{fig:burgers-ood-band-onestep-tau20}
	\end{subfigure}
	\hfill
	\begin{subfigure}[t]{0.48\textwidth}
		\centering
		\includegraphics[
		width=\linewidth,
		height=0.145\textheight,
		keepaspectratio
		]{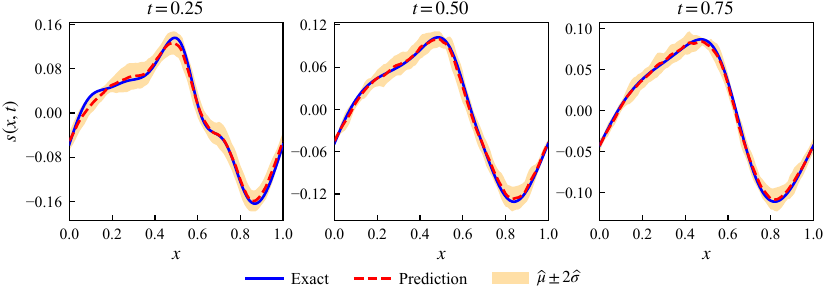}
		\caption{two-step MV-DeepONet, \(\tau=20\).}
		\label{fig:burgers-ood-band-twostep-tau20}
	\end{subfigure}
	
	\par\vspace{0.10em}
	
	\begin{subfigure}[t]{0.48\textwidth}
		\centering
		\includegraphics[
		width=\linewidth,
		height=0.145\textheight,
		keepaspectratio
		]{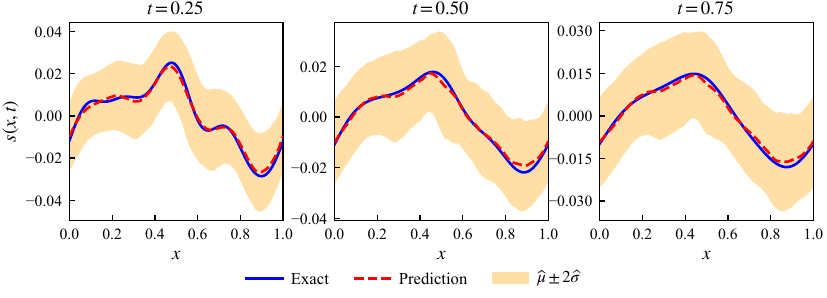}
		\caption{Prob-DeepONet, \(\tau=64\).}
		\label{fig:burgers-ood-band-onestep-tau64}
	\end{subfigure}
	\hfill
	\begin{subfigure}[t]{0.48\textwidth}
		\centering
		\includegraphics[
		width=\linewidth,
		height=0.145\textheight,
		keepaspectratio
		]{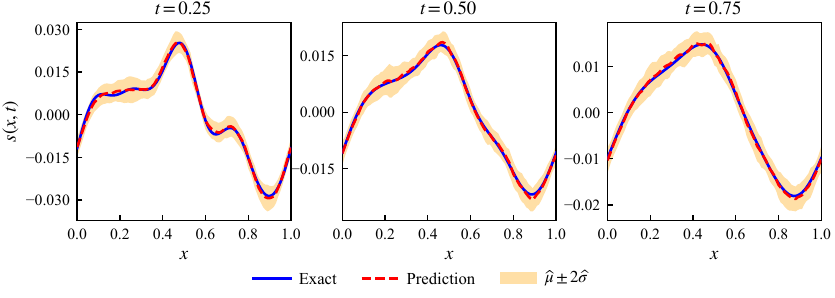}
		\caption{two-step MV-DeepONet, \(\tau=64\).}
		\label{fig:burgers-ood-band-twostep-tau64}
	\end{subfigure}
	
	\caption{Predictions and predictive uncertainty bands of Prob-DeepONet (left column) and two-step MV-DeepONet (right column) for representative out-of-distribution samples at different values of the random-field parameter \(\tau\).}
	\label{fig:burgers-ood-bands-comparison}
\end{figure}

\begin{figure}[!htbp]
	\ContinuedFloat
	\centering
	
	\captionsetup[subfigure]{
		skip=2pt,
		justification=centering,
		singlelinecheck=true
	}
	
	\setcounter{subfigure}{4}
	
	\begin{subfigure}[t]{0.48\textwidth}
		\centering
		\includegraphics[
		width=\linewidth,
		height=0.145\textheight,
		keepaspectratio
		]{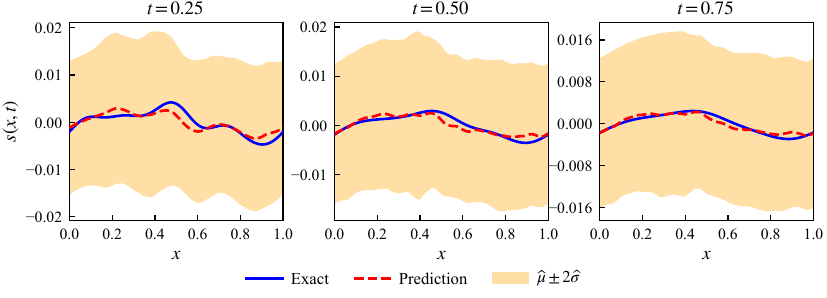}
		\caption{Prob-DeepONet, \(\tau=160\).}
		\label{fig:burgers-ood-band-onestep-tau160}
	\end{subfigure}
	\hfill
	\begin{subfigure}[t]{0.48\textwidth}
		\centering
		\includegraphics[
		width=\linewidth,
		height=0.145\textheight,
		keepaspectratio
		]{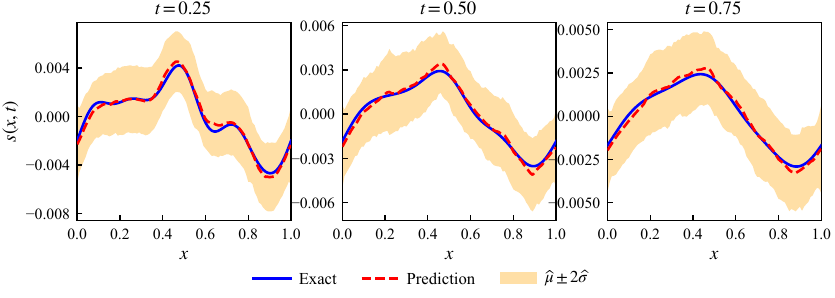}
		\caption{two-step MV-DeepONet, \(\tau=160\).}
		\label{fig:burgers-ood-band-twostep-tau160}
	\end{subfigure}
	
	\par\vspace{0.10em}
	
	\begin{subfigure}[t]{0.48\textwidth}
		\centering
		\includegraphics[
		width=\linewidth,
		height=0.145\textheight,
		keepaspectratio
		]{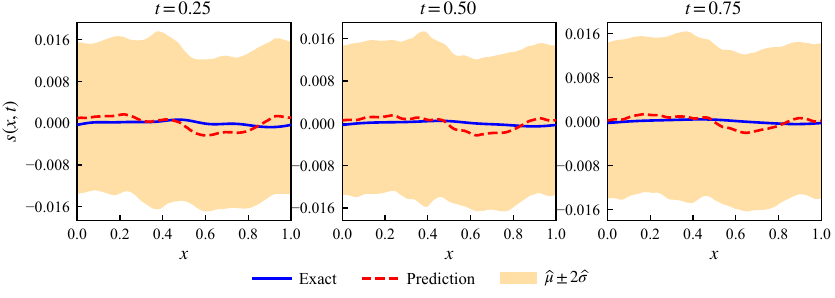}
		\caption{Prob-DeepONet, \(\tau=400\).}
		\label{fig:burgers-ood-band-onestep-tau400}
	\end{subfigure}
	\hfill
	\begin{subfigure}[t]{0.48\textwidth}
		\centering
		\includegraphics[
		width=\linewidth,
		height=0.145\textheight,
		keepaspectratio
		]{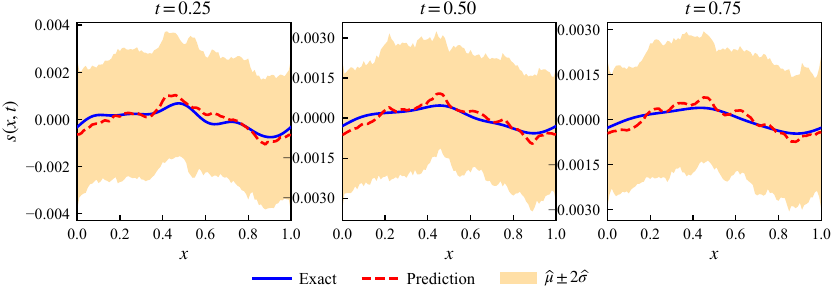}
		\caption{two-step MV-DeepONet, \(\tau=400\).}
		\label{fig:burgers-ood-band-twostep-tau400}
	\end{subfigure}
	
	\caption[]{Predictions and predictive uncertainty bands of Prob-DeepONet (left column) and two-step MV-DeepONet (right column) for representative out-of-distribution samples at different values of the random-field parameter \(\tau\) (continued).}
\end{figure}

\begin{table}[H]
	\centering
	\caption{Relative $L_2$ errors of Prob-DeepONet and two-step MV-DeepONet for out-of-distribution samples at different values of the random-field parameter $\tau$.}
	\label{tab:burgers-extrapolation-l2}
	\begin{tabular}{lcc}
		\toprule
		Random-field parameter \(\tau\)
		& Prob-DeepONet
		& two-step MV-DeepONet \\
		\midrule
		5   & \textbf{27.758\%}  & 37.269\% \\
		8   & \textbf{17.160\%}  & 23.467\% \\
		12  & 15.349\%  & \textbf{13.473\%} \\
		16  & 18.630\%  & \textbf{12.212\%} \\
		20  & 19.809\%  & \textbf{12.352\%} \\
		25  & 20.117\%  & \textbf{12.988\%} \\
		32  & 21.091\%  & \textbf{14.034\%} \\
		40  & 23.625\%  & \textbf{15.542\%} \\
		50  & 23.754\%  & \textbf{17.696\%} \\
		64  & 26.339\%  & \textbf{20.356\%} \\
		80  & 32.651\%  & \textbf{22.854\%} \\
		100 & 39.915\%  & \textbf{25.476\%} \\
		128 & 45.414\%  & \textbf{28.776\%} \\
		160 & 54.637\%  & \textbf{32.315\%} \\
		200 & 75.009\%  & \textbf{36.890\%} \\
		256 & 121.548\% & \textbf{44.173\%} \\
		320 & 190.620\% & \textbf{54.036\%} \\
		400 & 293.932\% & \textbf{68.931\%} \\
		\midrule
		Mean error
		& 59.2977\%
		& \textbf{27.3800\%} \\
		\bottomrule
	\end{tabular}
\end{table}

\FloatBarrier

\subsubsection{Analysis of Output Covariance Recovery Accuracy}
\label{subsec:burgers-cov-recovery}

This part further quantifies the output covariance recovery accuracy for the Burgers equation using the normalized Frobenius covariance error defined in Eq.~\eqref{eq:metric-cov-frobenius}. 

\vspace{4pt}
\noindent\textbf{Covariance Error Decomposition and Quantitative Accuracy Analysis}

\begin{figure}[!htbp]
	\centering
	\includegraphics[width=0.5\textwidth]
	{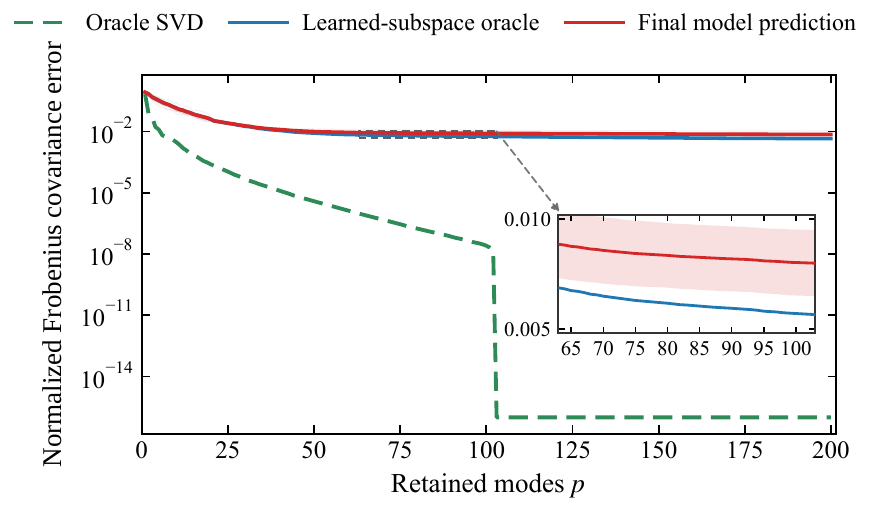}
	\caption{Normalized Frobenius covariance error as a function of the number of retained modes.}
	\label{fig:burgers-cov-decomp}
\end{figure}

Figure~\ref{fig:burgers-cov-decomp} shows how the covariance errors vary with the number of retained modes \(p\), while Table~\ref{tab:burgers-cov-thresholds} reports the minimum numbers of modes required to reach different target normalized Frobenius errors.

\begin{table}[!htbp]
	\centering
	\caption{Minimum numbers of retained modes required to achieve different target normalized Frobenius covariance-error levels for the three covariance-recovery curves.}
	\label{tab:burgers-cov-thresholds}
	\begin{tabular}{
			>{\centering\arraybackslash}p{0.24\textwidth}
			>{\centering\arraybackslash}p{0.20\textwidth}
			>{\centering\arraybackslash}p{0.25\textwidth}
			>{\centering\arraybackslash}p{0.25\textwidth}
		}
		\toprule
		Target error
		& Oracle SVD
		& Learned-subspace oracle
		& Final model prediction \\
		\midrule
		Error \(\leq 20.00\%\) & \(p\geq 2\) & \(p\geq 8\)  & \(p\geq 8\)  \\
		Error \(\leq 10.00\%\) & \(p\geq 2\) & \(p\geq 13\) & \(p\geq 13\) \\
		Error \(\leq 5.00\%\)  & \(p\geq 4\) & \(p\geq 19\) & \(p\geq 19\) \\
		Error \(\leq 2.00\%\)  & \(p\geq 4\) & \(p\geq 29\) & \(p\geq 30\) \\
		Error \(\leq 1.00\%\)  & \(p\geq 6\) & \(p\geq 43\) & \(p\geq 50\) \\
		\bottomrule
	\end{tabular}
\end{table}

\textbf{(1) Oracle SVD curve.}
The Oracle SVD error falls below \(5\%\) for \(p\geq4\) and below \(1\%\) for \(p\geq6\), indicating strong low-rank compressibility of the reference covariance. Although its decay is slower than in the reaction--diffusion problem, low-rank truncation is not the main source of covariance recovery error.

\textbf{(2) Comparison between the Learned-subspace oracle curve and the Final model prediction curve.}
The two curves follow similar decreasing trends over the full range of
\(p\). Their errors decline from approximately \(2\times10^{-1}\) to
\(8\times10^{-3}\). The enlarged inset shows that the Final model
prediction remains approximately \(1\)--\(2\times10^{-3}\) above the
Learned-subspace oracle. Two observations follow:

\begin{enumerate}[leftmargin=2.2em,itemsep=0.4em]
	\item
	The contribution of \(T_{3}\) is clearly visible but remains smaller
	than that of \(T_{2}\). The branch approximation and optimization
	errors are not the bottleneck here.
	
	\item
	The covariance recovery error is primarily governed by \(T_{2}\). For \(M=101\times101=10{,}201\) and \(K=10^{3}\), the corresponding relative statistical-error scale is
	\begin{equation}
		\frac{
			\sqrt{M}K^{-1/2}
		}{
			\left\|
			\widehat{\boldsymbol{\Sigma}}
			\right\|_{F}
		}
		\approx
		1.74\times10^{-2}.
		\label{eq:burgers-stat-err}
	\end{equation}
	This estimate is larger than the observed \(T_{2}\approx8\times10^{-3}\). Since the asymptotic relation provides only an order estimate, it does not permit a reliable separation of the contributions of \(\Delta_T\) and finite-sample statistical error within \(T_{2}\).
	The final normalized Frobenius covariance error is approximately \(8\times10^{-3}\), which is comparable in magnitude to the mean relative \(L_2\) prediction error of \(1.59\%\) reported in Table~\ref{tab:burgers-l2-train-interp}. The method therefore achieves stable and practically useful covariance recovery despite the stronger nonlinearity of the Burgers problem.
\end{enumerate}

In summary, covariance recovery for the Burgers equation exhibits the hierarchy \(T_{1}<T_{2}\) and \(T_{3}\lesssim T_{2}\). The overall recovery error is nonetheless of the same order as the mean prediction error, showing that the observed covariance recovery behavior is consistent with the theoretical decomposition in Section~\ref{sec:cov-error-decomposition}.

\vspace{4pt}
\noindent\textbf{Recovery of Local Correlation Structures:
	Multi-Anchor Correlation Map Assessment}

\begin{figure}[!htbp]
	\centering
	\includegraphics[width=0.88\textwidth]
	{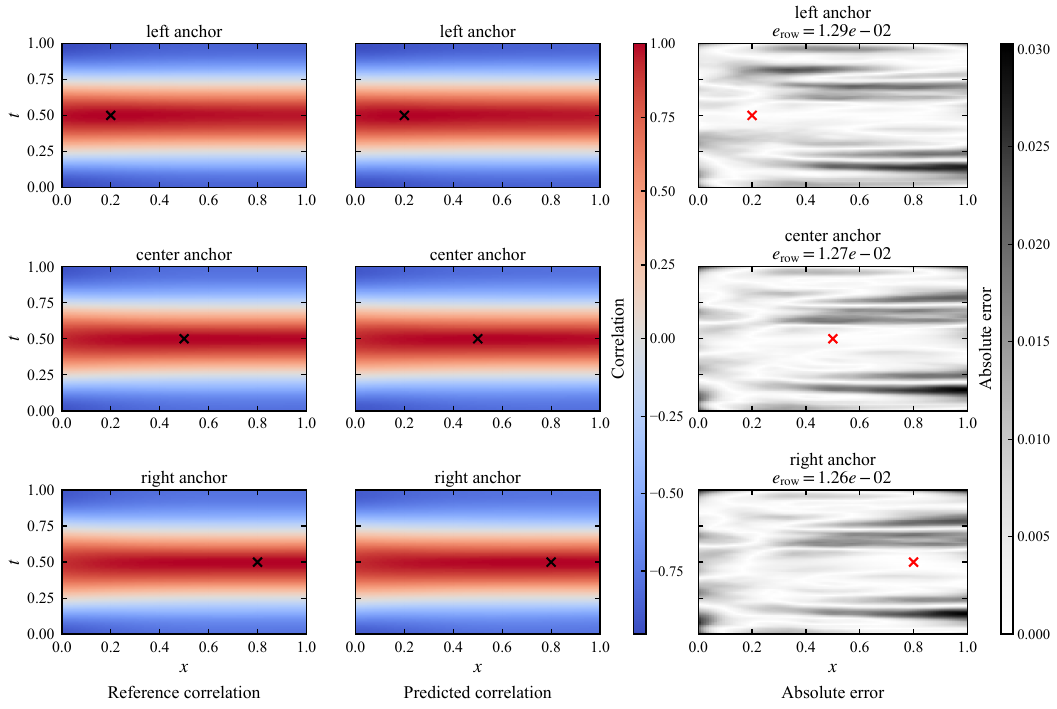}
	\caption{Multi-anchor correlation maps showing the reference correlations, predicted correlations, and corresponding absolute errors for three anchor points, where $e_{\mathrm{row}}$ denotes the relative $L_2$ error of the predicted correlation map for each anchor point.}
	\label{fig:burgers-anchor}
\end{figure}

For the Burgers equation, three anchor points are selected from the left, central, and right regions of the output domain \((x,t)\in[0,1]^2\). Fig.~\ref{fig:burgers-anchor} compares the corresponding reference, predicted, and absolute-error correlation maps.

\vspace{4pt}
\noindent\textbf{Overall pattern consistency.}
The predicted maps closely reproduce the reference patterns for all three anchors, with relatively small discrepancies in the absolute-error maps. The mean row-wise correlation error is approximately \(1.27\times10^{-2}\), confirming accurate recovery of the local correlation structure.

\vspace{4pt}
\noindent\textbf{Physical interpretation of the correlation structure.}
Unlike the vertical bands observed in the reaction--diffusion example, the Burgers correlation maps exhibit distinct horizontal bands centered near the anchor time \(t=t_a\). They reflect coordinated variations across spatial locations at similar stages of the nonlinear evolution, while the correlation gradually weakens as the temporal distance from the anchor increases. The predicted maps accurately preserve this characteristic structure.


\subsubsection{Prediction Interval Calibration and Coverage Analysis}
\label{subsec:burgers-cp-calibration}

CP is applied to calibrate the prediction intervals for the Burgers equation. The corresponding PICP and MPIW values before and after calibration are reported in Table~\ref{tab:burgers-cp-calibration}.

Before calibration, both methods exhibit substantial overcoverage. After CP calibration, their PICP values approach the nominal levels, and their interval widths decrease markedly. At the \(95\%\) nominal level, the MPIW ratio between Prob-DeepONet and two-step MV-DeepONet decreases from approximately \(4.1\) to \(2.3\). Two-step MV-DeepONet therefore retains considerably sharper intervals while achieving comparable coverage.

\begin{table}[H]
	\centering
	\caption{PICP and MPIW of Prob-DeepONet and two-step MV-DeepONet before and after CP calibration.}
	\label{tab:burgers-cp-calibration}
	
	\renewcommand{\arraystretch}{1.25}
	\setlength{\tabcolsep}{6pt}
	\footnotesize
	
	\sisetup{
		detect-all,
		table-number-alignment = center
	}
	
	\begin{tabular}{
			l
			c
			S[table-format=3.2]
			S[table-format=1.2e-2]
			S[table-format=3.2]
			S[table-format=1.2e-2]
		}
		\toprule
		\multirow{2}{*}{Method}
		& \multirow{2}{*}{\makecell{Nominal\\coverage (\%)}}
		& \multicolumn{2}{c}{Before CP calibration}
		& \multicolumn{2}{c}{After CP calibration} \\
		
		\cmidrule(lr){3-4}
		\cmidrule(lr){5-6}
		
		&
		& {PICP (\%)}
		& {MPIW}
		& {PICP (\%)}
		& {MPIW} \\
		\midrule
		
		\multirow{3}{*}{Prob-DeepONet}
		& 90
		& 99.96
		& 2.57e-2
		& 90.52
		& 4.68e-3 \\
		
		& 95
		& 99.98
		& 3.06e-2
		& 95.39
		& 6.66e-3 \\
		
		& 99
		& 100.00
		& 4.02e-2
		& 99.11
		& 1.21e-2 \\
		
		\midrule
		
		\multirow{3}{*}{\makecell[l]{two-step\\MV-DeepONet}}
		& 90
		& 99.14
		& 6.30e-3
		& 90.10
		& 2.12e-3 \\
		
		& 95
		& 99.44
		& 7.51e-3
		& 94.88
		& 2.94e-3 \\
		
		& 99
		& 99.72
		& 9.87e-3
		& 98.82
		& 5.54e-3 \\
		
		\bottomrule
	\end{tabular}
\end{table}

\FloatBarrier

\subsection{Darcy Equation}
\label{sec:darcy-introduction}

The Darcy equation is a classical elliptic model for steady incompressible flow through porous media, with applications in groundwater flow, seepage analysis, and contaminant transport modeling~\cite{Rubin2003}. Its dependence on a spatially heterogeneous permeability field makes it a standard benchmark for uncertainty propagation in random media. We consider the two-dimensional steady Darcy equation on the unit square \(\Omega=[0,1]^2\):
\begin{equation}
	\left\{
	\begin{aligned}
		-\nabla\cdot\left(
		c(\mathbf{s})\nabla u(\mathbf{s})
		\right)
		&= f(\mathbf{s}),
		\qquad \mathbf{s}\in\Omega,\\
		u(\mathbf{s})
		&= 0,
		\qquad \mathbf{s}\in\partial\Omega,
	\end{aligned}
	\right.
	\label{eq:darcy-equation}
\end{equation}
where \(\mathbf{s}=(x,y)\) denotes the spatial coordinate, \(u(\mathbf{s})\) is the pressure field, \(c(\mathbf{s})\) is the permeability field, and \(f(\mathbf{s})=50\) is the prescribed source term. Taking the random permeability $c(\mathbf{s})$ as the input and the solution $u(\mathbf{s})$ as the output response gives the target operator
\begin{equation}
	\mathcal{G}:c(\mathbf{s})\mapsto u(\mathbf{s}).
\end{equation}

\vspace{4pt}
\par\noindent\textbf{Source of Uncertainty.}
The uncertainty arises from the permeability field \(c(\mathbf{s})\), which governs the flow response but cannot generally be observed throughout the domain. It is therefore commonly modeled as a spatially correlated random field~\cite{Godoy2018}. Its variability induces spatially structured variations and correlations in the pressure field, making this problem suitable for evaluating field prediction and output covariance recovery under input uncertainty.

\vspace{4pt}
\par\noindent\textbf{Data Generation.}
The permeability is modeled as the lognormal random field \(c(\mathbf{s})=\exp(a(\mathbf{s}))\), where \(a(\mathbf{s})\) is a zero-mean Gaussian process with the squared-exponential covariance kernel
\begin{equation}
	\mathcal{C}(\mathbf{s}_1,\mathbf{s}_2)
	=
	\sigma^2
	\exp\left(
	-\frac{
		\left\|\mathbf{s}_1-\mathbf{s}_2\right\|_2^2
	}{
		2\ell^2
	}
	\right).
\end{equation}
Here, \(\sigma^2=1.0\) denotes the variance and \(\ell=0.25\) the correlation length. The latent field \(a(\mathbf{s})\) is generated using a truncated Karhunen--Loève expansion with the leading \(100\) modes, after which \(c(\mathbf{s})\) is obtained by exponentiation. For each permeability realization, Eq.~\eqref{eq:darcy-equation} is solved using the FEniCS finite element solver. The pressure field is then interpolated onto a uniform \(64\times64\) grid as the output label. A total of \(20{,}000\) input--output pairs are generated, comprising \(16{,}000\) training samples, \(2{,}000\) calibration samples, and \(2{,}000\) test samples.

\FloatBarrier

\subsubsection{Prediction Accuracy and Generalization Performance}
\label{subsec:darcy-accuracy}

\noindent\textbf{Training and In-Distribution Accuracy}
\label{subsec:darcy-train-interp}
\par\vspace{0.25em}

\begin{figure}[H]
	\centering
	\captionsetup{
		skip=2pt
	}
	
	\includegraphics[
	width=0.5\linewidth
	]{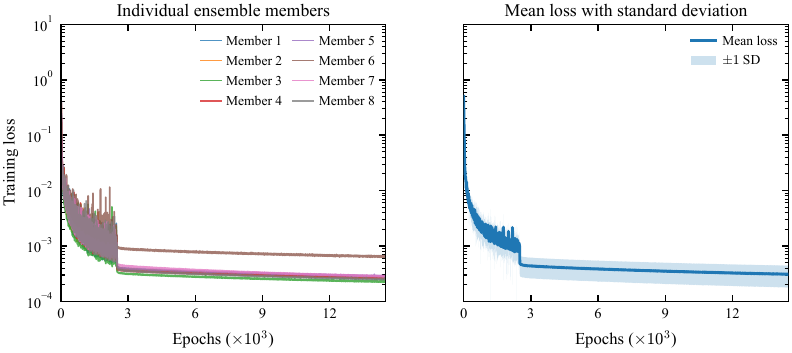}
	\caption{Training loss curves of Prob-DeepONet: individual ensemble members (left) and the ensemble mean with a $\pm 1$ standard deviation band (right).}
	\label{fig:darcy-trainloss-onestep}
\end{figure}

\begin{figure}[H]
	\centering
	
	\captionsetup{
		skip=2pt
	}
	\captionsetup[subfigure]{
		skip=1pt,
		justification=centering,
		singlelinecheck=true
	}
	
	\begin{subfigure}[t]{0.48\linewidth}
		\centering
		\begin{tikzpicture}
			\node[
			draw=black,
			dashed,
			line width=0.8pt,
			inner sep=3pt
			] {
				\includegraphics[
				width=\linewidth
				]{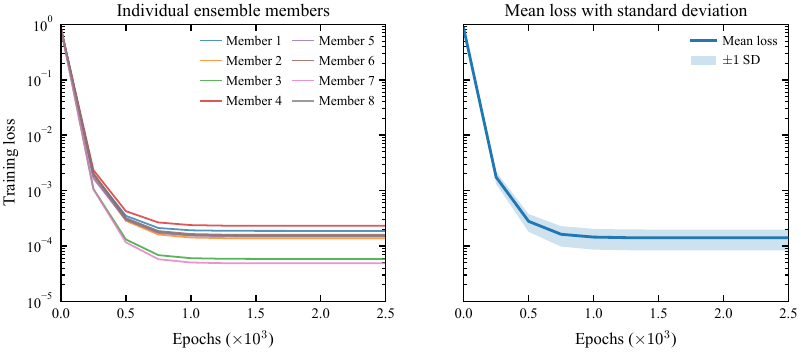}
			};
		\end{tikzpicture}
		\caption{Trunk network}
		\label{fig:darcy-trainloss-twostep-trunk}
	\end{subfigure}
	\hfill
	\begin{subfigure}[t]{0.48\linewidth}
		\centering
		\begin{tikzpicture}
			\node[
			draw=black,
			dashed,
			line width=0.8pt,
			inner sep=3pt
			] {
				\includegraphics[
				width=\linewidth
				]{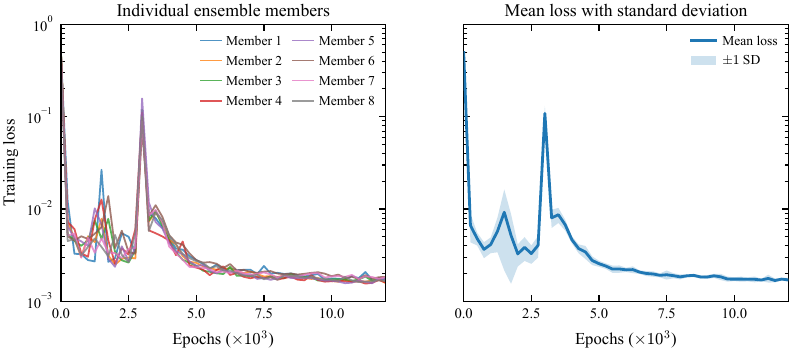}
			};
		\end{tikzpicture}
		\caption{Branch network}
		\label{fig:darcy-trainloss-twostep-branch}
	\end{subfigure}
	\caption{Training loss curves of two-step MV-DeepONet for the trunk (a) and branch (b) networks: individual ensemble members (left) and the ensemble mean with a $\pm 1$ standard deviation band (right).}
	\label{fig:darcy-trainloss-twostep}
\end{figure}

Figures~\ref{fig:darcy-trainloss-onestep}
and~\ref{fig:darcy-trainloss-twostep} show stable convergence for both
methods. The corresponding errors are summarized in
Table~\ref{tab:darcy-interp-l2}.

\vspace{-0.4em}

\vspace{-0.4em}

\begin{table}[H]
	\centering
	\captionsetup{
		skip=3pt
	}
	\caption{Mean and maximum relative $L_2$ errors on the training set and the in-distribution test set.}
	\label{tab:darcy-interp-l2}
	
	\small
	\renewcommand{\arraystretch}{1.10}
	
	\begin{tabular}{lcc}
		\toprule
		Metric
		& Prob-DeepONet
		& two-step MV-DeepONet \\
		\midrule
		
		Mean relative \(L_2\) error (training)
		& \(\textbf{1.0936\%}\)
		& \(1.3587\%\) \\
		
		Maximum relative \(L_2\) error (training)
		& \(\textbf{7.0382\%}\)
		& \(8.9639\%\) \\
		
		Mean relative \(L_2\) error (in-distribution)
		& \(3.4922\%\)
		& \(\textbf{1.7909\%}\) \\
		
		Maximum relative \(L_2\) error (in-distribution)
		& \(16.1542\%\)
		& \(\textbf{14.1562\%}\) \\
		
		\bottomrule
	\end{tabular}
\end{table}

\vspace{0.2em}

\noindent
Although Prob-DeepONet achieves lower training errors, two-step
MV-DeepONet reduces the mean in-distribution test error from \(3.4922\%\) to
\(1.7909\%\) and the maximum in-distribution test error from \(16.1542\%\) to
\(14.1562\%\), demonstrating better in-distribution performance.

\FloatBarrier

\vspace{4pt}
\noindent\textbf{Out-of-Distribution Generalization Performance}
\label{subsec:darcy-extrapolation}

For the Darcy equation, OOD performance is evaluated by varying the correlation length \(\ell\) of the random permeability field from its training value \(\ell=0.25\). The relative \(L_2\) errors are summarized in Table~\ref{tab:darcy-ood-l2}.

Prob-DeepONet exhibits large errors throughout the OOD range, with a mean error of \(69.98\%\). Two-step MV-DeepONet reduces the mean error to \(3.57\%\), approximately \(1/20\) of that of Prob-DeepONet, and remains below \(10\%\) except at \(\ell=0.105\).
Fig.~\ref{fig:darcy-ood-bands-comparison} compares representative results for rough and smooth OOD permeability fields.

As shown in Fig.~\ref{fig:darcy-ood-bands-comparison}, the bands produced by Prob-DeepONet exhibit irregular, locally jagged boundaries and often fail to cover the reference profiles, particularly under more severe OOD shifts. In contrast, the bands of two-step MV-DeepONet are smoother and better aligned with the pressure-field structure. They remain narrow for smooth OOD inputs and broaden for rougher cases while generally maintaining coverage of the reference solutions.
\begin{table}[!htbp]
	\centering
	\caption{Relative $L_2$ errors of Prob-DeepONet and two-step MV-DeepONet for out-of-distribution samples at different input correlation lengths.}
	\label{tab:darcy-ood-l2}
	\small
	\renewcommand{\arraystretch}{1.08}
	\begin{tabular}{ccc}
		\toprule
		Correlation length \(\ell\)
		& Prob-DeepONet
		& two-step MV-DeepONet \\
		\midrule
		0.06  & 78.321\% & \textbf{7.141\%} \\
		0.075 & 77.828\% & \textbf{8.765\%} \\
		0.09  & 71.807\% & \textbf{6.959\%} \\
		0.105 & 59.811\% & \textbf{10.691\%} \\
		0.12  & 69.914\% & \textbf{5.246\%} \\
		0.14  & 46.733\% & \textbf{4.429\%} \\
		0.165 & 46.040\% & \textbf{3.860\%} \\
		0.19  & 79.887\% & \textbf{2.651\%} \\
		0.22  & 40.665\% & \textbf{3.070\%} \\
		0.29  & 84.468\% & \textbf{3.084\%} \\
		0.34  & 83.440\% & \textbf{0.701\%} \\
		0.40  & 80.956\% & \textbf{1.401\%} \\
		0.47  & 76.157\% & \textbf{1.801\%} \\
		0.55  & 74.610\% & \textbf{0.743\%} \\
		0.65  & 77.185\% & \textbf{0.844\%} \\
		0.78  & 72.010\% & \textbf{1.200\%} \\
		0.93  & 70.730\% & \textbf{0.571\%} \\
		1.10  & 69.353\% & \textbf{1.097\%} \\
		\midrule
		Mean error
		& 69.98\%
		& \textbf{3.57\%} \\
		\bottomrule
	\end{tabular}
\end{table}

\FloatBarrier

\begin{figure}[!htbp]
	\centering
	
	\captionsetup[subfigure]{
		skip=2pt,
		justification=centering,
		singlelinecheck=true
	}
	
	\begin{subfigure}[t]{0.48\textwidth}
		\centering
		\includegraphics[
		width=\linewidth
		]{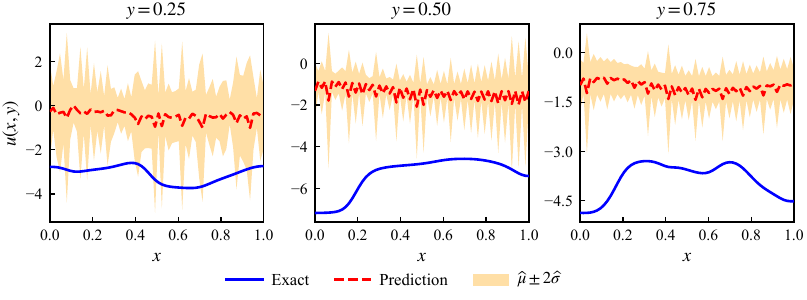}
		\caption{Prob-DeepONet, \(\ell=0.06\).}
		\label{fig:darcy-ood-band-onestep-L006}
	\end{subfigure}
	\hfill
	\begin{subfigure}[t]{0.48\textwidth}
		\centering
		\includegraphics[
		width=\linewidth
		]{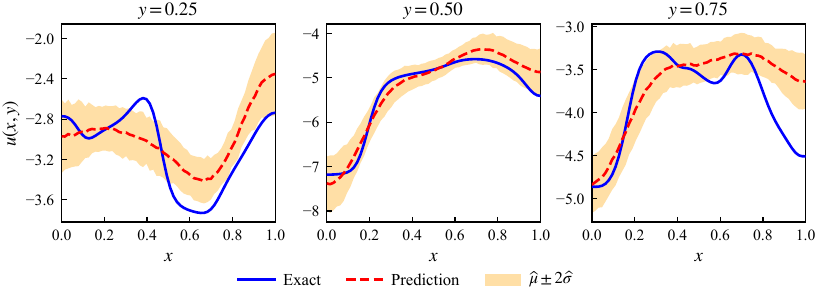}
		\caption{two-step MV-DeepONet, \(\ell=0.06\).}
		\label{fig:darcy-ood-band-twostep-L006}
	\end{subfigure}
	
	\par\vspace{-0.10em}
	
	\begin{subfigure}[t]{0.48\textwidth}
		\centering
		\includegraphics[
		width=\linewidth
		]{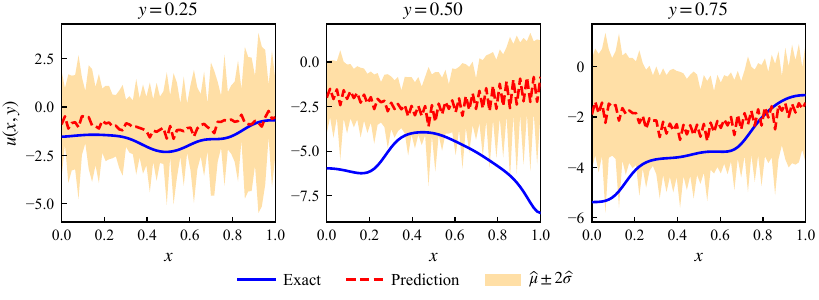}
		\caption{Prob-DeepONet, \(\ell=0.105\).}
		\label{fig:darcy-ood-band-onestep-L0105}
	\end{subfigure}
	\hfill
	\begin{subfigure}[t]{0.48\textwidth}
		\centering
		\includegraphics[
		width=\linewidth
		]{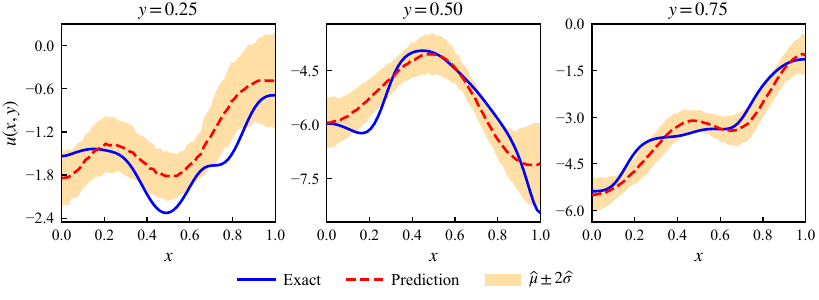}
		\caption{two-step MV-DeepONet, \(\ell=0.105\).}
		\label{fig:darcy-ood-band-twostep-L0105}
	\end{subfigure}
	
	\par\vspace{-0.10em}
	
	\begin{subfigure}[t]{0.48\textwidth}
		\centering
		\includegraphics[
		width=\linewidth
		]{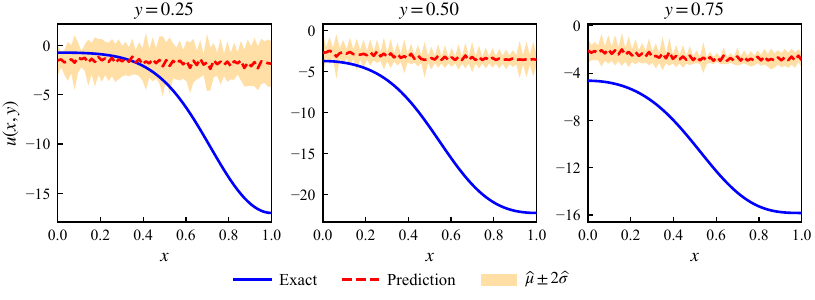}
		\caption{Prob-DeepONet, \(\ell=0.40\).}
		\label{fig:darcy-ood-band-onestep-L040}
	\end{subfigure}
	\hfill
	\begin{subfigure}[t]{0.48\textwidth}
		\centering
		\includegraphics[
		width=\linewidth
		]{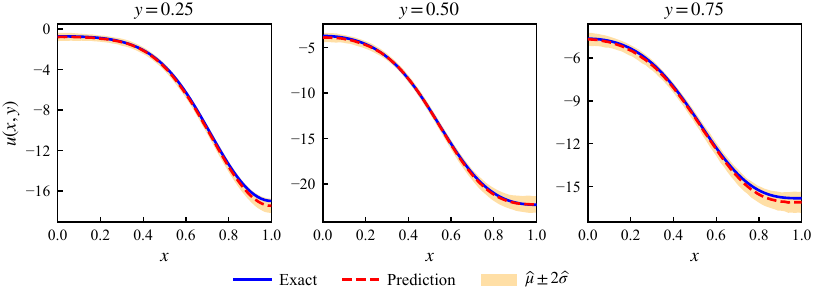}
		\caption{two-step MV-DeepONet, \(\ell=0.40\).}
		\label{fig:darcy-ood-band-twostep-L040}
	\end{subfigure}
	
	\caption{Predictions and predictive uncertainty bands of Prob-DeepONet
		(left column) and two-step MV-DeepONet (right column) for
		representative out-of-distribution samples at different input
		correlation lengths.}
	\label{fig:darcy-ood-bands-comparison}
\end{figure}

\begin{figure}[!htbp]
	\ContinuedFloat
	\centering
	
	\captionsetup[subfigure]{
		skip=2pt,
		justification=centering,
		singlelinecheck=true
	}
	
	\setcounter{subfigure}{6}
	
	\begin{subfigure}[t]{0.48\textwidth}
		\centering
		\includegraphics[
		width=\linewidth
		]{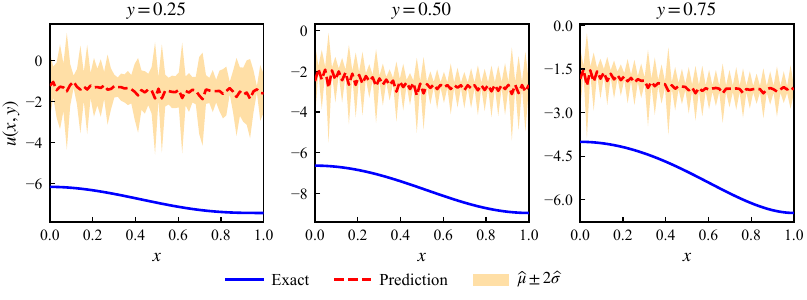}
		\caption{Prob-DeepONet, \(\ell=0.93\).}
		\label{fig:darcy-ood-band-onestep-L093}
	\end{subfigure}
	\hfill
	\begin{subfigure}[t]{0.48\textwidth}
		\centering
		\includegraphics[
		width=\linewidth
		]{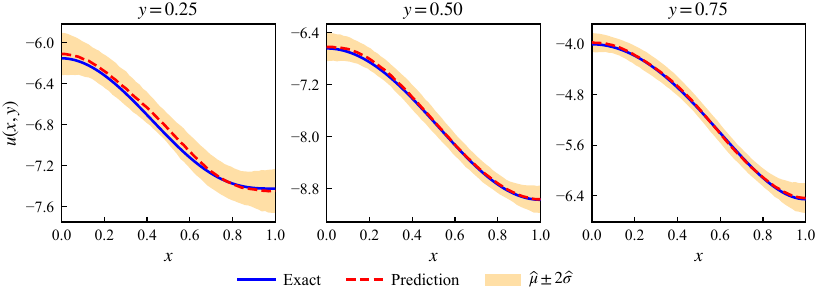}
		\caption{two-step MV-DeepONet, \(\ell=0.93\).}
		\label{fig:darcy-ood-band-twostep-L093}
	\end{subfigure}
	
	\caption[]{Predictions and predictive uncertainty bands of Prob-DeepONet
		(left column) and two-step MV-DeepONet (right column) for
		representative out-of-distribution samples at different input
		correlation lengths (continued).}
\end{figure}

\FloatBarrier

\subsubsection{Analysis of Output Covariance Recovery Accuracy}
\label{subsec:darcy-cov-recovery}

This section further quantifies the output covariance recovery accuracy for the Darcy equation.

\vspace{4pt}
\par\noindent\textbf{Covariance Error Decomposition and Quantitative Accuracy Analysis}

\begin{figure}[!htbp]
	\centering
	\includegraphics[width=0.5\textwidth]
	{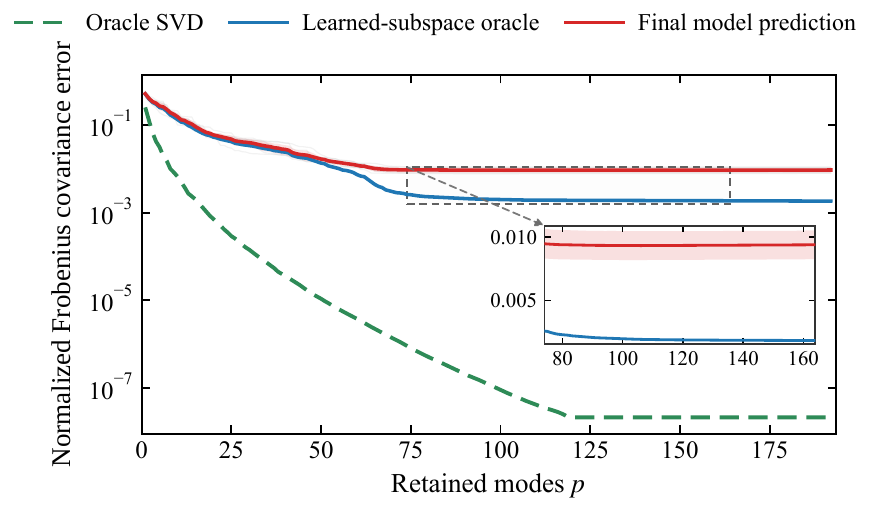}
	\caption{Normalized Frobenius covariance error as a function of the
		number of retained modes.}
	\label{fig:darcy-cov-decomp}
\end{figure}

Figure~\ref{fig:darcy-cov-decomp} shows the three normalized Frobenius covariance-error curves as a function of the number of retained modes \(p\), and Table~\ref{tab:darcy-cov-thresholds} reports the minimum numbers of modes required to reach different target normalized Frobenius errors. The main observations are as follows.

\textbf{(1) Oracle SVD curve.}
As shown in Table~\ref{tab:darcy-cov-thresholds}, \(T_{1}\) falls below \(5\%\) for \(p\geq4\) and below \(1\%\) for \(p\geq9\), indicating strong low-rank compressibility of the Darcy output covariance. Thus, the low-rank truncation error \(T_{1}\) is not the main factor limiting covariance recovery.

\textbf{(2) Comparison between the Learned-subspace oracle curve and the Final model prediction curve.}
The two curves exhibit similar decreasing trends over the full range of \(p\), with the Final model prediction curve remaining only slightly above the Learned-subspace oracle curve. This leads to two main observations:

\begin{enumerate}[leftmargin=2.2em,itemsep=0.3em,topsep=0.3em]
	\item
	The contribution of \(T_{3}\) is visible but small. The gap between
	the two curves remains much smaller than their overall error levels,
	indicating that branch prediction of the coefficient covariance is
	not the primary bottleneck.
	
	\item
	The dominant contribution comes from \(T_{2}\). For \(M=64\times64=4096\) and \(K=2\times10^{3}\), the corresponding relative finite-sample statistical-error scale is
	\begin{equation}
		\frac{
			\sqrt{M}K^{-1/2}
		}{
			\left\|
			\widehat{\boldsymbol{\Sigma}}
			\right\|_{F}
		}
		\approx
		6.2\times10^{-5}.
	\end{equation}
	This value is far smaller than the observed \(T_{2}\approx2\times10^{-3}\), indicating that \(T_{2}\) is governed mainly by the trunk-subspace learning error. Compared with the first two examples, the Darcy problem requires more retained modes to achieve the same error levels, reflecting the greater complexity of the covariance structure induced by a two-dimensional random medium. Nevertheless, the final covariance error decreases steadily to the \(10^{-2}\) level as \(p\) increases, which is comparable to the mean relative \(L_2\) error of \(1.7909\%\) reported for two-step MV-DeepONet on the in-distribution test set.
	
\end{enumerate}

Overall, covariance recovery for the Darcy problem follows the hierarchy \(T_{1}<T_{2}\) and \(T_{3}\lesssim T_{2}\), while the final recovery error is
governed primarily by trunk-subspace learning. The resulting covariance recovery accuracy is comparable to the mean prediction accuracy, demonstrating the effectiveness of the proposed method for recovering output covariance in problems involving two-dimensional random media.

\begin{table}[!htbp]
	\centering
	\caption{Minimum numbers of retained modes required to achieve different target normalized Frobenius covariance-error levels for the three covariance-recovery curves.}
	\label{tab:darcy-cov-thresholds}
	\small
	\renewcommand{\arraystretch}{1.18}
	\begin{tabular}{lccc}
		\toprule
		Target error
		& Oracle SVD
		& \makecell{Learned-subspace\\oracle}
		& \makecell{Final model\\prediction} \\
		\midrule
		Error \(\leq 20.00\%\)
		& \(p\geq2\)
		& \(p\geq8\)
		& \(p\geq8\) \\
		
		Error \(\leq 10.00\%\)
		& \(p\geq3\)
		& \(p\geq13\)
		& \(p\geq15\) \\
		
		Error \(\leq 5.00\%\)
		& \(p\geq4\)
		& \(p\geq22\)
		& \(p\geq25\) \\
		
		Error \(\leq 2.00\%\)
		& \(p\geq6\)
		& \(p\geq43\)
		& \(p\geq47\) \\
		
		Error \(\leq 1.00\%\)
		& \(p\geq9\)
		& \(p\geq56\)
		& \(p\geq66\) \\
		\bottomrule
	\end{tabular}
\end{table}


\vspace{4pt}
\par\noindent\textbf{Recovery of Local Correlation Structures:
	Multi-Anchor Correlation Map Assessment}

\begin{figure}[!htbp]
	\centering
	\includegraphics[width=0.85\textwidth]
	{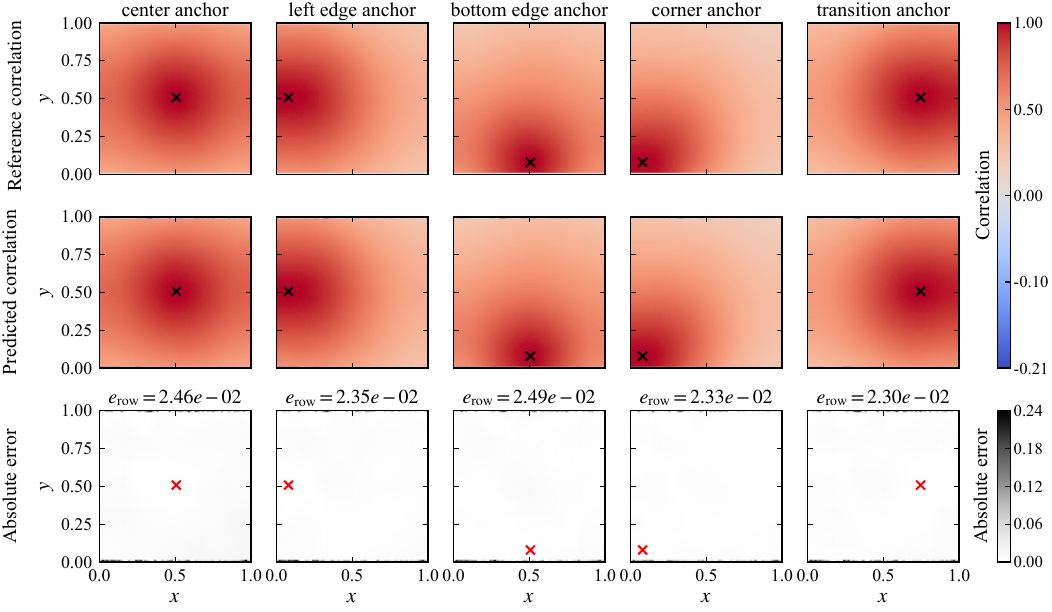}
	\caption{Multi-anchor correlation maps showing the reference correlations, predicted correlations, and corresponding absolute errors for five anchor points, where $e_{\mathrm{row}}$ denotes the relative $L_2$ error of the predicted correlation map for each anchor point.}
	\label{fig:darcy-anchor}
\end{figure}

Five anchor points \((x_a,y_a)\) are selected in the physical domain \((x,y)\in[0,1]^2\), representing the center, the left and bottom boundaries, a corner, and a transition region. Fig.~\ref{fig:darcy-anchor} compares the corresponding reference, predicted, and absolute-error correlation maps.

\vspace{4pt}
\noindent\textbf{Overall pattern consistency.}
The predicted maps reproduce the reference patterns for all five anchors, with a mean row-wise correlation error of approximately \(2.38\times10^{-2}\). This result confirms accurate recovery of the local correlation structure of the Darcy pressure field.

\vspace{4pt}
\noindent\textbf{Physical interpretation of the correlation structure.}
In contrast to the banded patterns of the preceding spatiotemporal examples, the Darcy maps exhibit localized two-dimensional correlation regions whose geometry depends on the anchor position. The interior anchor produces an approximately symmetric pattern with correlation decaying away from the anchor. Near a boundary or corner, this pattern becomes asymmetric and is confined by the domain geometry. These patterns show that the off-diagonal structure of the output covariance captures the spatial coupling induced by shared permeability perturbations through the elliptic operator. The strength, spatial extent, and geometry of the correlation patterns depend not only on spatial distance but also on the spatial correlation of the input field, the nonlocal action of the elliptic operator, and the boundary conditions.

\subsubsection{Prediction Interval Calibration and Coverage Analysis}
\label{subsec:darcy-cp-calibration}

CP is further applied to the prediction intervals for the Darcy problem. Table~\ref{tab:darcy-cp-calibration} summarizes the PICP and MPIW values before and after calibration.

\begin{table}[!htbp]
	\centering
	\caption{PICP and MPIW of Prob-DeepONet and two-step MV-DeepONet before and after CP calibration.}
	\label{tab:darcy-cp-calibration}
	
	\renewcommand{\arraystretch}{1.25}
	\setlength{\tabcolsep}{6pt}
	\footnotesize
	
	\sisetup{
		detect-all,
		table-number-alignment = center
	}
	
	\begin{tabular}{
			l
			c
			S[table-format=3.2]
			S[table-format=1.2e-1]
			S[table-format=3.2]
			S[table-format=1.2e-1]
		}
		\toprule
		\multirow{2}{*}{Method}
		& \multirow{2}{*}{\makecell{Nominal\\coverage (\%)}}
		& \multicolumn{2}{c}{Before CP calibration}
		& \multicolumn{2}{c}{After CP calibration} \\
		
		\cmidrule(lr){3-4}
		\cmidrule(lr){5-6}
		
		&
		& {PICP (\%)}
		& {MPIW}
		& {PICP (\%)}
		& {MPIW} \\
		\midrule
		
		\multirow{3}{*}{Prob-DeepONet}
		& 90
		& 99.14
		& 9.89e-2
		& 89.80
		& 4.63e-2 \\
		
		& 95
		& 99.57
		& 1.18e-1
		& 94.86
		& 5.93e-2 \\
		
		& 99
		& 99.84
		& 1.55e-1
		& 99.01
		& 9.55e-2 \\
		
		\midrule
		
		\multirow{3}{*}{\makecell[l]{two-step\\MV-DeepONet}}
		& 90
		& 98.37
		& 8.23e-3
		& 89.46
		& 4.69e-3 \\
		
		& 95
		& 99.21
		& 9.81e-3
		& 94.63
		& 5.92e-3 \\
		
		& 99
		& 99.76
		& 1.29e-2
		& 98.22
		& 9.12e-3 \\
		
		\bottomrule
	\end{tabular}
\end{table}

Both methods exhibit pronounced overcoverage before calibration. CP moves their PICP values close to the nominal levels and substantially reduces the interval
widths. Notably, the MPIW ratio between Prob-DeepONet and two-step MV-DeepONet remains large both before and after calibration, decreasing from approximately $12.0$ to $10.0$ at the $95\%$ nominal coverage level. Thus, two-step MV-DeepONet achieves similar coverage with substantially narrower intervals.

\FloatBarrier

\subsection{Hypersonic Blunt-Body Aerothermal Modeling}
\label{sec:aerothermal}

Hypersonic aerothermal prediction involves coupled multiscale and multiphysics processes in a realistic engineering setting. During atmospheric reentry or hypersonic cruise, shock compression and viscous dissipation convert kinetic energy into internal energy, producing high temperatures within the shock layer. The resulting chemical reactions and vibrational relaxation may occur on time scales comparable to the flow time scale, leading to thermochemical nonequilibrium~\cite{Passiatore2022,Williams2025}.

A representative configuration is considered here: hypersonic laminar flow past a hemispherical blunt body. The numerical configuration follows MacLean \textit{et al.}~\cite{MacLean2013} and Dsouza \textit{et al.}~\cite{Dsouza2024}. MacLean \textit{et al.} compared LENS-XX surface heat-transfer measurements with Data-Parallel Line Relaxation (DPLR) predictions, whereas Dsouza \textit{et al.} reproduced the same case in Ansys Fluent using a two-temperature
nonequilibrium model and validated the predictions against the experimental data. This example is used to evaluate structured output covariance recovery under realistic aerothermal conditions. 

\vspace{4pt}
\noindent\textbf{Governing Equations and Constitutive Relations.}
Under the assumption of laminar flow, the three-dimensional
thermochemical nonequilibrium Navier--Stokes equations are written in
conservation form as follows~\cite{Passiatore2022}:
\begin{equation}
	\begin{aligned}
		&\frac{\partial \rho}{\partial t}
		+\frac{\partial (\rho u_j)}{\partial x_j}
		=0,
		\\[2pt]
		&\frac{\partial (\rho u_i)}{\partial t}
		+\frac{\partial
			\left(\rho u_i u_j+p\delta_{ij}\right)}
		{\partial x_j}
		=
		\frac{\partial \tau_{ij}}{\partial x_j},
		\\[2pt]
		&\frac{\partial (\rho E)}{\partial t}
		+\frac{\partial
			\left[(\rho E+p)u_j\right]}
		{\partial x_j}
		=
		\frac{\partial (u_i\tau_{ij})}{\partial x_j}
		-\frac{\partial
			\left(q_j^{TR}+q_j^{V}\right)}
		{\partial x_j}
		-\frac{\partial}{\partial x_j}
		\left(
		\sum_{n=1}^{N_S}
		\rho_n u_{nj}^{D}h_n
		\right),
		\\[2pt]
		&\frac{\partial \rho_n}{\partial t}
		+\frac{\partial (\rho_n u_j)}{\partial x_j}
		=
		-\frac{\partial
			\left(\rho_n u_{nj}^{D}\right)}
		{\partial x_j}
		+\dot{\omega}_n,
		\qquad n=1,\ldots,N_S-1,
		\\[2pt]
		&\frac{\partial (\rho e_V)}{\partial t}
		+\frac{\partial (\rho e_V u_j)}{\partial x_j}
		=
		\frac{\partial}{\partial x_j}
		\left(
		-q_j^{V}
		-\sum_{m=1}^{N_M}
		\rho_m u_{mj}^{D}e_{Vm}
		\right)
		+\sum_{m=1}^{N_M}
		\left(
		Q_{TVm}+\dot{\omega}_m e_{Vm}
		\right).
	\end{aligned}
	\label{eq:aerothermal-governing-equations}
\end{equation}
Here, $E$ is the specific total energy, while $q_j^{TR}$ and
$q_j^{V}$ denote the translational--rotational and vibrational heat
fluxes, respectively. For species $n$, $\rho_n$, $u_{nj}^{D}$,
$h_n$, and $\dot{\omega}_n$ denote the partial density, diffusion
velocity, specific enthalpy, and production rate, respectively. The
quantities $e_V$ and $Q_{TVm}$ represent the specific vibrational
energy of the mixture and the energy exchange between the
translational and vibrational modes of molecular species $m$,
respectively.
Closure of the governing equations is provided by the Park
two-temperature model~\cite{Park1989} and a finite-rate chemistry model~\cite{Gupta1990,Gnoffo1989}. The former
accounts for thermal nonequilibrium between the
translational--rotational and vibrational--electronic energy modes,
whereas the latter describes species production and depletion through
finite-rate chemical reactions.

\vspace{4pt}
\noindent\textbf{Sources of Uncertainty.}
Variations in freestream density, temperature, pressure, and velocity
can alter the shock-layer structure and wall heat flux. Previous
uncertainty quantification studies of hypersonic aerothermodynamics
have treated these quantities as uncertain inputs, with freestream
velocity and density identified as important contributors to surface-heating uncertainty~\cite{Rataczak2024,Capriati2022,Lu2022}. Since wall heat flux is a
critical design load for thermal protection systems, the present study
represents freestream uncertainty through the Mach number, while the
remaining freestream quantities are determined consistently from the
prescribed atmospheric and flight conditions.

\vspace{4pt}
\noindent\textbf{Computational Fluid Dynamics (CFD) Solver Validation.}
Since the dataset is generated by CFD simulations, the solver is first
validated against the Mach \(12.4\) hemispherical blunt-body experiment
of MacLean \textit{et al.}~\cite{MacLean2013}. The configuration uses
the \(7.6\)-cm-diameter stainless-steel test article shown in their
Fig.~1(a), with mesh refinement near the stagnation region and within
the boundary layer~\cite{MacLean2013}.
The simulation employs an 11-species, 21-reaction thermochemical
nonequilibrium model together with a two-temperature formulation. The
freestream and boundary conditions are summarized in
Table~\ref{tab:aerothermal_params}.

\begin{figure}[!htbp]
	\centering
	\includegraphics[width=0.3\linewidth]
	{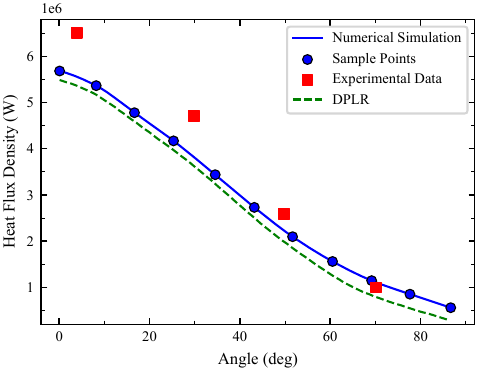}
	\caption{Comparison of wall heat flux distributions obtained from the present CFD simulation, the LENS-XX experiment, and NASA DPLR. The angle $\theta$ is measured along the hemispherical surface from the stagnation point, with $\theta=180s_{\mathrm{arc}}/(\pi R)$, where $s_{\mathrm{arc}}$ is the surface arc length and $R$ is the hemisphere radius. Thus, $\theta=0^\circ$ denotes the stagnation point and $\theta=90^\circ$ the hemisphere--cylinder junction.}
	\label{fig:aerothermal_validation}
\end{figure}

As shown in Fig.~\ref{fig:aerothermal_validation}, the predicted wall
heat flux agrees well with the experimental measurements and NASA DPLR
results~\cite{MacLean2013,Wright1998}. The simulation captures both the
stagnation-point heat-flux peak and its downstream decay, consistent
with previous numerical results for the same configuration~\cite{Dsouza2024}.

\begin{table}[!htbp]
	\centering
	\caption{Parameter settings for the high-fidelity aerothermal
		simulation.}
	\label{tab:aerothermal_params}
	
	\small
	\renewcommand{\arraystretch}{1.18}
	\setlength{\tabcolsep}{5pt}
	
	\begin{tabular}{
			@{}
			p{0.22\textwidth}
			p{0.25\textwidth}
			p{0.45\textwidth}
			@{}
		}
		\toprule
		Category
		& Parameter
		& Value \\
		\midrule
		
		\multirow{3}{*}{\makecell[l]{Physical\\configuration}}
		& Geometry
		& Stainless steel hemisphere, $D=7.6~\mathrm{cm}$ \\
		
		& Flow regime
		& Laminar \\
		
		& Thermochemical model
		& 11 chemical species, 21 elementary reactions, and a
		two-temperature model \\
		
		\midrule
		
		\multirow{5}{*}{\makecell[l]{Freestream\\conditions}}
		& Mach number
		& $Ma=12.4$ \\
		
		& Temperature
		& $T_{\infty}=535~\mathrm{K}$ \\
		
		& Pressure
		& $p_{\infty}=178.1~\mathrm{Pa}$ \\
		
		& Velocity
		& $U_{\infty}=5732~\mathrm{m\,s^{-1}}$ \\
		
		& Density
		& $\rho_{\infty}
		=1.16\times10^{-3}~\mathrm{kg\,m^{-3}}$ \\
		
		\midrule
		
		\multirow{3}{*}{\makecell[l]{Boundary\\conditions}}
		& Inlet
		& Pressure far-field \\
		
		& Outlet
		& Pressure outlet \\
		
		& Wall
		& Isothermal wall, $T_{w}=300~\mathrm{K}$ \\
		
		\bottomrule
	\end{tabular}
\end{table}

\vspace{4pt}
\noindent\textbf{Data Generation and Network Architecture}

Based on the validated CFD model, the dataset construction and network design for this example are described below. 
Taking the freestream conditions as the source of uncertainty, Latin
hypercube sampling is used to sample the freestream Mach number over
$10.4\leq Ma\leq16.4$. A total of 179
independent samples are generated and divided into 125 training
samples, 24 test samples, and 30 calibration samples. The trunk network
takes the wall-node coordinates as input and learns basis functions for
representing the wall heat flux distribution.

During the first training stage, a conventional mean squared error loss
may smooth the stagnation-point heat flux peak. Since wall heat flux is
related to the near-wall temperature gradients through
\begin{equation}
	q_{w}
	=
	q_{\mathrm{tr}}+q_{\mathrm{ve}}
	=
	k_{\mathrm{tr}}\nabla T_{\mathrm{tr}}
	+
	k_{\mathrm{ve}}\nabla T_{\mathrm{ve}},
\end{equation}
a physics-informed weighting based on the
translational--rotational temperature gradient is introduced into the
trunk reconstruction loss:
\begin{equation}
	\left(
	\omega^{*},
	\mathbf{A}^{*}
	\right)
	=
	\arg\min_{\omega,\mathbf{A}}
	\left\|
	\mathbf{W}_{\nabla T}
	\odot
	\left(
	\Phi_{\omega}\mathbf{A}
	-
	\mathbf{S}
	\right)
	\right\|_{F}^{2}.
	\qquad
	\mathbf{W}_{\nabla T}
	=
	1+\varepsilon
	\tanh\left(
	\frac{
		\left|\nabla T_{\mathrm{tr}}\right|
	}{
		\operatorname{std}\left(
		\left|\nabla T_{\mathrm{tr}}\right|
		\right)
	}
	\right).
	\label{eq:aero-trunk-weighted-loss}
\end{equation}
The translational--rotational temperature gradient is adopted because
it is numerically more stable than the vibrational--electronic
temperature gradient and is more directly related to wall heat transfer
under the cold-wall condition. The hyperbolic tangent limits excessive
weight amplification, while \(\varepsilon\) controls the weighting
strength. Setting \(\varepsilon=1.0\) gives \(1\leq\mathbf W_{\nabla T}<2\). This weighted loss emphasizes the stagnation and post-shock regions without destabilizing the reconstruction loss. Wall heat flux depends not only on the freestream conditions but also strongly on the thermochemical state of the shock layer. A multi-branch architecture is therefore used~\cite{MIONet_Jin2022}. Branch~1 receives the freestream Mach number, whereas Branch~2 receives features extracted from the mass fraction fields of the \(11\) chemical species.

The species fields are high dimensional and strongly correlated because of mass conservation,
\(\sum_{n=1}^{N_S}Y_n=1\), and finite-rate chemical coupling. A variational
autoencoder (VAE) is therefore used to encode the species mass
fraction fields into a low-dimensional latent vector
$\boldsymbol{z}$, which serves as the effective input to Branch~2.
The outputs of the two branches are concatenated and processed by a
multilayer perceptron (MLP). The resulting representation is then used
by two-step MV-DeepONet to predict the wall heat flux distribution and
recover its covariance structure. Details of the VAE encoding procedure are provided in Appendix~A.

\subsubsection{Prediction Accuracy and In-Distribution Generalization Performance}
\label{sec:aero-accuracy-generalization}
Compared with the preceding three standard benchmarks based on partial differential equations (PDEs), the
hypersonic aerothermal example is more complex in both input
dimensionality and the structure of the output field. The wall heat flux distribution over the hemispherical
blunt body exhibits pronounced spatial gradients near the
stagnation point and within the shock-affected region. This section compares the Prob-DeepONet and
two-step MV-DeepONet in terms of training accuracy, in-distribution
generalization, and uncertainty bands for challenging samples.

\vspace{4pt}
\noindent\textbf{Training and In-Distribution Generalization Accuracy}

\begin{figure}[!htbp]
	\centering
	\includegraphics[
	width=0.5\linewidth
	]{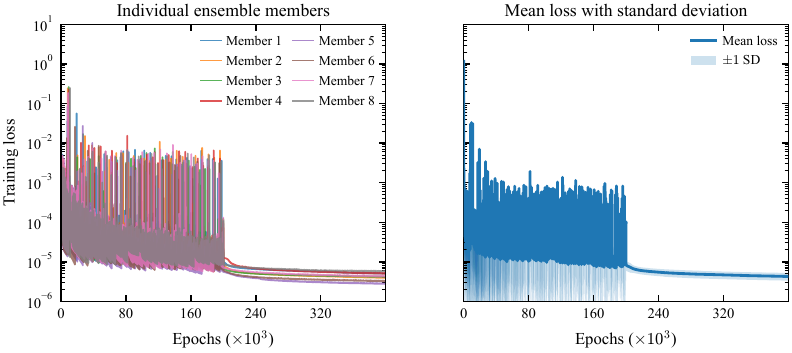}
	\caption{Training loss curves of Prob-DeepONet: individual ensemble members (left) and the ensemble mean with a $\pm 1$ standard deviation band (right).}
	\label{fig:aero-trainloss-onestep}
\end{figure}

\begin{figure}[!htbp]
	\centering
	
	\captionsetup[subfigure]{
		skip=2pt,
		justification=centering,
		singlelinecheck=true
	}
	
	\begin{subfigure}[t]{0.48\linewidth}
		\centering
		\begin{tikzpicture}
			\node[
			draw=black,
			dashed,
			line width=0.8pt,
			inner sep=3pt
			] {
				\includegraphics[
				width=\linewidth
				]{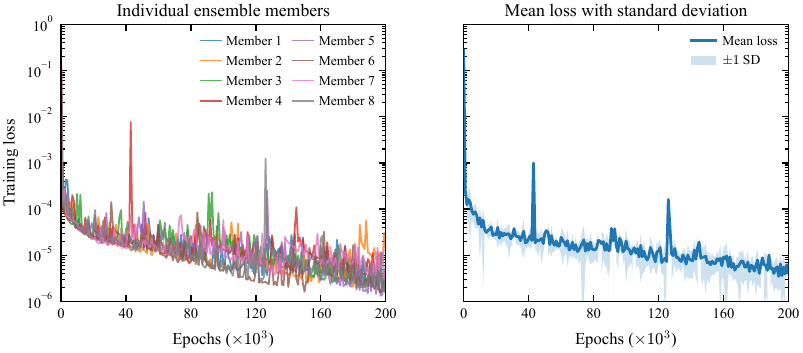}
			};
		\end{tikzpicture}
		\caption{Trunk network}
		\label{fig:aero-trainloss-twostep-trunk}
	\end{subfigure}
	\hfill
	\begin{subfigure}[t]{0.48\linewidth}
		\centering
		\begin{tikzpicture}
			\node[
			draw=black,
			dashed,
			line width=0.8pt,
			inner sep=3pt
			] {
				\includegraphics[
				width=\linewidth
				]{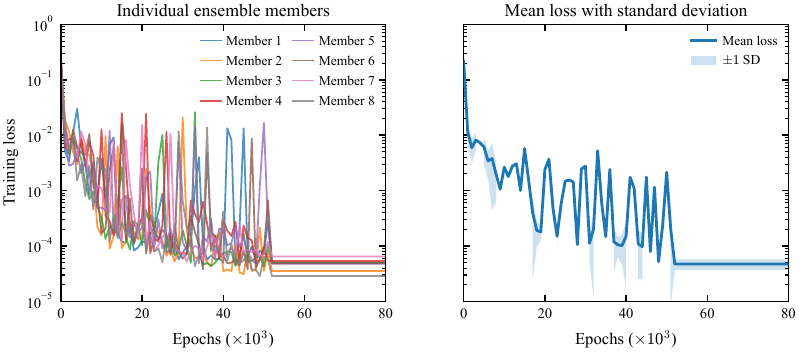}
			};
		\end{tikzpicture}
		\caption{Branch network}
		\label{fig:aero-trainloss-twostep-branch}
	\end{subfigure}
	\caption{Training loss curves of two-step MV-DeepONet for the trunk (a) and branch (b) networks: individual ensemble members (left) and the ensemble mean with a $\pm 1$ standard deviation band (right).}
	\label{fig:aero-trainloss-twostep}
\end{figure}

Figures~\ref{fig:aero-trainloss-onestep}
and~\ref{fig:aero-trainloss-twostep} show stable convergence for both
methods. Their training and in-distribution test errors are summarized
in Table~\ref{tab:aero-l2-train-test}.

\begin{table}[!htbp]
	\centering
	\caption{Mean and maximum relative $L_2$ errors on the training set and the in-distribution test set.}
	\label{tab:aero-l2-train-test}
	
	\small
	\renewcommand{\arraystretch}{1.2}
	\setlength{\tabcolsep}{8pt}
	
	\begin{tabular}{@{}lcc@{}}
		\toprule
		Metric
		& Prob-DeepONet
		& two-step MV-DeepONet \\
		\midrule
		
		Mean relative \(L_2\) error (training)
		& 0.5496\%
		& \textbf{0.2555\%} \\
		
		Maximum relative \(L_2\) error (training)
		& 0.8589\%
		& \textbf{0.3729\%} \\
		
		Mean relative \(L_2\) error (in-distribution test)
		& 3.7174\%
		& \textbf{2.1064\%} \\
		
		Maximum relative \(L_2\) error (in-distribution test)
		& 8.8195\%
		& \textbf{3.5545\%} \\
		
		\bottomrule
	\end{tabular}
\end{table}

\FloatBarrier

As reported in Table~\ref{tab:aero-l2-train-test}, two-step MV-DeepONet achieves lower mean and maximum relative \(L_2\) errors on both the training and in-distribution test sets. The increases from training to test errors are \(1.8509\) and \(3.1816\) percentage points, respectively, compared with \(3.1678\) and \(7.9606\) percentage points for Prob-DeepONet. These results demonstrate that two-step MV-DeepONet achieves higher fitting accuracy and more stable in-distribution generalization, particularly for the worst-case samples.

\vspace{4pt}
\noindent\textbf{Predictions and Uncertainty Bands for Challenging Samples}

\begin{figure}[!htbp]
	\centering
	
	\captionsetup[subfigure]{
		skip=2pt,
		justification=centering,
		singlelinecheck=true
	}
	
	\begin{subfigure}[t]{0.485\textwidth}
		\centering
		\figmaybe[
		width=\linewidth
		]{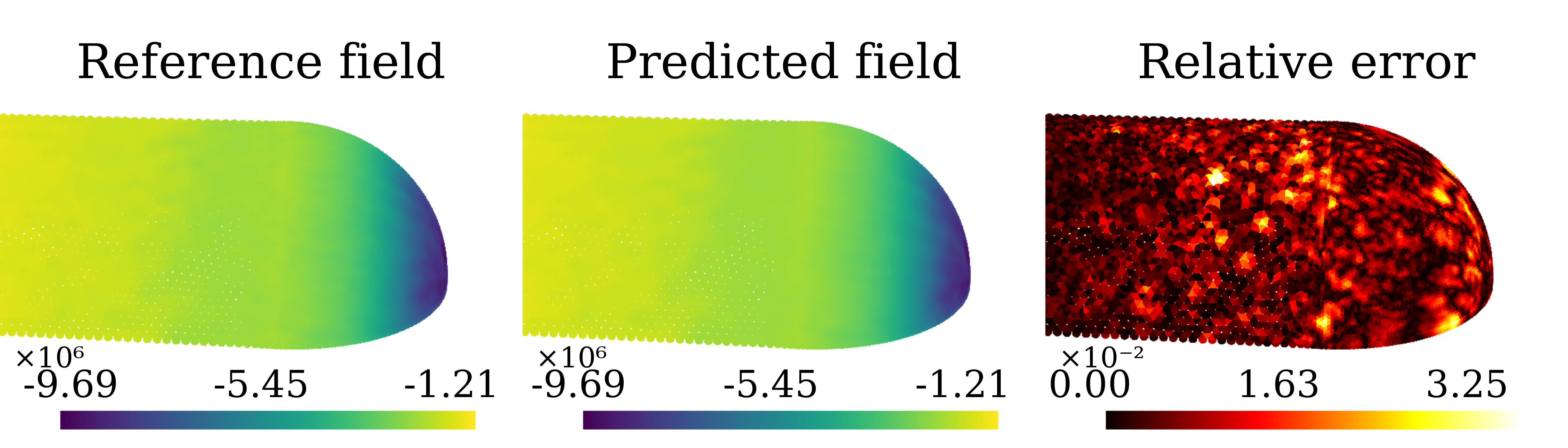}
		\caption{Prob-DeepONet: training sample with the largest
			relative \(L_2\) error.}
		\label{fig:aero-worst-onestep-train}
	\end{subfigure}%
	\hfill
	\begin{subfigure}[t]{0.485\textwidth}
		\centering
		\figmaybe[
		width=\linewidth
		]{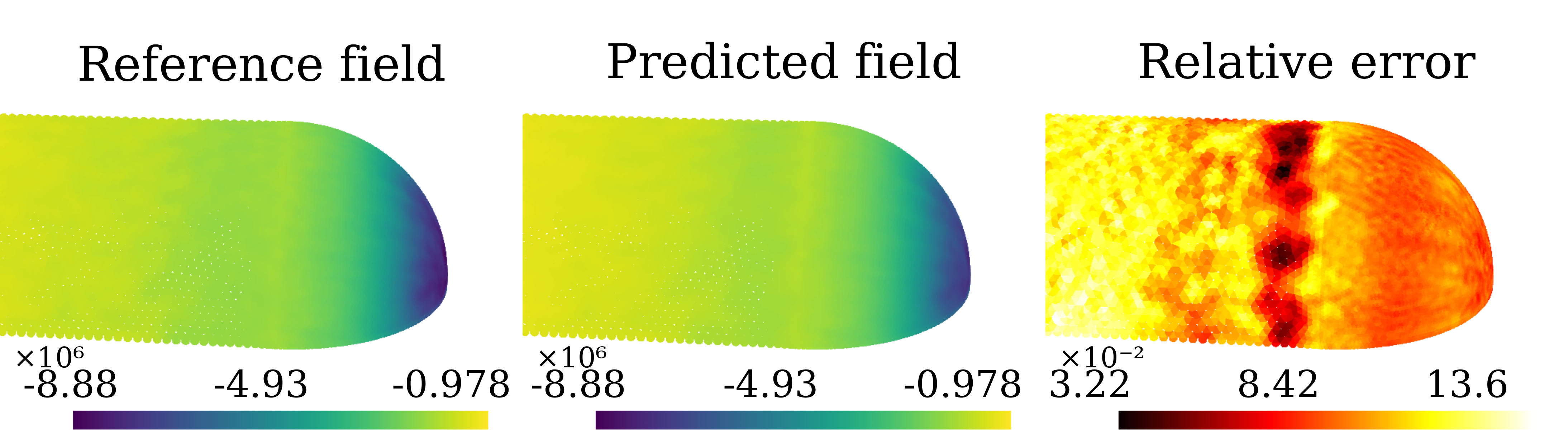}
		\caption{Prob-DeepONet: in-distribution test sample with the
			largest relative \(L_2\) error.}
		\label{fig:aero-worst-onestep-test}
	\end{subfigure}
	
	\par\vspace{0.6em}
	
	\begin{subfigure}[t]{0.485\textwidth}
		\centering
		\figmaybe[
		width=\linewidth
		]{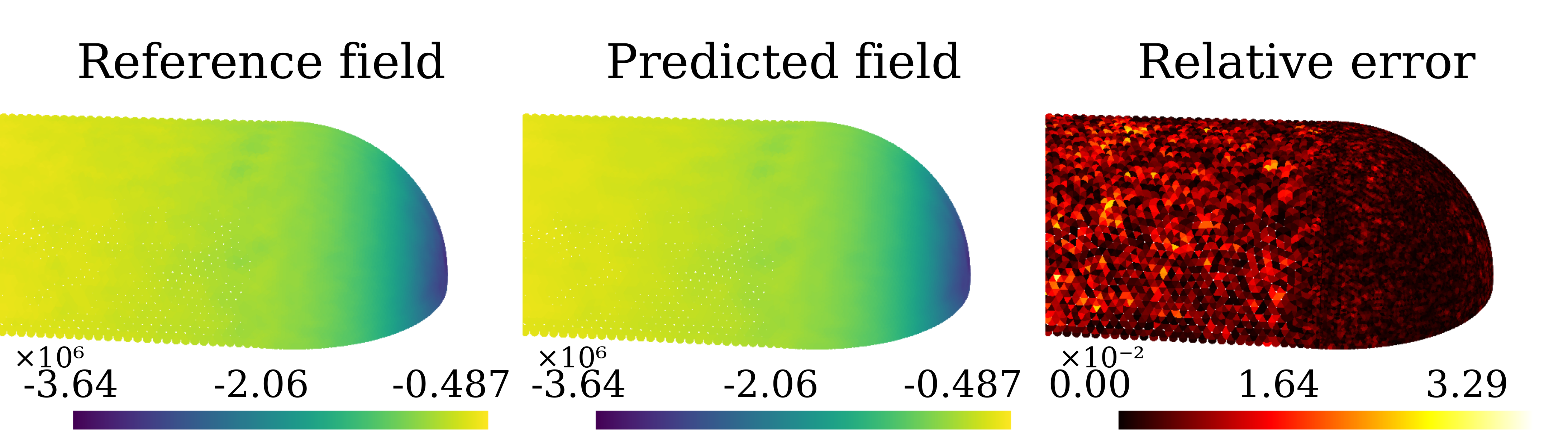}
		\caption{two-step MV-DeepONet: training sample with the largest
			relative \(L_2\) error.}
		\label{fig:aero-worst-twostep-train}
	\end{subfigure}%
	\hfill
	\begin{subfigure}[t]{0.485\textwidth}
		\centering
		\figmaybe[
		width=\linewidth
		]{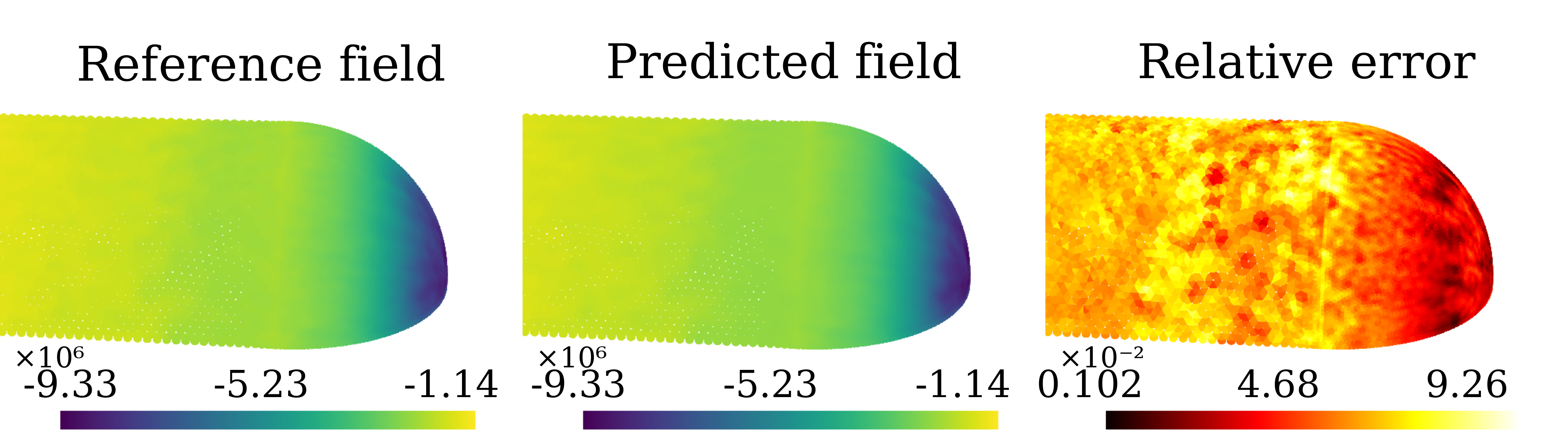}
		\caption{two-step MV-DeepONet: in-distribution test sample with
			the largest relative \(L_2\) error.}
		\label{fig:aero-worst-twostep-test}
	\end{subfigure}
	
	\caption{Reference and predicted wall heat flux fields, together with the corresponding relative error distributions, for the training and in-distribution test samples with the largest relative $L_2$ errors. The first and second rows correspond to Prob-DeepONet and two-step MV-DeepONet, respectively.}
	\label{fig:aero-worst-fields}
\end{figure}

\begin{figure}[!htbp]
	\centering
	\captionsetup[subfigure]{
		skip=2pt,
		justification=centering,
		singlelinecheck=true
	}
	
	\begin{subfigure}[t]{0.4\textwidth}
		\centering
		\figmaybe[width=\linewidth]{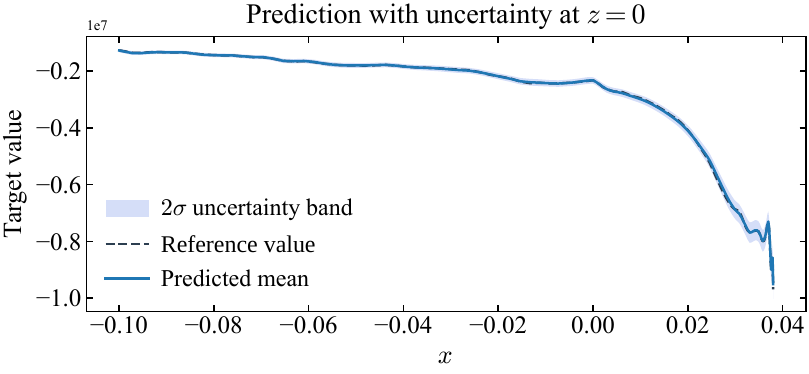}
		\caption{Prob-DeepONet: training sample with the largest
			relative \(L_2\) error.}
		\label{fig:aero-band-worst-onestep-train}
	\end{subfigure}
	\hspace{0.04\textwidth}
	\begin{subfigure}[t]{0.4\textwidth}
		\centering
		\figmaybe[width=\linewidth]{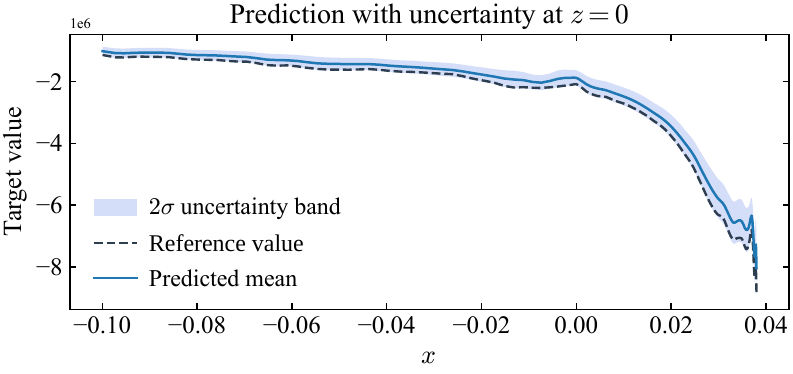}
		\caption{Prob-DeepONet: in-distribution test sample with the
			largest relative \(L_2\) error.}
		\label{fig:aero-band-worst-onestep-test}
	\end{subfigure}
	
	\par\vspace{0.5em}
	
	\begin{subfigure}[t]{0.4\textwidth}
		\centering
		\figmaybe[width=\linewidth]{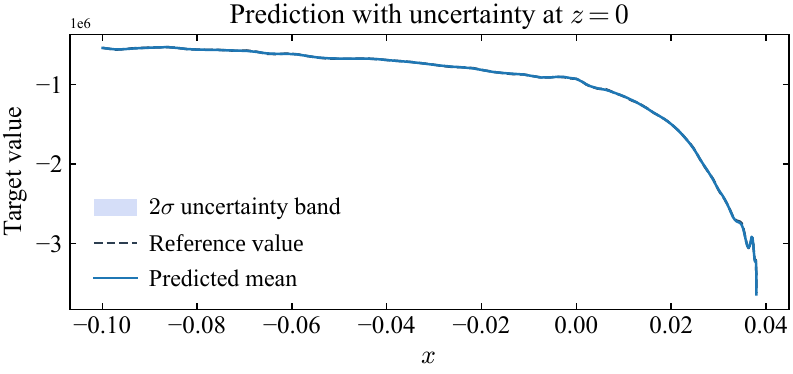}
		\caption{two-step MV-DeepONet: training sample with the largest
			relative \(L_2\) error.}
		\label{fig:aero-band-worst-twostep-train}
	\end{subfigure}
	\hspace{0.04\textwidth}
	\begin{subfigure}[t]{0.4\textwidth}
		\centering
		\figmaybe[width=\linewidth]{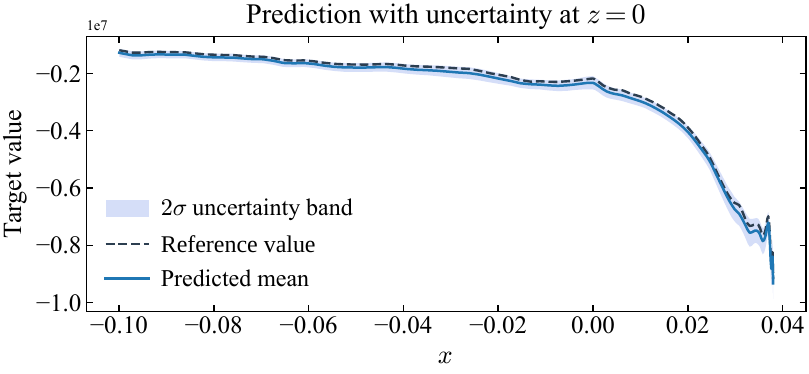}
		\caption{two-step MV-DeepONet: in-distribution test sample with
			the largest relative \(L_2\) error.}
		\label{fig:aero-band-worst-twostep-test}
	\end{subfigure}
	
	\caption{Predicted mean and uncertainty bands for the training and in-distribution test samples with the largest relative $L_2$ errors. The uncertainty bands correspond to $\pm2\sigma$. The first and second rows correspond to Prob-DeepONet and two-step MV-DeepONet, respectively.}
	\label{fig:aero-band-worst}
\end{figure}

To further examine the predictive behavior of the two methods on
challenging samples, the training and in-distribution test samples
with the largest relative $L_2$ errors under each method are selected
as stress-test cases. Fig.~\ref{fig:aero-worst-fields} presents
the reference fields, predicted fields, and relative error
distributions for these samples. Fig.~\ref{fig:aero-band-worst} shows the corresponding
predictive means and $\pm 2\sigma$ uncertainty bands along selected
spatial profiles.

Regarding the wall heat flux fields, both methods capture the dominant
variation along the wall, but their relative error distributions differ
markedly. For the most challenging in-distribution test sample, Prob-DeepONet
exhibits localized errors with a patch-like pattern. In contrast, the errors of two-step MV-DeepONet vary more smoothly along
the wall and generally increase from the stagnation region toward the downstream surface.

This difference can be explained by the different spaces in which
probabilistic modeling is performed. Prob-DeepONet evaluates the Gaussian NLL pointwise in the physical output space. For an input
\(\mathbf{u}\) and the \(i\)-th output location \(y_i\), its negative
log-likelihood can be written as
\begin{equation}
	\mathcal{L}_{\mathrm{NLL}}^{\mathrm{Prob}}(\mathbf{u})
	=
	\frac{1}{2M}
	\sum_{i=1}^{M}
	\left[
	\log \hat{\sigma}_{s,i}^{2}(\mathbf{u})
	+
	\frac{
		\left[
		s_i(\mathbf{u})
		-
		\hat{\mu}_{s,i}(\mathbf{u})
		\right]^2
	}{
		\hat{\sigma}_{s,i}^{2}(\mathbf{u})
	}
	\right],
\end{equation}
where \(\hat{\mu}_{s,i}\) and
\(\hat{\sigma}_{s,i}^{2}\) denote the predictive mean and variance
at \(y_i\), respectively. The gradient of this loss with respect to the
predictive mean is
\begin{equation}
	\frac{
		\partial \mathcal{L}_{\mathrm{NLL}}^{\mathrm{Prob}}
	}{
		\partial \hat{\mu}_{s,i}
	}
	=
	\frac{1}{M}
	\frac{
		\hat{\mu}_{s,i}-s_i
	}{
		\hat{\sigma}_{s,i}^{2}
	}.
\end{equation}
Thus, the mean residual at each location is weighted by the
corresponding pointwise variance. When the predicted variance at a
particular location is large, the gradient and training penalty
associated with the mean residual at that location may be reduced,
so that the local mean error may remain insufficiently corrected. This
pointwise mean--variance tradeoff may consequently produce relatively
isolated error patches in local regions that are difficult to fit.

By contrast, two-step MV-DeepONet evaluates the negative
log-likelihood in the modal coefficient space. For the \(m\)-th target
modal coefficient \(c_m^{*}(\mathbf{u})\), the loss is written as
\begin{equation}
	\mathcal{L}_{\mathrm{NLL}}^{\mathrm{2step}}(\mathbf{u})
	=
	\frac{1}{2p}
	\sum_{m=1}^{p}
	\left[
	\log \hat{\sigma}_{c,m}^{2}(\mathbf{u})
	+
	\frac{
		\left[
		c_m^{*}(\mathbf{u})
		-
		\hat{\mu}_{c,m}(\mathbf{u})
		\right]^2
	}{
		\hat{\sigma}_{c,m}^{2}(\mathbf{u})
	}
	\right],
\end{equation}
and its gradient with respect to the predicted coefficient mean is
\begin{equation}
	\frac{
		\partial \mathcal{L}_{\mathrm{NLL}}^{\mathrm{2step}}
	}{
		\partial \hat{\mu}_{c,m}
	}
	=
	\frac{1}{p}
	\frac{
		\hat{\mu}_{c,m}-c_m^{*}
	}{
		\hat{\sigma}_{c,m}^{2}
	}.
\end{equation}
It follows that the predicted variance in two-step MV-DeepONet
adjusts the weighting of modal coefficient mean errors rather than
that of local residuals at individual output locations. Although a
larger coefficient variance may likewise reduce the training penalty
on the corresponding coefficient mean error, this tradeoff occurs at the modal-coefficient level rather than in the physical output space.

To further illustrate this distinction, let
\begin{equation}
\delta\boldsymbol{\mu}_{c}(\mathbf{u})
=
\hat{\boldsymbol{\mu}}_{c}(\mathbf{u})
-
\mathbf{c}^{*}(\mathbf{u})
\end{equation}
denote the coefficient mean error. When only the coefficient regression
error is considered, the corresponding mean error in the physical
output space is
\begin{equation}
	\delta\boldsymbol{\mu}_{s}(\mathbf{u})
	=
	\widetilde{\mathbf{Q}}^{*}
	\delta\boldsymbol{\mu}_{c}(\mathbf{u}).
\end{equation}
At any two output locations \(y_i\) and \(y_j\), this relation gives
\begin{equation}
	\delta\mu_s(y_i,\mathbf{u})
	=
	\sum_{m=1}^{p}\widetilde{q}_{m}^{*}(y_i)\delta\mu_{c,m}(\mathbf{u}),
	\qquad
	\delta\mu_s(y_j,\mathbf{u})
	=
	\sum_{m=1}^{p}\widetilde{q}_{m}^{*}(y_j)\delta\mu_{c,m}(\mathbf{u}).
\end{equation}
These expressions show that each coefficient mean error
\(\delta\mu_{c,m}\) is distributed across multiple output locations
according to the spatial profile of the corresponding shared basis
function \(\widetilde{q}_{m}^{*}\). If the shared basis functions vary
smoothly over the output domain, the resulting mean-prediction error
field also tends to vary coherently across locations, rather than
forming isolated patch-like concentrations.

Consequently, transferring probabilistic modeling to the coefficient
space and combining it with the two-step training strategy shifts the mean--variance tradeoff from the
pointwise level in the physical output space to the low-dimensional
modal level. This provides a structural explanation for the smoother
and more spatially coordinated error distributions produced by
two-step MV-DeepONet. The gradual increase in relative error away from
the stagnation region may also be amplified by the decreasing magnitude
of the reference wall heat flux.

For the training samples, the predictive means of both methods closely
follow the reference profiles, and their uncertainty bands remain
narrow. For the most challenging in-distribution test samples, the
uncertainty bands widen as the prediction discrepancies increase,
particularly in the region where the wall heat flux changes rapidly.
The uncertainty band of Prob-DeepONet expands over a relatively broad
portion of the profile, whereas that of two-step MV-DeepONet varies more
coherently with the spatial variation of the wall heat flux field.

\FloatBarrier

\subsubsection{Analysis of Output Covariance Recovery Accuracy}
\label{sec:aero-cov}

This part further quantifies the output covariance recovery accuracy for the aerothermal example.

\vspace{4pt}
\noindent\textbf{Covariance Error Decomposition and Quantitative Accuracy Analysis}

\begin{figure}[!htbp]
	\centering
	\includegraphics[width=0.5\linewidth]
	{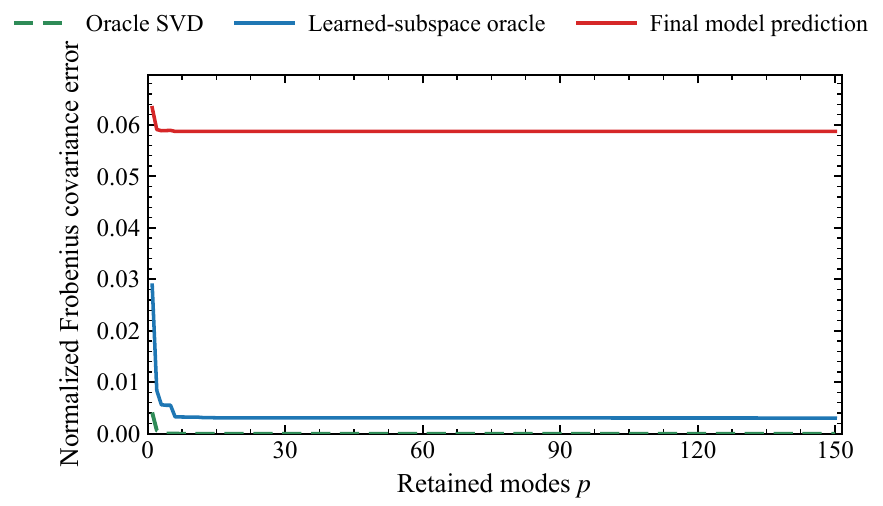}
	\caption{Normalized Frobenius covariance error as a function of the number of retained modes.}
	\label{fig:aero-cov-decomp}
\end{figure}

Figure~\ref{fig:aero-cov-decomp} shows how the normalized Frobenius covariance errors vary with the number of retained modes \(p\), while Table~\ref{tab:aero-cov-modes} reports the minimum numbers of modes required to reach different target normalized Frobenius errors. The main observations are as follows.

\textbf{(1) Oracle SVD curve (\(T_{1}\)).}
The Oracle SVD error falls below \(1\%\) with only one retained mode,
below \(10^{-4}\) for \(p\geq3\), and reaches \(1.55\times10^{-7}\) at \(p=10\). The reference covariance therefore has an approximately rank-one structure. Its physical origin and spectral characteristics are further discussed in Appendix~B. Consequently,
\(T_{1}\) is negligible compared with the model-dependent errors.

\begin{table}[!htbp]
	\centering
	\caption{Minimum numbers of retained modes required to achieve different target normalized Frobenius covariance-error levels for the three covariance-recovery curves.}
	\label{tab:aero-cov-modes}
	
	\small
	\renewcommand{\arraystretch}{1.2}
	\setlength{\tabcolsep}{7pt}
	
	\begin{tabular}{@{}cccc@{}}
		\toprule
		Target error
		& Oracle SVD
		& \makecell[c]{Learned-subspace\\oracle}
		& \makecell[c]{Final model\\prediction} \\
		\midrule
		
		Error \(\leq 20.00\%\)
		& \(p\geq1\)
		& \(p\geq1\)
		& \(p\geq1\) \\
		
		Error \(\leq 10.00\%\)
		& \(p\geq1\)
		& \(p\geq1\)
		& \(p\geq1\) \\
		
		Error \(\leq 5.00\%\)
		& \(p\geq1\)
		& \(p\geq1\)
		& Not reached \\
		
		Error \(\leq 2.00\%\)
		& \(p\geq1\)
		& \(p\geq2\)
		& Not reached \\
		
		Error \(\leq 1.00\%\)
		& \(p\geq1\)
		& \(p\geq2\)
		& Not reached \\
		
		\bottomrule
	\end{tabular}
\end{table}

\FloatBarrier

\textbf{(2) Comparison between the Learned-subspace oracle curve and the Final model prediction curve.}
The two curves remain clearly separated over the full range of \(p\).
For \(p\geq5\), the Learned-subspace oracle error stabilizes at
approximately \(3.0\times10^{-3}\), whereas the Final model prediction
error remains near \(6.6\times10^{-2}\). This behavior differs from
that observed in the preceding PDE examples and leads to two
observations:

\begin{enumerate}[leftmargin=2.2em,itemsep=0.4em]
	\item
	The final covariance recovery error is dominated by \(T_{3}\).
	The Learned-subspace oracle error is only approximately \(0.3\%\),
	indicating that the trunk network accurately captures the dominant
	output subspace. By contrast, the large separation from the Final
	model prediction curve identifies branch coefficient-covariance
	regression as the principal bottleneck. This result is consistent with the approximately rank-one covariance
	structure. Since most covariance energy is concentrated in the
	dominant mode, bias in its predicted coefficient variance is
	transferred directly to the Frobenius covariance error. The limited
	training set and the complex inputs combining latent species-field
	features with freestream parameters may further increase the
	difficulty of variance regression. Accordingly, the Final model
	prediction reaches the \(10\%\) error level but not the \(5\%\)
	level, whereas the other two curves fall below \(1\%\) for
	\(p\geq2\).
	
	\item
	For \(M=13{,}056\) and \(K=125\), the corresponding relative statistical-error scale is
	\begin{equation}
		\frac{
			\sqrt{M}K^{-1/2}
		}{
			\left\|
			\widehat{\boldsymbol{\Sigma}}
			\right\|_{F}
		}
		\approx
		6.2\times10^{-2}.
	\end{equation}
	This estimate is more than one order of magnitude larger than the observed \(T_{2}\approx3.0\times10^{-3}\), showing that the order bound is loose in this case. It therefore does not permit a reliable ranking of the two contributions within $T_{2}$.
\end{enumerate}

Overall, covariance recovery in the aerothermal example follows the hierarchy \(T_{1}\ll T_{2}\ll T_{3}\) in contrast to the \(T_{2}\)-dominated behavior observed in the preceding PDE examples. These results show that when covariance energy is concentrated in a single dominant mode, recovery accuracy becomes particularly sensitive to the predicted variance of the corresponding coefficient.

\vspace{4pt}
\noindent\textbf{Recovery of Local Correlation Structures:
	Multi-Anchor Correlation Map Assessment}

\begin{figure}[!htbp]
	\centering
	\includegraphics[width=0.70\linewidth]
	{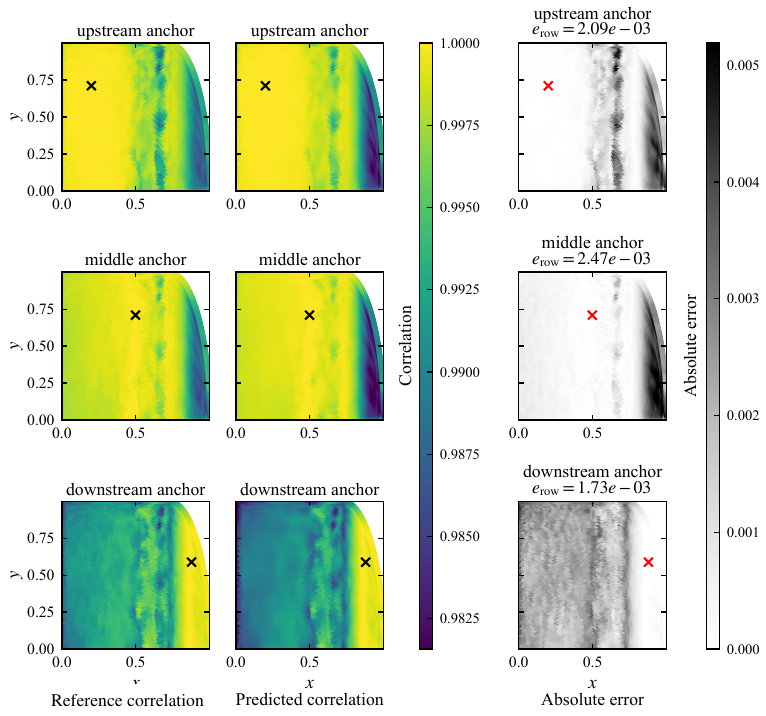}
	\caption{Multi-anchor correlation maps showing the reference correlations, predicted correlations, and corresponding absolute errors for three anchor points, where $e_{\mathrm{row}}$ denotes the relative $L_2$ error of the predicted correlation map for each anchor point.}
	\label{fig:aero-corr-anchor}
\end{figure}

For the aerothermal example, three anchor points are selected from the upstream, middle, and downstream regions of the wall along the flow direction. Fig.~\ref{fig:aero-corr-anchor} compares the corresponding reference, predicted, and absolute-error correlation maps.

\vspace{4pt}
\noindent\textbf{Overall pattern consistency.}
For all three anchors, the predicted and reference correlation maps show strong agreement in their overall spatial patterns. The spatial extent of the principal correlation region, the streamwise variation in correlation strength, and the transition between regions of relatively higher and lower correlation are all well reproduced. The absolute errors remain small over most of the wall, with a mean row-wise correlation error of approximately \(2.10\times10^{-3}\).
Although the normalized Frobenius covariance error is approximately \(6.6\%\) and dominated by \(T_{3}\), the local correlation error is only about \(0.2\%\). This difference indicates that the Frobenius error primarily reflects a bias in the magnitude of the dominant modal variance rather than an error in the recovered spatial correlation structure. The model therefore recovers the spatial correlation pattern accurately while slightly misestimating the overall covariance scale.

\vspace{4pt}
\noindent\textbf{Physical interpretation of the correlation structure.}
The aerothermal maps exhibit globally coherent positive correlations along the wall. This behavior is consistent with the physical mechanism of the aerothermal problem. Since the freestream Mach number is the only sampled random input, its perturbation acts as a common forcing on the coupled flow and thermochemical response. Over broad regions of the wall, the heat flux responds to variations in \(Ma\) in the same direction. Consequently, increases or decreases in \(Ma\) produce coordinated changes in wall heat flux, giving rise to the extensive positive correlations observed in the correlation maps.

This structure differs from the vertical, horizontal, and localized correlation patterns observed in the preceding examples. These differences reflect the distinct mechanisms through which uncertainty enters and propagates in the corresponding physical systems.

\FloatBarrier

\subsubsection{Prediction Interval Calibration and Coverage Analysis}
\label{subsec:aerothermal-cp-calibration}

CP is applied to calibrate the prediction intervals for the aerothermal example. Table~\ref{tab:aerothermal-cp-calibration} reports the corresponding PICP and MPIW values.

Before calibration, both methods produce conservative intervals with PICP values close to \(100\%\). At the \(90\%\) and \(95\%\) nominal levels, CP moves the coverage close to the prescribed targets and reduces the interval widths. At the \(95\%\) level, the MPIW of Prob-DeepONet is approximately \(1.74\) times that of two-step MV-DeepONet before calibration and \(1.65\) times after calibration. At the \(99\%\) level, the calibration effect is less uniform: Prob-DeepONet retains \(100\%\) coverage, while the MPIW of two-step MV-DeepONet increases despite achieving a PICP of \(99.14\%\). Nevertheless, its calibrated interval remains narrower than that of Prob-DeepONet. Overall, two-step MV-DeepONet provides comparable coverage with sharper intervals across the three nominal levels.

\begin{table}[!htbp]
	\centering
	\caption{PICP and MPIW of Prob-DeepONet and two-step MV-DeepONet before and after CP calibration.}
	\label{tab:aerothermal-cp-calibration}
	
	\small
	\renewcommand{\arraystretch}{1.18}
	\setlength{\tabcolsep}{5pt}
	
	\begin{tabular*}{\textwidth}{
			@{\extracolsep{\fill}}
			l
			c
			c
			c
			c
			c
			@{}
		}
		\toprule
		\multirow{2}{*}{Method}
		&
		\multirow{2}{*}{
			\makecell[c]{Nominal\\coverage (\%)}
		}
		&
		\multicolumn{2}{c}{Before CP calibration}
		&
		\multicolumn{2}{c}{After CP calibration}
		\\
		
		\cmidrule(lr){3-4}
		\cmidrule(l){5-6}
		
		&
		&
		PICP (\%)
		&
		MPIW ($\mathrm{W\,m^{-2}}$)
		&
		PICP (\%)
		&
		MPIW ($\mathrm{W\,m^{-2}}$)
		\\
		\midrule
		
		\multirow{3}{*}{Prob-DeepONet}
		&
		90
		&
		99.63
		&
		$4.0469\times10^{5}$
		&
		90.58
		&
		$2.0484\times10^{5}$
		\\
		
		&
		95
		&
		99.99
		&
		$4.8219\times10^{5}$
		&
		95.27
		&
		$2.9337\times10^{5}$
		\\
		
		&
		99
		&
		100.00
		&
		$6.3373\times10^{5}$
		&
		100.00
		&
		$4.8961\times10^{5}$
		\\
		
		\midrule
		
		\multirow{3}{*}{two-step MV-DeepONet}
		&
		90
		&
		99.05
		&
		$2.3197\times10^{5}$
		&
		90.63
		&
		$1.4897\times10^{5}$
		\\
		
		&
		95
		&
		99.99
		&
		$2.7641\times10^{5}$
		&
		94.79
		&
		$1.7783\times10^{5}$
		\\
		
		&
		99
		&
		100.00
		&
		$3.6326\times10^{5}$
		&
		99.14
		&
		$4.0642\times10^{5}$
		\\
		
		\bottomrule
	\end{tabular*}
\end{table}

\FloatBarrier

\section{Conclusions}
\label{sec:conclusions}

This paper developed a two-step MV-DeepONet framework for forward uncertainty propagation driven by random input fields, with particular emphasis on recovering the total predictive covariance and its structured conditional component. The framework was designed to relax the pointwise conditional-independence assumption of Prob-DeepONet by incorporating a two-step training strategy and transferring Gaussian probabilistic modeling from the physical output space to the low-dimensional modal coefficient space. The probabilistic modal coefficients jointly affect multiple output locations through the shared basis functions, thereby inducing a generally non-diagonal conditional predictive covariance without explicitly learning or storing the full high-dimensional covariance matrix. To further improve practical uncertainty estimation, deep ensembles and conformal prediction were incorporated to quantify model uncertainty and construct calibrated prediction intervals. More fundamentally, a Frobenius-norm error decomposition and corresponding upper bound were derived, identifying low-rank covariance compressibility, trunk-subspace approximation, finite-sample statistical error, and coefficient-space covariance estimation as the principal factors governing covariance recovery.

Numerical experiments on several PDE problems and a representative engineering problem validated the effectiveness of the proposed model. The results showed that although Prob-DeepONet fit the training data well, it exhibited larger generalization errors, wider prediction intervals, and more limited recovery of spatial correlation structures in challenging test and extrapolation scenarios, consistent with the limitations of its pointwise conditional covariance representation. Compared with the baseline, the two-step MV-DeepONet yielded more stable mean predictions, spatially structured uncertainty estimates, and improved covariance recovery. After calibration, both models approached the nominal coverage levels, while the proposed model achieved comparable coverage with narrower prediction intervals, indicating more compact and informative uncertainty estimates.

Future work may proceed in two directions. First, since the effectiveness of the proposed model depends partly on the low-rank compressibility of the output covariance and the retained modal dimension \(p\), adaptive mode-selection strategies may be investigated to balance subspace expressiveness and coefficient-space probabilistic regression complexity. Second, this paper mainly adopted Gaussian probabilistic modeling in the modal coefficient space. For cases where the conditional coefficient-space predictive distributions exhibit non-Gaussian features, such as skewness, heavy tails, or multimodality, future work may adopt non-Gaussian coefficient-space distributions to improve robustness in scenarios involving strong nonlinearity, high-dimensional random inputs, and complex spatially correlated outputs.

\section*{Declaration of Generative AI and AI-assisted technologies in the writing process}
During the preparation of this work the authors used ChatGPT in order to improve readability and language. After using this tool, 
the authors reviewed and edited the content as needed and take full responsibility for the content of the publication.

\section*{CRediT authorship contribution statement}
\textbf{Yupei Nie}: Idea of this paper, Simulations and calculations, Analysis, Code, Conceptualization, Validation, Software, Wrote this paper.  \textbf{Lei Wang}: Research direction, Idea of this paper, Algorithm, Total design scheme, Analysis, Conceptualization, Wrote this paper. \textbf{Jiasen Liu}: Numerical simulation of the aerothermal case. 

\section*{Declaration of competing interest}
The authors declare that they have no known competing financial interests or personal relationships that could have appeared to influence the work reported in this paper.

\section*{Data availability statements}
The data that support the findings of this study are available from the corresponding author upon reasonable request. 

\section*{Acknowledgments}
We express our sincere thanks to all the members of our discussion group for their valuable comments. 

\appendix

\renewcommand{\thesection}{\Alph{section}}
\renewcommand{\thesubsection}{\thesection.\arabic{subsection}}
\renewcommand{\theequation}{\thesection.\arabic{equation}}
\renewcommand{\thefigure}{\thesection.\arabic{figure}}
\renewcommand{\thetable}{\thesection.\arabic{table}}

\newcommand{\appendixsection}[1]{%
	\refstepcounter{section}%
	\setcounter{subsection}{0}%
	\setcounter{equation}{0}%
	\setcounter{figure}{0}%
	\setcounter{table}{0}%
	\section*{Appendix~\thesection.\ #1}%
	\addcontentsline{toc}{section}{Appendix~\thesection.\ #1}%
}

\appendixsection{VAE-Based Reduction for the Aerothermal Case Study in
	Section~\ref{sec:aerothermal}}
\label{app:vae}

The aerothermal flow field contains mass fraction distributions for
eleven chemical species. Directly using these fields
as operator-learning inputs would substantially increase the input
dimension, while the large differences in magnitude among species
would complicate network training. A three-dimensional VAE is therefore employed to compress the species fields
into a latent vector \(\boldsymbol{z}\), which provides a compact input
representation for the subsequent aerothermal prediction model. The
VAE is used solely for dimensionality reduction. Its data processing,
training objective, and reconstruction performance are described below.

\subsection{Data Processing and Network Architecture}
\label{app:vae-data-architecture}

The original flow fields are defined on an unstructured mesh. The 11 species mass fraction fields are therefore
first mapped onto a regular three-dimensional voxel grid using inverse
distance weighting. The grid has dimensions
$96\times26\times26$, and each sample contains 11 channels
corresponding to the mass fraction distributions of the chemical
species.

Because some voxels lie outside the valid flow domain, a binary mask
\(\mathbf{m}\) is introduced to identify the valid physical region, and
the reconstruction error is evaluated only over this region. 
To improve training stability in the presence of large differences in magnitude among species, the mass fractions are first transformed
using $\log_{10}$ and then standardized channelwise over the valid masked region:
\begin{equation}
	x_{\mathrm{std}}
	=
	\frac{
		\log_{10}
		\left(
		\max
		\left(
		x_{\mathrm{raw}},
		\varepsilon_{\log}
		\right)
		\right)
		-
		\mu_c
	}{
		s_c
	}.
	\label{eq:vae-standardization}
\end{equation}
Here, \(\mu_c\) and \(s_c\) denote the mean and standard deviation of
the \(c\)th channel over the valid region, and
\(\varepsilon_{\log}\) prevents the logarithm of zero.

\begin{figure}[!htbp]
	\centering
	\captionsetup{skip=2pt}
	
	\figplaceholder[
	width=0.85\linewidth,
	trim=0 18bp 0 50bp,
	clip
	]{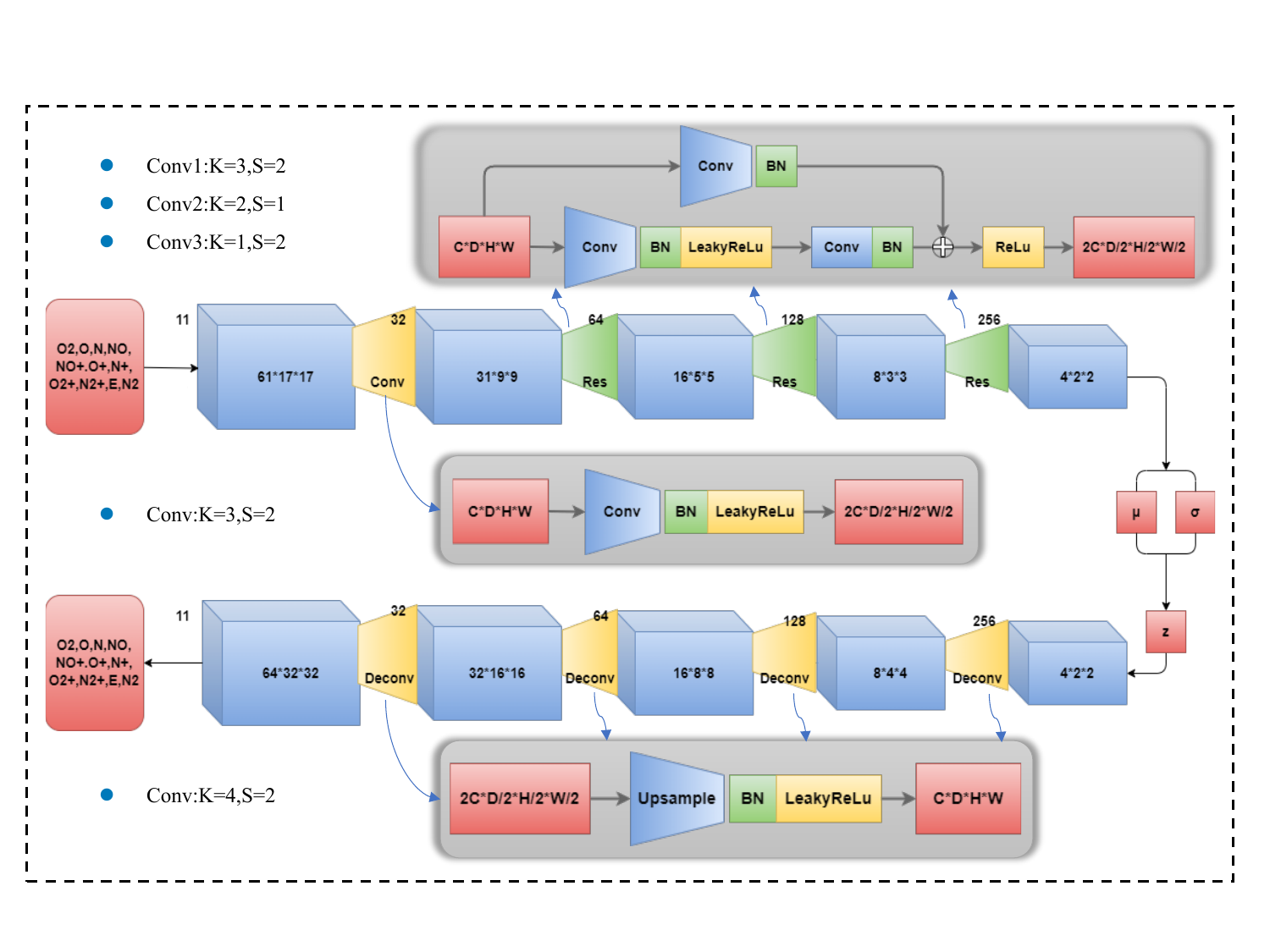}
	
	\caption{Schematic of the three-dimensional convolutional VAE architecture with residual blocks.}
	\label{fig:vae-architecture}
\end{figure}

The VAE consists of an encoder, a latent-variable sampling layer, and
a decoder. The encoder uses three-dimensional convolutions,
anisotropic downsampling, and residual blocks to map the standardized
multichannel voxel field to the latent distribution parameters
$\boldsymbol{\mu}$ and $\log\boldsymbol{\sigma}^{2}$. The latent
vector is then generated using the reparameterization trick,
\begin{equation}
	\boldsymbol{z}
	=
	\boldsymbol{\mu}
	+
	\boldsymbol{\sigma}
	\odot
	\boldsymbol{\epsilon},
	\qquad
	\boldsymbol{\epsilon}
	\sim
	\mathcal{N}
	\left(
	\boldsymbol 0,
	\boldsymbol I
	\right).
	\label{eq:vae-reparameterization}
\end{equation}
The decoder reconstructs the three-dimensional multichannel field
through upsampling, convolutional layers, and residual refinement
blocks. The complete architecture is illustrated in
Fig.~\ref{fig:vae-architecture}.


\subsection{Training Objective and Hyperparameter Selection}
\label{app:vae-training}

The VAE is trained using a combination of the reconstruction loss and the KL
regularization term:
\begin{equation}
	\mathcal{L}_{\mathrm{total}}
	=
	\mathcal{L}_{\mathrm{recon}}
	+
	\lambda_{\mathrm{KL}}\mathcal{L}_{\mathrm{KL}},
	\label{eq:vae-total-loss}
\end{equation}
where \(\lambda_{\mathrm{KL}}\) controls the strength of the Kullback–Leibler (KL) regularization.

The standard reconstruction loss is evaluated in the transformed
network input space and is defined as
\begin{equation}
	\mathcal{L}_{\mathrm{recon}}
	=
	\frac{1}{N_b C}
	\sum_{n=1}^{N_b}
	\sum_{c=1}^{C}
	\frac{
		\left\|
		\mathbf{m}\odot
		\left(
		\hat{\mathbf{x}}_{n}^{(c)}
		-
		\mathbf{x}_{n}^{(c)}
		\right)
		\right\|_{2}^{2}
	}{
		\left\|\mathbf{m}\right\|_{1}
	},
	\label{eq:vae-reconstruction-loss}
\end{equation}
where \(\mathbf{x}_{n}^{(c)}\) is the \(\log_{10}\)-transformed and
standardized input field, and \(\hat{\mathbf{x}}_{n}^{(c)}\) is
the corresponding decoder reconstruction in the same transformed
space.

The KL regularization term is given by
\begin{equation}
	\mathcal{L}_{\mathrm{KL}}
	=
	-\frac{1}{2N_b d_z}
	\sum_{n=1}^{N_b}
	\sum_{j=1}^{d_z}
	\left(
	1+\log\sigma_{n,j}^{2}
	-\mu_{n,j}^{2}
	-\sigma_{n,j}^{2}
	\right),
	\label{eq:vae-kl-loss}
\end{equation}
where \(\mu_{n,j}\) and \(\sigma_{n,j}^{2}\) are the mean and variance
of the \(j\)-th latent variable for the \(n\)-th sample, respectively.

\begin{figure}[htbp]
	\centering
	\figplaceholder[width=0.87\linewidth]
	{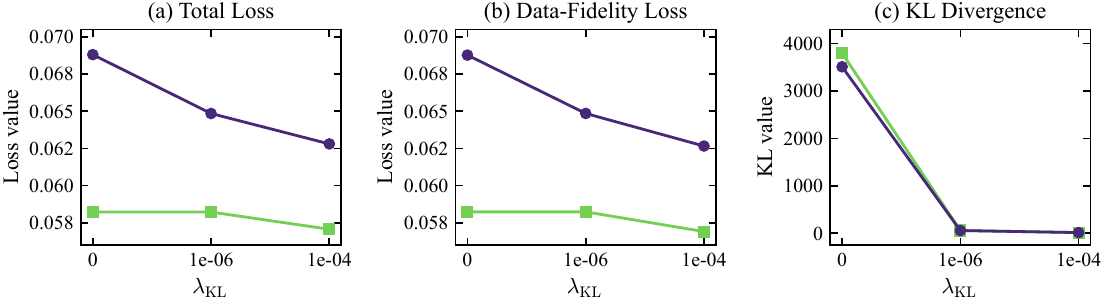}
	\caption{Hyperparameter sensitivity analysis for the VAE. The
		purple and green curves correspond to latent dimensions
		\(d_z=56\) and \(d_z=64\), respectively.}
	\label{fig:vae-hparam}
\end{figure}

A sensitivity study examines \(d_z\in\{56,64\}\) and \(\lambda_{\mathrm{KL}}\in\{0,10^{-6},10^{-4}\}\).
Fig.~\ref{fig:vae-hparam} compares the final loss components and
identifies \(d_z=64\) and \(\lambda_{\mathrm{KL}}=10^{-4}\) as the
best trade-off between reconstruction accuracy and latent-space
regularization.

\FloatBarrier

\subsection{Reconstruction of Species Mass Fraction Fields}
\label{app:vae-reconstruction}

\begin{figure}[!htbp]
	\centering
	
	\captionsetup[subfigure]{
		skip=2pt,
		justification=centering,
		singlelinecheck=true
	}
	
	
	\begin{subfigure}[t]{0.49\textwidth}
		\centering
		\includegraphics[
		width=\linewidth,
		keepaspectratio
		]{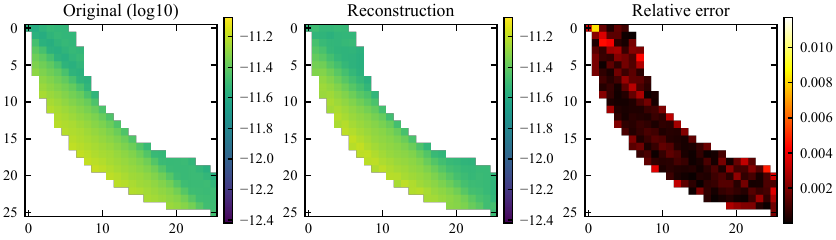}
		\caption{Electron species \(e^{-}\).}
		\label{fig:vae-reconstruction-e}
	\end{subfigure}
	\hfill
	\begin{subfigure}[t]{0.49\textwidth}
		\centering
		\includegraphics[
		width=\linewidth,
		keepaspectratio
		]{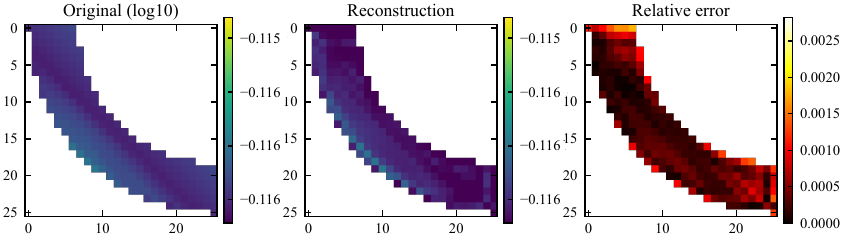}
		\caption{Molecular nitrogen species \(N_{2}\).}
		\label{fig:vae-reconstruction-n2}
	\end{subfigure}
	
	\par\vspace{0.5em}
	
	
	\begin{subfigure}[t]{0.49\textwidth}
		\centering
		\includegraphics[
		width=\linewidth,
		keepaspectratio
		]{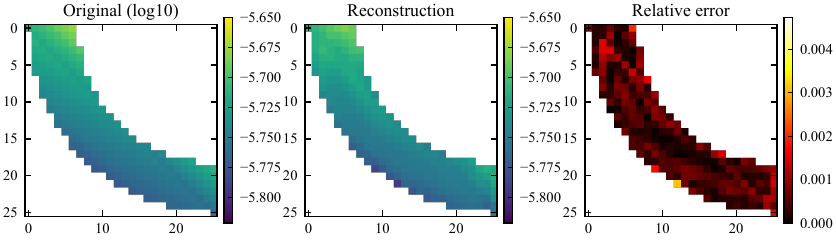}
		\caption{Nitric oxide ion species \(NO^{+}\).}
		\label{fig:vae-reconstruction-no-plus}
	\end{subfigure}
	\hfill
	\begin{subfigure}[t]{0.49\textwidth}
		\centering
		\includegraphics[
		width=\linewidth,
		keepaspectratio
		]{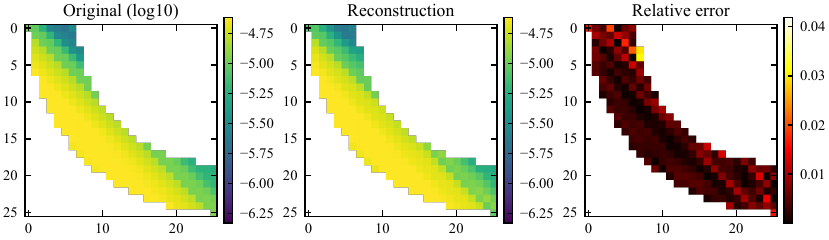}
		\caption{Atomic oxygen species \(O\).}
		\label{fig:vae-reconstruction-o}
	\end{subfigure}
	
	\caption{Representative slice reconstruction results for the species mass fraction fields. Within each subfigure, the original field, reconstructed field, and corresponding relative-error field are shown from left to right.}
	\label{fig:vae-reconstruction-representative}
\end{figure}

\FloatBarrier

The reconstruction performance is illustrated using four representative
species: \(e^{-}\), \(N_{2}\), \(NO^{+}\), and \(O\). These species
cover electronic, molecular, ionic, and atomic components of the
multispecies flow field. Fig.~\ref{fig:vae-reconstruction-representative}
compares the original fields, reconstructed fields, and relative-error
fields. Results for the remaining species exhibit similar
reconstruction behavior and are omitted for conciseness.

Overall, the principal concentration and gradient patterns are well
preserved, with errors mainly confined to regions of sharp variation. 
This confirms that the \(d_z=64\) latent representation preserves the
essential information in the multispecies flow field.
The VAE therefore provides a compact representation of the species mass
fraction fields while maintaining satisfactory reconstruction accuracy,
thereby reducing the input dimension and computational cost of the
subsequent operator-learning model.

\FloatBarrier



\appendixsection{Approximate Rank-One Covariance Structure for the Aerothermal Case Study in Section~\ref{sec:aerothermal}}
\label{app:rank-one-covariance}

The aerothermal output covariance exhibits an approximately rank-one structure, which helps explain the distinct error hierarchy observed in the main text. This appendix provides a first-order explanation for this structure based on the one-dimensional effective stochastic input.

\subsection{Effective Stochastic Dimension of the Input}

In the aerothermal example, the freestream temperature, pressure, and density are fixed, and only the freestream Mach number varies randomly over \(Ma\in[10.4,16.4]\). For a prescribed \(Ma\), the species mass fraction fields are deterministic CFD outputs, and the VAE latent vector \(\boldsymbol{z}\) used by Branch~2 is therefore uniquely determined by \(Ma\). Consequently, although the multi-branch model receives several input components, all samples lie on a one-dimensional manifold parameterized by \(Ma\). The effective stochastic dimension is therefore one, unlike the random-field inputs used in the preceding PDE examples. This one-dimensional stochastic structure provides the basis for the approximately rank-one covariance obtained from the first-order Taylor expansion below.

\subsection{First-Order Derivation of the Approximately Rank-One Covariance}

Let the physical operator \(\mathcal{G}\) map the input function
$\mathbf{u}$ to the wall heat flux field
$\mathbf{s}=\mathcal{G}(\mathbf{u})\in\mathbb{R}^{M}$. For a perturbation
$\delta\mathbf{u}$ about the mean input $\bar{\mathbf{u}}$, a first-order
Taylor expansion gives
\begin{equation}
	\mathbf{s}
	=
	\mathcal{G}\left(
	\bar{\mathbf{u}}+\delta\mathbf{u}
	\right)
	\approx
	\mathcal{G}\left(
	\bar{\mathbf{u}}
	\right)
	+
	\nabla \mathcal{G}\left(
	\bar{\mathbf{u}}
	\right)
	\delta\mathbf{u},
	\label{eq:appB-first-order-output}
\end{equation}
where $\nabla \mathcal{G}(\bar{\mathbf{u}})$ denotes the linearized sensitivity
of the physical operator at $\bar{\mathbf{u}}$. Let
$\boldsymbol{\mu}_{s}=\mathbb{E}[\mathbf{s}]$. The corresponding output
covariance is approximated by
\begin{equation}
	\boldsymbol{\Sigma}
	=
	\mathbb{E}\!\left[
	(\mathbf{s}-\boldsymbol{\mu}_{s})
	(\mathbf{s}-\boldsymbol{\mu}_{s})^{\top}
	\right]
	\approx
	\nabla \mathcal{G}(\bar{\mathbf{u}})
	\operatorname{Cov}(\delta\mathbf{u})
	\nabla \mathcal{G}(\bar{\mathbf{u}})^{\top}.
	\label{eq:app-cov-prop}
\end{equation}

Under these conditions, the composite input function is completely
determined by the single scalar \(Ma\), namely
\(\mathbf{u}=\mathbf{u}(Ma)\). Let
\[
Ma=\overline{Ma}+\delta Ma,
\]
where \(\overline{Ma}=\mathbb{E}[Ma]\) denotes the mean Mach number and
\(\delta Ma\) is the corresponding zero-mean fluctuation, satisfying
\(\mathbb{E}[\delta Ma]=0\). Expanding this vector-valued function about
\(\overline{Ma}\) gives
\begin{equation}
	\delta\mathbf{u}
	=
	\mathbf{u}(\overline{Ma}+\delta Ma)
	-
	\mathbf{u}(\overline{Ma})
	\approx
	\left.
	\frac{\partial\mathbf{u}}{\partial Ma}
	\right|_{\overline{Ma}}
	\delta Ma
	=:
	\mathbf v\,\delta Ma
	\label{eq:appB-input-perturbation}
\end{equation}
where \(\mathbf v=\left.\dfrac{\partial\mathbf{u}}{\partial Ma}\right|_{\overline{Ma}}\)
is a deterministic vector in the input space and \(\delta Ma\) is a
scalar random variable. The first-order approximation to the input
covariance is therefore
\begin{equation}
	\operatorname{Cov}(\delta\mathbf{u})
	\approx
	\mathbb{E}\!\left[
	(\mathbf v\,\delta Ma)
	(\mathbf v\,\delta Ma)^{\top}
	\right]
	=
	\operatorname{Var}(Ma)\,
	\mathbf v\mathbf v^{\top}.
	\label{eq:appB-input-covariance}
\end{equation}

Substituting this result into Eq.~\eqref{eq:app-cov-prop} and defining
\(\mathbf w=\nabla \mathcal{G}(\bar{\mathbf{u}})\mathbf v\), we obtain
\begin{equation}
	\boldsymbol{\Sigma}
	\approx
	\operatorname{Var}(Ma)
	\left[
	\nabla \mathcal{G}(\bar{\mathbf{u}})\mathbf v
	\right]
	\left[
	\nabla \mathcal{G}(\bar{\mathbf{u}})\mathbf v
	\right]^{\top}
	=
	\operatorname{Var}(Ma)\,
	\mathbf w\mathbf w^{\top}.
	\label{eq:app-rank1}
\end{equation}

The first-order approximation to the output covariance is therefore a
scalar multiple of an outer product. Its rank is at most one and
equals one whenever $\mathbf w\neq\mathbf{0}$. Hence, the
one-dimensional input uncertainty directly induces an approximately
rank-one output covariance after first-order propagation through the
physical operator. This result is consistent with the observation in the main text that
the Oracle SVD curve ($T_{1}$) approaches zero with only a few retained
modes.

The vector
\begin{equation}
	\mathbf w
	=
	\nabla \mathcal{G}
	\left(
	\bar{\mathbf{u}}
	\right)
	\mathbf v
	=
	\left.
	\frac{\partial\mathbf{s}}{\partial Ma}
	\right|_{\overline{Ma}}
	\label{eq:appB-mach-sensitivity}
\end{equation}
has a direct physical interpretation. Its $i$th component represents
the local sensitivity of the wall heat flux at the $i$th wall grid
point to a change in the freestream Mach number. Larger values of
$\lvert w_i\rvert$ indicate regions with greater sensitivity to $Ma$,
such as the region near the stagnation point.

To first order, the perturbation of the wall heat flux field satisfies
\begin{equation}
	\delta\mathbf{s}
	\approx
	\mathbf w\,\delta Ma.
	\label{eq:appB-output-perturbation}
\end{equation}
Hence, all wall locations are driven by the same scalar fluctuation
and differ only in their local sensitivities, producing coherent
variations across the surface.

\subsection{Interpretation of the Evaluation Metrics}

Equation~\eqref{eq:app-rank1} directly gives the following first-order
approximation to the correlation coefficient between any two wall
locations $i$ and $j$:
\begin{equation}
	\rho_{ij}
	=
	\frac{
		\operatorname{Cov}(s_i,s_j)
	}{
		\sqrt{
			\operatorname{Var}(s_i)
			\operatorname{Var}(s_j)
		}
	}
	\approx
	\frac{
		w_iw_j
	}{
		\lvert w_i\rvert\lvert w_j\rvert
	},
	\qquad
	w_iw_j\neq0.
	\label{eq:appB-correlation-coefficient}
\end{equation}
Under these conditions, the wall heat flux
exhibits a generally positive sensitivity to the freestream Mach
number over the entire surface. Therefore, \(w_i>0\) and \(\rho_{ij}\approx1\) are expected. 
This agrees with Fig.~\ref{fig:aero-corr-anchor}, where the reference correlations are
predominantly within \(0.98\)--\(1.00\).

The approximately rank-one structure also explains the sensitivity of
the normalized Frobenius error. For
\(\boldsymbol{\Sigma}\approx
\lambda_1\mathbf{u}_1\mathbf{u}_1^\top\), nearly all covariance energy is
concentrated in one dominant mode. Any relative error in this direction therefore contributes directly to the normalized Frobenius error and cannot be diluted by energy in other directions. By contrast, for a higher-rank compressible
covariance, the energy is distributed across multiple directions, so
an error in a single direction is naturally averaged and has a smaller
relative effect.
Therefore, for an approximately rank-one reference covariance, spatial correlation recovery provides a complementary metric that assesses whether the correlation structure is accurately recovered independently of the overall covariance magnitude.

Finally, this physical structure also highlights a limitation of the
pointwise conditional covariance representation adopted in Prob-DeepONet.
For an individual input \(\mathbf{u}\), its conditional predictive
covariance is
\begin{equation}
	\widehat{\boldsymbol{\Sigma}}_{\mathrm{Prob}}(\mathbf{u})
	=
	\operatorname{diag}
	\left(
	\hat{\sigma}_{s,1}^{2}(\mathbf{u}),
	\ldots,
	\hat{\sigma}_{s,M}^{2}(\mathbf{u})
	\right).
	\label{eq:appB-prob-covariance-second}
\end{equation}
Thus, the conditional predictive component contains no off-diagonal
dependence. The total predictive covariance additionally includes the
covariance of the conditional means across random input realizations and
therefore need not be diagonal. Nevertheless, when the reference response
exhibits strong coherent cross-location correlations, restricting the
conditional predictive covariance to a diagonal form limits the explicit
representation of structured dependence in this component of the total
predictive covariance.

By contrast, two-step MV-DeepONet models uncertainty in the modal
coefficient space and maps the coefficient-space covariance to the
physical output space through the shared basis functions, thereby allowing
a generally non-diagonal conditional predictive covariance. This
representation is better suited to capturing the strongly correlated
response structure observed in the present aerothermal example.

\bibliographystyle{elsarticle-num}
\bibliography{references}

\end{document}